\documentclass[12pt]{extarticle}
\usepackage[top=1cm,bottom=1.7cm,left=1.8cm,right=1.8cm]{geometry}
\usepackage{amssymb,amsthm,amsmath,amsxtra,dsfont,appendix}
\usepackage[nottoc]{tocbibind}
\usepackage{graphicx,tikz,esint}
\usepackage{hyperref}

\numberwithin{equation}{section} 

\renewcommand{\limsup}[1]{\underset{{#1}}{\overline{\operatorname{lim}}}}

\newcommand{\supp}{\operatorname{supp}}
\newcommand{\tr}{\operatorname{Tr}}
\newcommand{\Var}{\operatorname{Var}}
\renewcommand{\P}[1]{\mathbb{P}\left[#1\right]}
\newcommand{\E}[2]{\mathbb{E}_{#1}\left[#2\right]}

\newcommand{\1}{\mathds{1}}
\newcommand{\A}{\mathfrak{A}}
\newcommand{\B}{\operatorname{B}}
\renewcommand{\c}{\mathrm{c}}
\newcommand{\C}{\mathbb{C}}
\newcommand{\Cu}{\operatorname{C}}
\newcommand{\CD}{\operatorname{CD}}
\newcommand{\F}{\mathfrak{F}}
\newcommand{\G}{\operatorname{G}}
\newcommand{\h}{\operatorname{H}}
\newcommand{\I}{\operatorname{I}}
\newcommand{\K}{\mathrm{K}} 
\newcommand{\M}{\operatorname{M}}
\newcommand{\N}{\mathbb{N}}
\newcommand{\No}{\mathcal{N}}
\newcommand{\R}{\mathbb{R}}
\newcommand{\s}{\underline{s}}
 
\newcommand{\Sy}{\mathbb{S}}
\newcommand{\T}{\mathbb{T}} 
\newcommand{\U}{\mathcal{U}}  
\newcommand{\X}{\mathfrak{X}} 
\newcommand{\Z}{\mathbb{Z}}

\newtheorem{theorem}{Theorem}[section]
\newtheorem{definition}[theorem]{Definition}
\newtheorem{proposition}[theorem]{Proposition}
\newtheorem{corollary}[theorem]{Corollary}

\newtheorem{lemma}[theorem]{Lemma}
\newtheorem{remark}[theorem]{Remark}

\title{ \vspace{2cm} 
Gaussian and non-Gaussian fluctuations for mesoscopic linear statistics in determinantal processes}
\date{}
\author{Kurt Johansson\thanks{KTH Royal Institute of Technology, Department of Mathematics. Supported by the grant KAW 2010.0063 from the Knut and Alice Wallenberg Foundation.}
\and
Gaultier Lambert\thanks{KTH Royal Institute of Technology, Department of Mathematics, glambert@kth.se. Supported by the grant KAW 2010.0063 from the Knut and Alice Wallenberg Foundation.}}

\begin{document}

\maketitle

\vspace{1.5cm}

\begin{abstract} \normalsize
We study mesoscopic linear statistics for a class of determinantal point processes which interpolate between Poisson and Gaussian Unitary Ensemble (GUE) statistics. These processes are obtained by modifying the spectrum of the correlation kernel of the GUE eigenvalue process. An example of such a system comes from considering the distribution of non-colliding Brownian motions in a cylindrical geometry, or a grand canonical ensemble of free fermions in a quadratic well at positive temperature. When the scale of the modification of the spectrum of the GUE kernel, related to the size of the cylinder or the temperature, is different from the scale in the mesoscopic linear statistic, we get a central limit theorem (CLT) of either Poisson or GUE type. 
On the other hand, in the critical regime where the scales are the same, we get a non-Gaussian process in the limit.
Its distribution is characterized by explicit but complicated formulae for the cumulants of  smooth linear statistics. These results rely on an asymptotic sine-kernel approximation of the GUE kernel which is valid at all mesoscopic scales, and a generalization of cumulant computations of Soshnikov for the sine process. Analogous determinantal processes on the circle are also considered with similar results.
\end{abstract}

\vspace{1.5cm}



\clearpage

\
\vspace{1cm}

\tableofcontents

\clearpage

\section{Introduction and results}  \label{sect:intro_0}

\subsection{Introduction} \label{sect:intro_1}

Recently there has been a lot of discussion about universality of  random matrices statistics at  mesoscopic or intermediate scales. 
For instance, the proofs of the local semicircle law and the Wigner-Dyson-Gaudin-Metha conjecture, see~\cite{EY_12,BEYY_14,Erdos_14} for further references, the work~\cite{EK_14a, EK_14b} on random band matrices and the so-called Anderson metal-insulator transition, or  the CLT for linear statistics of orthogonal polynomial ensembles~\cite{BD_15}. One motivation to investigate these models comes from E$.$ Wigner's fundamental observation that the spectral statistics of complicated quantum systems exhibit universal patterns. On the other hand, eigenvalues of quantum systems which are classically integrable are expected to be described by Poisson statistic~\cite{Gutzwiller}. Therefore, it is natural to investigate the transition from Poisson to random matrix statistics at intermediate scales.
There are many ways to interpolate between the two point processes, for instance using Dyson's Brownian motion, one gets a determinantal process called the {\it deformed} Gaussian Unitary Ensemble~\cite{Johansson_07}. For this model, the transition has been investigated using mesoscopic linear statistics in~\cite{DJ_14}. The authors proved central limit theorems whose fluctuations depend on the scale the test function samples the eigenvalues.
In this paper we will study the fluctuations of another general class of determinantal processes which interpolate between Poisson and GUE statistics that we call {\it modified} GUEs (see definition \ref{modified_GUE}). Instead of adding an independent matrix to a GUE matrix, we directly {\it modify} or mollify the spectrum of the correlation kernel of the process. This has the effect of introducing some extra disorder in the system while keeping the determinantal structure.
Our main motivation to study such ensembles comes from the so-called  MNS ensemble which was introduced by Moshe, Neuberger and Shapiro in \cite{MNS_94}, see also \cite{Johansson_07,DDMS_14},  and corresponds to the joint probability density function of the positions of a Grand-Canonical system of Free fermions at positive temperature confined in a one dimensional harmonic potential.
In general, it  is also of interest to investigate fluctuations  of determinantal processes whose correlation kernels are not necessarily reproducing.

\subsection{The Modified Gaussian Unitary Ensembles} \label{sect:intro_2}


Let $\X$ be a Polish space equipped with some reference measure $d\mu$. In the sequel, we will only be interested  in the two cases $\X=\R$ or the unit circle $\T$ equipped with the Lebesgue measure.
A point process is a random measure on $\X$ of the form $\Xi = \sum \delta_{X_i}$. The support of the measure $\Xi$ is the random object of interest, it is called a  point configuration $\{X_i\}$ and we suppose that it has neither double points nor accumulation points. 
Point processes are usually described by  their correlation functions  $\rho_k(x_1,\dots,x_k)$ which are characterized by
\begin{equation} \label{point_process}
 \E{}{\prod_i \big(1+g(X_i)\big)}= \sum_{k=1}^\infty \frac{1}{k!} \int_{\X^k} \prod_i g(x_i) \rho_k(x_1,\dots,x_k) d\mu(x_1)\cdots d\mu(x_k)  \end{equation}
for any measurable bounded function $g:\X\to\C$ with compact support.
%
%
 A point process is called determinantal if its correlation functions exist and satisfy the identity
\begin{equation*} \rho_k(x_1,\dots,x_k)= \underset{k\times k}{\det} [ K(x_i,x_j) ]  .\end{equation*}
Hence a determinantal process is characterized by its correlation kernel  $K:\X\times\X\to\R$ and we will denote by $\mathbb{E}_K$ the corresponding probability measure on the space of point configurations. One generally assumes that $K$ defines an integral operator $\K$ on $L^2(\X,d\mu)$ which is locally of trace class and then the RHS of equation (\ref{point_process}) is a Fredholm determinant:
\begin{equation}\label{Fredholm_det}
\E{K}{\prod_i (1+g(x_i))}= \det[\I+ \K g]_{L^2(d\mu)}  .
\end{equation} 
In most cases, the operator $\K$ is self-adjoint (although there are natural examples of non-Hermitian determinantal processes, such as the deformed GUE studied in~\cite{DJ_14}) and its kernel $K$ defines a determinantal process if and only if all the eigenvalues  of the operator $\K$ lie in $[0,1]$. 
These facts are well-known and we refer to~\cite{ HKPV_06,Johansson_05,Soshnikov_00b} for different introductions to the theory of (determinantal) point processes and to the survey~\cite{Borodin_11} for an overview of some applications.
%
%
In this paper we will investigate examples of determinantal processes with correlation kernels of the general form
\begin{equation} \label{kernel_0} K(x,y)= \sum_{k=0}^\infty p^N_k \vartheta_k(x)\vartheta_k(y)  ,\end{equation} 
where $(\vartheta_k )_{k=0}^\infty$ is an orthonormal basis in $L^2(\X, d\mu)$ and the  {\it spectrum} $0 \le p_k^N \le 1$.
A classical example is the correlation kernel of the GUE eigenvalue process:
 \begin{equation} \label{GUE_kernel}
K_0^N(x,y)= \sum_{k=0}^{N-1} \varphi_k(x)\varphi_k(y)  ,
 \end{equation}
 where  $ \varphi_k(x)=\sqrt{\frac{\pi\sqrt{N}}{\sqrt{2}}}h_k\left( x\frac{\pi\sqrt{N}}{\sqrt{2}} \right) e^{- \left( x\frac{\pi\sqrt{N}}{\sqrt{2}} \right)^2/2} $
and $h_k$ is the normalized Hermite polynomial with respect to the weight $e^{-x^2}$ on $\R$. 
The parameter~$N \in\mathbb{N}$ is the dimension of the matrix and the variance of the entries is scaled so that the eigenvalue density at the origin is of order $N$ as $N\to\infty$. 
For some background on the GUE process, we refer to~\cite[Chap.~2-5]{Pastur_Shcherbina} and to the appendix~\ref{sect:Asymptotic} for a collection of standard facts on the asymptotics of the Hermite polynomials and the GUE kernel. In the following, our main interest is in determinantal processes whose correlation kernels are modifications of $K_0^N$, in the sense that, instead of taking the {\it spectrum} $p_k^N= \1_{k<N}$, we assume that $k\mapsto p_k^N$ is a  function which decays from 1 to 0.
The following class of functions will be called {\bf shapes}:
\begin{equation} \label{class_F}
 \F=\{ \Psi: \R\to[0,1]\  |\ \Psi' \le 0  \text{ is Riemann integrable},
 \Psi \in L^1(0,\infty) \text{ and}\  (1-\Psi)\in L^1(-\infty,0) \}  \end{equation}
and, for any $\Psi\in\F$, we let
\begin{equation} \label{B^2}
\B^2_\Psi = \int_\R  \Psi(x)\big(1-\Psi(x)\big) dx .
\end{equation}


\begin{definition} \label{modified_GUE} 
A {\bf modified GUE} is a determinantal process on $\R$ $($with respect to the Lebesgue measure$)$  whose correlation kernel is 
\begin{equation}  \label{kernel}
K^N_{\Psi,\alpha}(x,y)= \sum_{k=0}^\infty \Psi\left(\frac{k-N}{\tau N^\alpha}\right)  \varphi_k(x) \varphi_k(y) .
\end{equation}
where  $\Psi\in\F$,  $N>0$, $\alpha\in(0,1)$ and $\tau>0$.
In the sequel,  the parameter $\alpha$ is  called the {\bf modification scale},  $\tau$ is called the {\bf temperature} and we assume that the shape $\Psi$ is normalized so that $\B^2_\Psi =1$.  \end{definition}

Note that he  kernel $K^N_{\Psi,\alpha}$ is not reproducing, so that the total number of particles in the process, denoted by~$\#$, is random. 
Moreover,  the point process is scaled so that its density at the origin is of order  $N$ (global scaling) and a  simple computation yields $\E{K^N_{\Psi,\alpha}}{\#} \sim N$. According to~\cite[Thm.~7]{HKPV_06}, another correlation kernel for the modified GUE  is given by
\begin{equation} \label{Bernoulli}
K(x,y)= \sum_{k=0}^\infty I^N_k \varphi_k(x)\varphi_k(y)   , 
\end{equation}   
where $I^N_k$ are independent Bernoulli random variables with $\E{}{I^N_k}=\Psi\left(\frac{k-N}{\tau N^\alpha}\right)$. Hence, the modified ensembles are {\it more random} than the GUE and  the amount of {\it extra randomness} is estimated by
\begin{align*}
 \Var_{K^N_{\Psi,\alpha}}[\#] 
  &= \sum_{k=0}^\infty \Var[I_k^N] 
 = \sum_{k=0}^\infty \Psi\left(\frac{k-N}{\tau N^\alpha}\right) \left( 1- \Psi\left(\frac{k-N}{\tau N^\alpha}\right)\right)  \\
  &\sim \tau N^{\alpha} \B^2_\Psi  \ ,   \end{align*}
by a Riemann sum approximation. 
Heuristically, it means that the more the spectrum of the correlation kernel is modified, the more disorder is forced into the system. So we expect that, for large modifications, the modified ensembles behave like the Poisson process rather than like the GUE.

\begin{remark}\label{rk:shape}
According to formula \eqref{kernel},  for any $\sigma>0$, the shapes $\Psi$ and $\Psi_\sigma(t)=\Psi(t/\sigma)$ define the same modified GUE at different temperatures.
Moreover, by formula $(\ref{B^2})$, 
$ \B^2_{\Psi_\sigma} =  \sigma \B^2_\Psi$ and the  condition $\B^2_\Psi=1$ fixes the temperature $\tau>0$ so that
\begin{equation} \label{variance_modification}  
\Var_{K^N_{\Psi,\alpha}}[\#] = \tau N^{\alpha} + \underset{N\to\infty}{o}(1)  \ . \end{equation}
\end{remark}


Our interest in determinantal processes with correlation kernels of the form (\ref{kernel}) is mainly motivated by the following example that we call the {\bf MNS ensemble}. In~\cite{MNS_94}, motivated by the physics of disordered systems,  Moshe, Neuberger and Shapiro  introduced an ensemble of unitary invariant Hermitian matrices whose eigenvalue distribution interpolates between the GUE and the Poisson process. This model was rigorously analyzed in~\cite{Johansson_07} and it was proved that its Grand Canonical version is a determinantal process with correlation kernel on $\R$ given by
\begin{equation} \label{kernel_1}
 K^N_{\psi,\alpha}(x,y)= \sum_{k=0}^\infty \frac{1}{1+e^{(k-N)/\tau N^\alpha}}\varphi_k(x)\varphi_k(y) \ ,\end{equation} 
for some $\tau>0$ and $0<\alpha<1$. So the MNS ensemble is a modified GUE with shape $\psi(t)=(1+e^{t})^{-1}$.  Moreover, this model has two natural interpretations.
%
%
First, since the rescaled Hermite functions $(\varphi_k)_{k\ge 0}$ are the eigenfunctions of the Schr\H{o}dinger operator $-\Delta + \frac{\pi^2 N}{2} x^2 $ on $\R$, the MNS process describes a grand canonical system of free fermions at positive temperature confined in a quadratic external well. Note that the probability that the k$^{\text{th}}$ state of this harmonic oscillator is occupied is equal to the Fermi factor $(1+e^{(k-N)/T})^{-1}$ where $T=N^\alpha\tau$ is the temperature of the system. Thus, the  Gaussian Unitary Ensemble corresponds to the ground state of such a system with $N$ fermions. Namely, taking the temperature to zero (i.e. the limit $\tau\to 0$ in (\ref{kernel_1})), one recovers the GUE kernel $K_0^N$ given by (\ref{GUE_kernel}).
On the other hand,  for large temperature (i.e. taking $\tau \to \infty $), the kernel degenerates to that of a Poisson process on $\R$. Therefore, at a heuristic level, the MNS ensemble interpolates between Poisson and random matrix (GUE) statistics. We shall prove that such a {\it transition} occurs for smooth mesoscopic linear statistics of the process.
%
%
The kernel \eqref{kernel_1} also occurs in connection with the KPZ equation (Kardar-Parisi-Zhang), where it is related to the crossover distribution for the height function, see \cite{SS10a, SS10b, ACQ11, DDMS_14}. We are not aware of any connection between the present work and the KPZ equation. 
%
%
Second, in ~\cite{Johansson_07}, it was shown that the MNS process also describes a system of Brownian particles moving on a cylinder and conditioned not to collide
(by rotation invariance, the distribution of the particles is stationary). We have seen that the parameter $N$ is the expected number of particles and one can check that  the length $\beta$ of the cylinder where the particles are diffusing is related to the temperature $T$ of the Free Fermions by $\beta = 4 \sinh(T^{-1})$. This particle system is expected to behave like Dyson's Brownian motion, \cite{AGZ,DJ_14}, and this provides another heuristic description of the transition. Namely, at small scales, the particles remain roughly independent, while when $\beta$ gets large the trajectories start regularizing because of the non-colliding constraints until eventually their joint distribution obeys the law of the GUE eigenvalues.  \\

In general, one can still think of  $\Psi\left(\frac{k-N}{\tau N^\alpha}\right)$ as the probability that the k$^{\text{th}}$ state of a quantum system is occupied. Then, $1-\Psi$ corresponds to the distribution function of a probability measure on $\R$ and we denote by  $\Phi=-\Psi'$ the corresponding density. For technical reasons, it will be simpler to consider the following subclass of shapes,
\begin{equation} \label{class_F'}
\F^*=\{ \Psi \in \F\ |\  \Phi(x)=   - \Psi'(x) \le e^{-c|x|} \text{ for some } c>0  \} ,  \end{equation}
which contains the MNS shape $\psi(t)=(1+e^{t})^{-1}$. 

\clearpage


In random matrix theory, it is well-known that one can analyze the eigenvalue processes at different scales. The global or {\it macroscopic} scale refers to the size of the whole process. 
On the other hand, the local or {\it microscopic} scale is that of individual eigenvalues, i.e.~the gaps between consecutive eigenvalues are of order 1. At this ultimate scale, in the Hermitian case, universality means that the rescaled point process converges in the bulk to the celebrated {\it sine process}. Any scale in between is called {\it mesoscopic}. 
In other words, a {\it mesoscopic} random variable is a function of the point process which depends on a growing fraction of the total number of particles. A typical example of such observables is the following class of linear statistics.


\begin{definition} \label{linear_statistic}
Given a point process $\Xi$ with density of order $N$ at the origin and a function $f: \X \to \C$ with compact support, for any  $0 < \delta<1$, we define 
$$ f_\delta(x) = f(x N^{\delta}) \ , $$ 
and we call a {\bf mesoscopic linear statistic} the random variable
 \begin{equation} \label{LS_1}  \Xi f_\delta =\sum_i f(X_i N^\delta) \ .
 \end{equation}
 In the following, the parameter $\delta$ is called the {\bf scale}.
\end{definition}
We will investigate the asymptotic distribution of mesoscopic linear statistics of the modified GUEs. Note that since the density at the origin is of order $N$, we have
\begin{equation} \label{LS_mean}
 \E{K_{\Psi,\alpha}^N}{\Xi f_\delta}  \sim N^{1-\delta} \int_\R f(x)  dx \ .
 \end{equation}
 If $\delta<1$, this expectation is diverging as $N\to\infty$ and it is natural to consider centered linear statistics instead:
 \begin{equation} \label{LS_2}  \tilde\Xi f_\delta =\sum_i f(X_i N^\delta) -   \E{K_{\Psi,\alpha}^N}{\Xi f_\delta} \ .
 \end{equation}   

For any   random variable $Z$ with a well-defined Laplace transform, its cumulants $C^n[Z]$ are given by the power series
\begin{equation} \label{CGF}
\log\E{}{e^{t Z}}= \sum_{n=1}^\infty \Cu^n[Z] \frac{t^n}{n!} \ . \end{equation}

For determinantal processes, it turns out that there are explicit formulae,  in terms of the correlation kernel $K$, for the cumulants  of a linear statistic $\Xi f = \sum f(X_i)$. Taking $g(x)=e^{t f(x)}-1$ for some function $f\in C(\X)$ with compact support and $t\in\R$ in  equation (\ref{Fredholm_det}), we see that
\begin{equation*}\E{K}{e^{t \Xi f }}= \det\left[\I+\K(e^{t f(x)}-1)\right] .
\end{equation*} 
Since the operator $\K$ is assumed to be locally trace-class, the RHS of this equation is a Fredholm determinant and taking logarithm (see for instance~\cite[chap.~3]{Simon_05}), we obtain
\begin{align*} \log\E{K}{e^{t \Xi f}}&= \tr \left[ \log\left(\I+\K(e^{t f(x)}-1)\right) \right] \\
&=\sum_{l=1}^\infty \frac{(-1)^{l+1}}{l} \tr\left[\big(\K(e^{t f(x)}-1)\big)^l\right] .
\end{align*} 

If we expand $e^{t f(x)}-1=\displaystyle\sum f(x)^n\frac{t^n}{n!}$ and use linearity of $\tr$, we deduce that the cumulants of the random variable $\Xi f$ are given by
 \begin{equation} \label{cumulant_1} 
 \Cu^n_K[\Xi f]= \sum_{l=1}^n \frac{(-1)^{l+1}}{l} \sum_{\begin{subarray}{c} m_1, \dots, m_l \ge 1 \\ m_1+ \cdots + m_l = n  \end{subarray}  } \frac{n!}{m_1!\cdots m_l!}
  \tr[f^{m_1}\K f^{m_2}\cdots f^{m_l} \K]  \ ,
  \end{equation}
where we interpret $f^{m_j}$ as multiplication operators acting on $L^2(\X,d\mu)$.  In particular, we have
 \begin{equation} \label{cumulant_4} 
\tr[f^{m_1}K f^{m_2}\cdots f^{m_\ell} K] =
  \int_{\X^\ell} f(x_1)^{m_1}K(x_1,x_2)\cdots f(x_\ell)^{m_\ell} K(x_\ell,x_1) d\mu(x_1)\cdots d\mu(x_\ell) , 
  \end{equation}
  so that, provided the precise asymptotics of the correlation kernel $K$ is available, we can deduce  from formula  \eqref{cumulant_1}  the limit law of the linear statistic $\Xi f$. 
For instance, we get a CLT with variance $\sigma^2$ if for any $n>2$, 
  \begin{equation*}
  \lim_{N\to\infty}  \Cu^n_K[\Xi f]= 0   
 \hspace{.7cm} \text{and} \hspace{.7cm}
  \lim_{N\to\infty}  \Cu^2_K[\Xi f]=  \sigma^2 .    
\end{equation*}

A {\it composition} of a number $n\in \N$ is a tuple ${\bf m}=(m_1,m_2,\dots, m_\ell)$ of positive integers such that $|{\bf m}|= m_1+\cdots+ m_\ell =n$, where $\ell=\ell({\bf m})$ is called the length of ${\bf m}$. 
Using the notation
\begin{equation} \label{M} 
 { n \choose {\bf m}} = \frac{n!}{\prod_j {\bf m}_j !} 
 \hspace{.7cm} \text{and} \hspace{.7cm}
\M({\bf m})= \frac{(-1)^{\ell+1}}{\ell} { n \choose {\bf m}}  , 
\end{equation} 
it will be convenient to rewrite formula (\ref{cumulant_1}) as
 \begin{equation} \label{cumulant_2} 
 \Cu^n_K[\Xi f]= \sum_{|{\bf m }|=n} \M({\bf m})  \tr[f^{m_1}K f^{m_2}\cdots f^{m_\ell} K]  .
  \end{equation}

\subsection{Main results} \label{sect:intro_3}


In this section, we summarize the main results of sections~\ref{sect:Modified_GUE} and~\ref{sect:Critical} about the asymptotics of sufficiently smooth linear statistics of the modified GUEs. These results are summarized in the diagram of figure \ref{fig:phase} below. \\

For any function $f \in L^1 \cap L^2(\R)$, we define its Fourier transform
  \begin{equation*} \hat{f}(u)= \int_\R f(x) e^{-i2\pi x u} dx \ . \end{equation*}
 We will consider the following spaces of test functions:\\
- $H^{1/2}(\R)$ denotes the Sobolev space of real-valued $L^2$-functions equipped with the norm
 \begin{equation}  \label{norm_1/2}
  \| f \|^2_{H^{1/2}}=  \int_\R \left| \hat{f}(u) \right|^2 |u| du  
  =\frac{1}{4\pi^2} \iint_{\R^2}  \left| \frac{f(x)-f(y)}{x-y} \right|^2 dxdy \ . \end{equation}
 - $H^{1}(\R)$ denotes the Sobolev space of real-valued $L^2$-functions equipped with the norm
 \begin{equation}  \label{norm_1}
  \| f \|^2_{H^1}=  \int_\R\left| \hat{f}(u) \right|^2 |u|^2du  
  = \frac{1}{4\pi^2} \int_\R \left| f'(x) \right|^2 dx \ . \end{equation} 
 Modulo constants, the spaces  $H^{1/2}(\R)$ and $H^{1}(\R)$ are complete normed. 
 Moreover, we denote by $C_0(\R), H^{1/2}_0(\R)$, etc,  the corresponding subspaces of compactly supported functions. 
 


\begin{theorem} \label{thm:Poisson_CLT}
 For any parameters $0<\delta<\alpha <1$ and for any bounded function $f\in H^{1/2}_0(\R)$,
 \begin{equation} \label{variance_1}  
 \Var_{K_{\Psi,\alpha}^N}[\Xi f_\delta ] = \frac{\tau}{2} N^{\alpha-\delta} \int_\R f(x)^2 dx +\underset{N\to\infty}{o  }(N^{\alpha-\delta}) . 
 \end{equation}
 This asymptotics implies the following classical central limit theorem as $N\to\infty$,
\begin{equation} \label{Poisson_CLT}  
N^{- \frac{\alpha-\delta}{2}}\ \tilde \Xi f_\delta\  
\Rightarrow\ \No\left(0, \tau \|f\|_{L^2}^2/2\right)  . 
\end{equation} 
\end{theorem}
\proof The asymptotic expansion of the variance is proved in section~\ref{sect:Variance}. The CLT (\ref{Poisson_CLT}) follows directly from the estimates (\ref{LS_mean}) and  (\ref{variance_1}) by applying Theorem 1 in~\cite{Soshnikov_01}. \qed\\

Hence, we will call the set $\{ \delta \in [0,1) : \delta<\alpha\}$ the {\bf Poisson scales} because the variance of any linear statistic is diverging in this regime.
Viewing the process $\tilde \Xi$, see \eqref{LS_2},   as a random distribution acting on $C^\infty_0(\R)$,  theorem
~\ref{thm:Poisson_CLT} implies that, once normalized, it converges at any scale $\delta<\alpha$ to a
 white noise with intensity $\tau/2$, i.e.~a centered Gaussian field $\Xi_\infty$ on the real line with covariance:
$$ \E{}{\Xi_\infty(f)\Xi_{\infty}(g) }= \frac{\tau}{2}  \langle f, g\rangle_{L^2(\R)}  . $$


At scales $ \delta\ge\alpha$, the variance remains bounded. Therefore, we expect a limiting process with non-trivial correlations.
Actually, by comparing linear statistics of the modified GUEs to that of the sine process, we will obtain the following CLT.

\begin{theorem}  \label{thm:MNS_CLT}  Let  $\Psi \in\F^*$ and  $f\in H^{1/2}_0(\R)$ be a bounded function. If $0<\alpha<\delta \le1$, then
\begin{equation*} \tilde\Xi f_\delta\  \Rightarrow\ \No\left(0, \|f\|_{H^{1/2}}^2\right)  , \end{equation*}   
 as $N\to\infty$, where the norm $\|f\|_{H^{1/2}}^2$ is given by formula $(\ref{norm_1/2})$.
\end{theorem}
\proof Section~\ref{sect:MNS_CLT}. \qed\\

Hence, we will call the  set $\{\delta \in (0,1) : \alpha<\delta\}$ the  {\bf GUE scales} by analogy with  theorem~\ref{thm:Szego_2} below.
The interpretation of theorem~\ref{thm:MNS_CLT} is that the centered modified GUEs  converge
 weakly at any scale $\delta>\alpha$ to a Gaussian process $\Xi_0$ on the real line with covariance 
$$ \E{}{\Xi_0(f)\Xi_0(g) }= \langle f, g\rangle_{H^{1/2}}   .$$
In contrast to the  white noise, the Gaussian process $\Xi_0$ is spatially correlated and self-similar as can be seen from equation (\ref{norm_1/2}).\\


Theorems ~\ref{thm:Poisson_CLT}  and ~\ref{thm:MNS_CLT} imply that the modified GUEs undergo a transition from Poisson to random matrix statistics when the mesoscopic scale $\delta$ is equal to the modification scale $\alpha$ of its correlation kernel. Our next question is what happens at the critical scale?
The first step is to investigate the variance of linear statistics. 

\begin{theorem} \label{thm:modified_variance}
For any shape $\Psi\in\F^*$, for any bounded function $f \in H^{1/2}_0(\R)$, and for any scale $0<\alpha<1$,
\begin{align} \notag 
\lim_{N\to\infty} \Var_{K_{\Psi,\alpha}^N}[\Xi f_\alpha ] 
&= 2 \tau' \int_\R \left| \hat{f}(u) \right|^2 \int_\R \Psi(t)\big( 1-\Psi(t+u/\tau') \big)\ dudt \\
& \label{modified_variance}
=2\tau'\int f(x)^2 dx \   
+\frac{1}{4\pi^2}  \iint \left| \frac{f(x)-f(y)}{x-y}\right|^2\left|\hat\Phi\big(\tau'(x-y)\big)\right|^2  dxdy  ,
\end{align}
where $\Phi=-\Psi'$ and the parameter $\tau'=\frac{\tau}{4}>0$.
\end{theorem}
\proof  Appendix \ref{A:variance}.   \qed\\

Since $\Phi$ is a probability distribution function, it is clear that formula (\ref{modified_variance}) interpolates between $\|f\|^2_{H^{1/2}}$ as $\tau \to 0$, respectively $2\tau' \|f\|_{L^2}^2$ as $\tau\to\infty$. In both cases, we recover the variances of theorem ~\ref{thm:Poisson_CLT} and theorem ~\ref{thm:MNS_CLT}   respectively. In analogy,  we would also expect Gaussian fluctuations when $\delta=\alpha$.
 Surprisingly, at this critical scale, the cumulants of  linear statistics of the modified ensembles have non-trivial limits. In order to formulate our main result, we need to introduce additional notation. Let
 $$\R_<^n=\{ x\in\R^n : x_1<\cdots<x_n \} ,$$
and for any $z\in\R$,
$$\R_z^n=\{x\in\R^n : x_1+\cdots+x_n=z\}  .$$ 
 
For any composition {\bf m} and for any $k \le \ell(\bf {m})$, we let
\begin{equation} \label{N_2}
\overline{ m}_k = m_1+  m_2+ \cdots + m_k  .
\end{equation}  
For any $u\in \R^{|{\bf m}|}$, we define
\begin{equation} 
\label{Lambda}
\Lambda^{\bf m}_{i, s}(u)= \sum_{j=1}^{{\overline{ m}_i}} u_j -  \sum_{j=1}^{\overline{ m}_{s}}u_j 
=\begin{cases} u_{\overline{ m}_{s}+1}+\cdots+ u_{\overline{ m}_i} &\text{if } s< i \\
-u_{\overline{ m}_{i}+1}-\cdots- u_{\overline{ m}_s} &\text{if } i<s \\
0 &\text{if } i=s
\end{cases}  ,
 \end{equation}
so that $\Lambda^{\bf m} = \big( \Lambda^{\bf m}_{i, s}\big) $ is a $\ell({\bf m})\times \ell({\bf m})$  antisymmetric matrix.\\

If $\sigma \in \Sy(n)$ is a permutation of $[n]=\{1,\dots,n\}$, we will use the shorthand notation $\sigma u= (u_{\sigma(1)},\dots, u_{\sigma(n)})$ and we define $\s(\sigma)=\arg\min(\sigma) \in [n]\times[n-1]\times \dots \times[1] $ as follows. For any $l=1,\dots, n$,  the number $\s_l(\sigma)$ is given implicitly by the relation
\begin{equation} \label{N_3} 
 \sigma(\s_l)= \min\{\sigma(j) : j=1,\dots, l\}  . \end{equation}

 For any $\tau>0$ and any composition ${\bf m}$ of $n \ge 2$, let
 \begin{equation} \label{G}
  G^{\bf m}_{\tau}(u,x) = \sum_{\sigma\in \Sy(n)} \max_{i \le \ell}\left\{\Lambda^{\bf m}_{i,\s_\ell}(u) -\tau(x_{\sigma(i)} - x_{\sigma(\s_\ell)}) \right\} ,
  \end{equation}
 where we used the shorthand notation $\ell=\ell({\bf m})$ and $\s_\ell=\s_{\ell({\bf m})}(\sigma)$.
Note that, since $\Lambda^{\bf m}_{\s_\ell, \s_\ell}(u)=0$,  the functions $G^{\bf m}_{\tau}(u,x)$ are  non-negative on $\R^n\times \R^n_<$. \\

 Finally, for any function $\Psi\in\F$, let  $\B^1_\Psi=0$ and for any $n \ge 2$,
\begin{equation}\label{B} 
\B^n_\Psi= \sum_{k=0}^{n-1}   b_k^n \int_\R  x\Phi(x) \Psi(x)^k\big(1-\Psi(x)\big)^{n-1-k} dx  ,
\end{equation}
where the coefficients $b^n_k$ are given by
\begin{equation}\label{b} 
b^n_k=  \sum_{l=1}^{ k+1} (-1)^{l+1} {n- l \choose k+1-l}
 \sum_{\begin{subarray}{c}  |{\bf m}|=n \\ \ell({\bf m})=l \end{subarray}} {n \choose {\bf m}}  .
  \end{equation}
Note that since  $b^2_0=-b^2_1=1$, by formula (\ref{B}),
$$ \B^2_\Psi = \int_\R x \Phi(x)\big(1- 2\Psi(x) \big) dx   .$$
Then, using that for any shape $\Psi \in \F$, we have 
$$\Phi(x)\big(1- 2\Psi(x) \big) = -\frac{d}{dx}\big\{\Psi(x)\big(1-\Psi(x)\big) \big\} ,$$
we recover formula \eqref{B^2} by integration by part. \\

\begin{theorem} \label{thm:MNS_C^n}
Let $0<\alpha<1$, $\Psi\in\F^*$, and consider the determinantal process with correlation kernel $K_{\Psi,\alpha}^N$. For any bounded function $f\in H^{1/2}_0(\R)$, the centered linear statistic $\tilde\Xi f_\alpha$  converges in distribution as the density $N\to \infty$ to a random variable denoted  $\Xi_{\Psi,\tau/4} f$.  Moreover, if  $f\in H^{1}_0(\R)$, then the cumulants of $\Xi_{\Psi,\tau} f$ are given by 
\begin{equation} \label{MNS_C^n}
 \Cu^n\big[\Xi_{\Psi,\tau} f \big]=  
2\tau \B^n_\Psi \int_{\R} f(t)^n dt\ - 2\sum_{|{\bf m}|=n} \M({\bf m})
\underset{\R^n_0} {\int du}  \underset{\R^n_<}{\int dx}\ \Re\left\{ \prod_{i=1}^n \hat{f}(u_i)\Phi(x_i) \right\} \G_\tau^{\bf m}(u,x)   \end{equation}
for all $n\ge 2$. 
\end{theorem}
 \proof Section~\ref{sect:C^n}. \qed


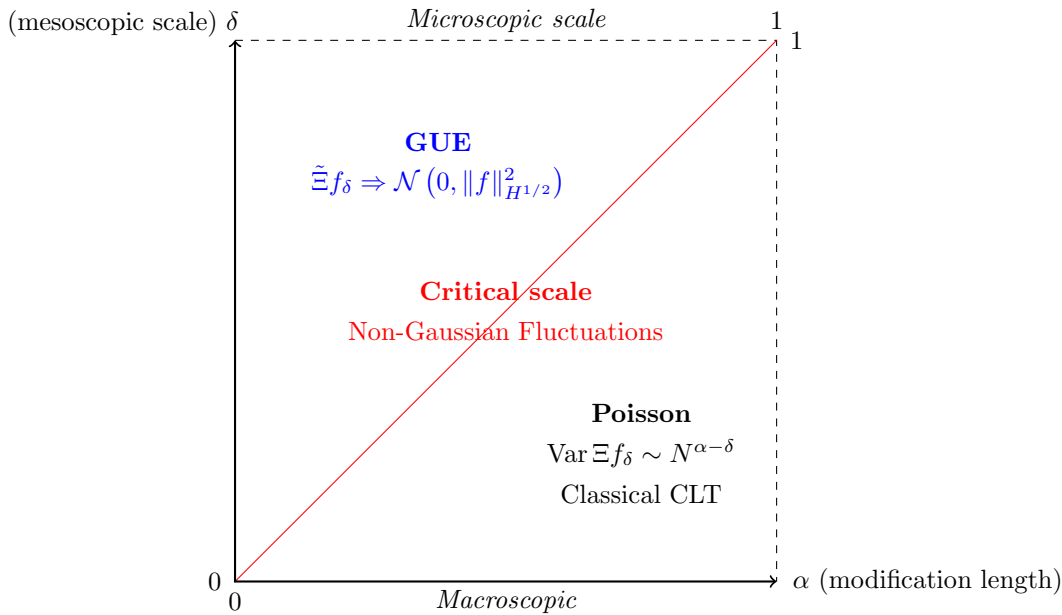
\begin{figure}[h]
\begin{center}
\begin{tikzpicture}[scale = 0.9]
\draw  [dashed]  (0,8) -- (8,8) -- (8,0);
\draw[red]  (0,0) -- (8,8) ;
\draw [thick, <->] (0,8) -- (0,0) -- (8,0);
\node at (-.3,0) {0};
\node at (8.3,8) {1};
\node at (8,8.3) {1};
\node at (0,-.3) {0};
\node at (3,6.5)[blue] {\bf GUE};
\node at (3,5.9)[blue] {$ \tilde\Xi f_\delta \Rightarrow\No\left(0, \|f\|^2_{H^{1/2}} \right)$};
\node at (6,2.5) {\bf Poisson};
\node at (6,1.9)  {$\Var \Xi f_\delta \sim   N^{\alpha-\delta}$};
\node at (6,1.3) { Classical CLT};
\node at (4,4.3)[red] {\bf Critical scale};
\node at (4,3.7)[red] {Non-Gaussian Fluctuations};
\node at (.2,8.25) [left] {(sampling scale) $\delta$ };
\node at (8.1,0) [right] {$\alpha$ (modification scale)};
\node at (4,-.3){\it Macroscopic scale};
\node at (4,8.3){\it Microscopic scale};
\end{tikzpicture}
\caption{Fluctuations of a Modified GUE  as a function of the scale $\delta$ and the modification scale~$\alpha$. \label{fig:phase}}
\end{center}
\end{figure}

It is a difficult problem to recover from equation (\ref{MNS_C^n}) the Laplace transform of the random variable $\Xi_{\Psi,\tau} f$. In fact, we can only infer a few properties of its distribution, such as the fact that it is not Gaussian and the dependence in the temperature $\tau$.

\begin{theorem} \label{thm:Normal}
For any shape $\Psi\in\F$ and any $\tau>0$, the random process $\Xi_{\Psi,\tau} $ which arise from the limit of a modified Ensemble at the critical scale $\delta=\alpha$  is not Gaussian.
\end{theorem}

The proof of theorem \ref{thm:Normal} is  rather complicated and divided into two parts. First, we show that, if  there exists $n>2$ such  that $\B^n_\Psi \neq 0$, then the random process  $\Xi_{\Psi,\tau} $ is not Gaussian; see proposition~\ref{thm:Gaussian} below. Surprisingly, this condition is satisfied by all modified GUEs, except the MNS Ensemble.

\begin{proposition} \label{thm:MNS_property}
The only  function $\Psi \in \F$ which satisfies the  conditions $ \B^n_\Psi = 0$ for all $n>2$ is the MNS shape  $\psi(t)=1/(1+e^t)$.
\end{proposition}
 \proof Section~\ref{sect:Poisson}. \qed\\
 
This special property  of the MNS Ensemble is maybe related to the fact that this point process originates from a {\it Grand Canonical}  model and it would be interesting to know whether it has any physical interpretation. 
 From a probabilistic perspective, the main consequence of proposition~\ref{thm:MNS_property}  is that we need a separate argument to show that the  process $\Xi_{\psi,\tau} $ is not Gaussian.
To this end, we  show in section~\ref{sect:C3C4} that the Schwartz function $ y(x)=2e^{-\epsilon \pi x^2}\cos(2\pi x)$ satisfies,  when $\epsilon$ is sufficiently small,
 \begin{equation} \label{G4'} 
 \Cu^4\big[\Xi_{\psi,\tau} y \big] \neq 0 . 
 \end{equation} 
We have to look at the $4^{\text{th}}$ cumulant since $ \Cu^3\big[\Xi_{\psi,\tau} f \big] = 0$ for any test function because of some symmetries; see proposition~\ref{thm:C3}. Let us also mention that there is nothing special about the function $ y(x)$ except that it is simple enough to provide a good example. \\

We have seen that the MNS  ensemble describes non-colliding Brownian motions in a cylindrical geometry. In this case, the diagram of figure \ref{fig:phase} shows that the particle statistics at a given mesoscopic scale~$\delta$ exhibit a sharp transition from Poisson to GUE at time $N^{-\delta}$. This situation is similar to that of the Dyson's Brownian motion which was investigated in~\cite{DJ_14}. However the transitions are different: for the MNS ensemble, there is no intermediate regime which depends on the test function $f$ and quite unexpectedly the critical fluctuations are not Gaussian. 
In both interpretations, either as non-colliding Brownian motions, or as a grand canonical ensemble of free fermions, it is not clear why this phenomenon occurs at a
certain relation between the sampling scale of the linear statistic and the size of cylinder, or  the temperature of the gas respectively. 
Actually, it would be very interesting to get another description than theorem~\ref{thm:MNS_C^n}  of the random field~$\Xi_{\psi,\tau}$ which arises at the critical scale in order to understand how the spatial correlations disappear in the transition from the $H^{1/2}$-Gaussian noise $\Xi_0$ to the white noise $\Xi_\infty$.
Another point of interest would be to understand how the processes $\Xi_{\Psi,\tau}$ depend on the shape $\Psi$ and why the MNS Ensemble appears to be special.\\

\subsection{Outline of the main ideas} \label{sect:counter_examples} 


For determinantal processes, a general strategy to obtain a CLT for linear statistics is to use formula (\ref{cumulant_1}) to show that  all cumulants of order $\ge 3$ converge to zero. This approach has been applied to many classical models in 1 or 2 dimensions using a wide range of techniques,  see e.g.~\cite{CL_95,Soshnikov_00a,Soshnikov_01, RV_07a, RV_07b, AHM_11,BD_14, BD_15}. In this paper, it is applied to the modified ensembles of definition \ref{modified_GUE}. To the authors' knowledge, the results of this paper at the critical scale provide the first example of  a determinantal point process for which the fluctuations of a mesoscopic linear statistic are not Gaussian in the limit. At the global scale, non-Gaussian limits can be obtained for unitary invariant Hermitian matrix models with several cuts, but the mechanism is different. In fact, it was recently proved in~\cite{L_15a} that,  at mesoscopic scales, such models are also described in the large $N$ limit by the $H^{1/2}$-Gaussian noise $\Xi_0$ appearing in theorem~\ref{thm:MNS_CLT}. 
To prove theorems ~\ref{thm:MNS_CLT} and~\ref{thm:MNS_C^n},  we use a perturbative approach which consists in comparing the correlation kernels of two processes to establish that a given linear statistic $\Xi f$ has the same limit for both ensembles as the density $N\to\infty$.  To this end, we will use the following definition.


\begin{definition} \label{def_equivalent_kernel} 
 Two families of kernels $(K^N)_{N>0}$ and $(L^N)_{N>0}$ defined on the same space $\X$ are {\bf asymptotically equivalent} $($we write $K^N \cong L^N)$ if, for any $\ell\in\N$ and any functions $f_1,\dots, f_\ell \in C_0^\infty(\X)$,
\begin{equation*}  \lim_{N\to\infty} \tr[f_1K^N \cdots f_\ell K^N] =  \lim_{N\to\infty}  \tr[f_1L^N \cdots f_\ell L^N]   .  \end{equation*}
\end{definition}

Definition \ref{def_equivalent_kernel} does not apply only to correlation kernels, but to all 
families of kernels which are locally trace-class.  However, if both $K^N$ and $L^N$ define determinantal processes and $K^N \cong L^N$, by formula (\ref{cumulant_2}),  these processes have the same limit as $N\to\infty$.
For instance, theorem ~\ref{thm:MNS_CLT} is proved  by showing that the kernel $K_{\Psi,\alpha}^N$  and the GUE kernel $K_0^N$ are asymptotically equivalent at any scale $\delta > \alpha $
 and using the CLT for the mesoscopic GUE  (see theorem~\ref{thm:Szego_2} proved in~\cite{FKS_13, BEYY_14,L_15a}).
On the other hand, at the critical scale, $\delta=\alpha$, the kernels $K_{\Psi,\alpha}^N$ are not asymptotically equivalent to any kernel which has been studied previously and we will need to compute the limits of the cumulants explicitly. 
If $\Psi\in\F$, we let $\Phi=-\Psi'$ and define $\xi_k \in(0,1)$ such that by the mean-value theorem:  $\frac{1}{\tau N^\alpha} \Phi\left(\frac{k +\xi_k}{\tau N^\alpha}\right)=\Psi\left(\frac{k}{\tau N^\alpha}\right)- \Psi\left(\frac{k+1}{\tau N^\alpha}\right)$ for any $k\in\Z$. Then, if  $N,\tau,\Gamma>0$, $\alpha \in (0,1)$, and $\eta$ is a non-decreasing function, we define the kernel

\begin{equation} \label{kernel_L}
L_{\Psi,\eta}^N(x,y)=\frac{1}{\tau N^{\alpha}} \sum_{|k| \le \Gamma N^\alpha } 
 \Phi\left(\frac{k+\xi_k}{\tau N^\alpha}\right) \frac{\sin\big[2\pi \eta(k)(x-y)\big]}{\pi(x-y)}  .
 \end{equation}
 
We show that the kernel $L_{\Psi,\eta}^N$ defines a determinantal process on $\R$, see  lemma~\ref{thm:kernel_L}, and the next proposition  implies that, at sufficiently small scales, mesoscopic linear statistics of the critical modified GUE with shape $\Psi$ and the determinantal process with kernel $L_{\Psi,\eta}^N$ have the same limit.

\begin{proposition} \label{thm:GUE_kernel_sine}  Let $\Psi\in\F^*$, $\eta(k)=N^{1-\alpha} \sqrt{1+k/N}/2$ and $\Gamma=(\log N)^2$.  For all  $1/3<\alpha<1$,   the rescaled correlation kernel of the modified GUE satisfies
$$N^{-\alpha}K_{\Psi,\alpha}^N(xN^{-\alpha},yN^{-\alpha}) \cong L_{\Psi,\eta}^N(x,y)$$
 in the sense of definition~$\ref{def_equivalent_kernel}$.
\end{proposition}
\proof Section~\ref{sect:MNS_Critical}. \qed\\

 Proposition~\ref{thm:GUE_kernel_sine}  combined with the analysis of the determinantal process with correlation kernel $L_{\Psi,\eta}^N$ performed in section~\ref{sect:MNS_Critical} imply theorem~\ref{thm:MNS_C^n} in the regime $1/3<\alpha<1$.
The main technical challenge of this argument  is to get the asymptotic expansion of the GUE kernel at mesoscopic scales. Namely we show in the appendix~\ref{sect:Asymptotic} that, if $M\sim N$, then  for any $\delta>0$, 
\begin{equation}  \label{sine_approximation_2}
N^{-\delta}K_0^M(xN^{-\delta},yN^{-\delta}) =  \frac{\sin\big[ \pi N^{1/2-\delta}\sqrt{M} (x-y)\big]}{\pi(x-y)} + \underset{N\to\infty}{O} \left( N^{1-3\delta} \right) 
\end{equation}
 uniformly for all $x,y$ in compact subsets of $\R$.
 Note that when $\alpha \le 1/3$, the error term in formula~\eqref{sine_approximation_2} does not converge to 0 and this gives the restriction in proposition~\ref{thm:GUE_kernel_sine}. This restriction comes from the fact the GUE kernel is not asymptotically translation-invariant at such scales because of the curvature of the density of the semicircle law and the proof of theorem~\ref{thm:MNS_C^n} in the general case relies on a different argument which basically consists in {\it unfolding the point process} to reduce again  to the case of the kernel  $L_{\Psi,\eta}^N$, see proposition~\ref{thm:MNS_trace}.
The advantage of this approach is that  the kernel $L_{\Psi,\eta}^N$ is  translation-invariant and we can compute the cumulants of its linear statistics by using the method  introduced in~\cite{Soshnikov_00a} to prove a CLT for mesoscopic linear statistics of the Circular Unitary Ensemble  (see also theorem 4 in~\cite{Soshnikov_01} for an application to the sine process, as well as some generalizations). 
In fact, taking the parameter $\tau \to 0$ in equation~(\ref{MNS_C^n}), we recover Soshnikov's formula for the mesoscopic sine process:
\begin{equation*}  \Cu^n\big[\Xi_{\Psi,0} f \big]=  
 2 \underset{\R^n_0} {\int du}\ \Re\left\{ \prod_i \hat{f}(u_i) \right\} \sum_{|{\bf m}|=n} \M({\bf m}) \G^{\bf m}_0(u)  ,
\end{equation*}
where 
$\displaystyle \G^{\bf m}_0(u)= \sum_{\pi \in \Sy(n)} \max\left\{ u_{\pi(1)}+\cdots+u_{\pi( m_1)},\dots,\ u_{\pi(1)}+\cdots+u_{\pi(m_1+\cdots +m_{\ell-1})},\ 0 \right\}$.

Then, the main combinatorial lemma of~\cite{Soshnikov_00a} implies that 
$\displaystyle \sum \M({\bf m}) \G^{\bf m}_0(u)=0$
 for any  $n > 2$, so that the process $\Xi_{\Psi,0}$ (which is independent of $\Psi$) is Gaussian. The details of the computations are given in the proof of proposition~\ref{thm:double_limit}. 
For the modified ensembles, we observed that there is no counterpart of the main combinatorial lemma, i.e. for generic points $u\in\R^n_0$ and $x\in \R^n_<$, 
$\displaystyle \sum \M({\bf m}) \G^{\bf m}_\tau(u,x) \neq 0$ for any $\tau>0$ and $n\neq 3$. Then we use this fact to prove equation (\ref{G4'}); see in particular lemma~\ref{thm:G4}.
%
%
The bottom line is  that the combinatorial structure behind the cumulants of the sine process, which corresponds to the continuous counterpart of the Strong Szeg\H{o} theorem, is very sensitive. In general, CLTs with bounded variance are due to some special correlation structures  which are rather sensitive under perturbation such as some small modification of the correlation kernel. In the remainder of this section, we provide a basic example which elaborates on this fact and illustrate how asymptotic normality breaks down. Before proceeding, we define the {\it circular} counterparts of the modified GUEs. These point processes  are of interest because asymptotic expansions are not required in order to apply Soshnikov's method and they retain the same  features as the modified GUEs. In the sequel, we let  $\T=[-\frac{1}{2},\frac{1}{2}]$ with the boundary points identified.

\begin{definition} \label{modified_CUE} 
A {\bf modified CUE} is a determinantal process on $\T$ $($wrt$.$ the Lebesgue measure$)$ whose correlation kernel $K^N_p$ is of the form 
\begin{equation}  \label{modified_CUE_1}
K^N_{p}(x,y)= \sum_{k\in\Z} p_k^N e^{i 2\pi k (y-x)} ,
\end{equation}
where $p_{-k}^N=p_k^N$ so that the corresponding integral operator  is self-adjoint on $L^2(\T)$.
\end{definition}

\begin{remark}
In this section, the spectrum of the kernel $K^N_p$ is arbitrary except for the constraints $p_k^N \in [0,1]$. However, when we refer to a modified CUE in sections~\ref{sect:CUE_CLT}, \ref{sect:CUE_critical} and~\ref{sect:Critical}, like in definition~\ref{modified_GUE}, it is understood that
\begin{equation} \label{modified_CUE_2}
p^N_k=\Psi\left(\frac{|k|-N}{\tau N^\alpha}\right) 
\end{equation}
for a shape $\Psi \in \F$ such that $\B^2_\Psi=1$, a modification scale $\alpha \in (0,1)$, and a temperature $\tau>0$.
\end{remark}

A special case is the so-called Dyson's Circular Unitary Ensemble (CUE) which has the correlation kernel
\begin{equation} \label{CUE_kernel} 
 K^N(x,y)= \sum_{|k|\le N} e^{i 2\pi k (x-y)}= \frac{\sin((2N+1)\pi (x-y))}{\sin(\pi(x-y))}  .  
 \end{equation}
This process describes the eigenvalues of a random matrix distributed according to the Haar measure on the group $\U(2N+1)$. The cumulants of its linear statistics were computed explicitly in~\cite{Soshnikov_00a} and a similar computation yields the following formula in the case of the modified CUEs.
 
 \begin{lemma} \label{thm:CUE_cumulants}
For any continuous function $f:\T\to\R$,
\begin{equation*}
\Cu^n_{K^N_p}[\Xi f]=  \sum_{u \in\Z^n_0} \prod_i \hat{f}(u_i) \sum_{|{\bf m}|=n} \M({\bf m}) \sum_{k\in\Z}\ \prod_{i=0}^{\ell(m)-1} p^N_{k+u_1+\cdots+u_{\overline{ m}_i}}  ,
 \end{equation*}
 where $\hat{f}(k)$ is the $k^{\text{th}}$ Fourier coefficient  of the function $f$. 
 \end{lemma}
 \proof Section~\ref{sect:Soshnikov}.  \qed\\

We can use this formula  to  investigate the behavior of (global) linear statistics under some very simple modification of the spectrum of the CUE correlation kernel. 
For instance, we can remove a single mode, i.e$.$  we let
 \begin{equation}  \label{mode_remove}
 p_{k}^N= \1_{|k|\le N} - \1_{|k|=N-m} 
  \end{equation}
for some $0<m\le N$. In this case, $K_p^N$ is still a projection kernel and,
by lemma~\ref{thm:CUE_cumulants}, the $3^{\text{rd}}$ cumulant of a linear statistic is given by
\begin{equation*}
\Cu^3_{K^N_p}[\Xi f]=  \sum_{u\in\Z^3_0} \prod_i \hat{f}(u_i) \left\{ \sum_{k\in\Z} p_k - \frac{3}{2}\sum_{k\in\Z} p_k(p_{k+u_1}+p_{k+u_1+u_2})+ 2  \sum_{k\in\Z} p_kp_{k+u_1}p_{k+u_1+u_2} \right\}  .
 \end{equation*}

We can symmetrize this expression using permutations of the $u_i$'s and the condition $u_1+u_2=-u_3$, this yields
\begin{equation} \label{CUE_C3}
\Cu^3_{K^N_p}[\Xi f]= \sum_{u\in\Z^3_0} \prod_i \hat{f}(u_i) \sum_{k\in\Z} p_k  \left\{ 1 - 3p_{k+u_1}+ 2 p_{k+u_1}p_{k-u_2} \right\}  .
 \end{equation}
 Let us consider the function $g_j(t)=2\cos(2\pi j t)+ a\cos(4\pi j t)$ for some parameters $ a\in\R$ and $j\in\Z_+$,  so that 
$$\hat{g_j}(u)= \delta_j(|u|)+  \frac{a}{2} \delta_{2j}(|u|) \ . $$
For this test function, the only {\it frequencies} $u\in\Z_0^3$ which contribute to (\ref{CUE_C3}) are given by all possible permutations of $(\pm j, \pm j, \mp 2j)$ and an elementary computation shows that
\begin{equation} \label{CUE_C3_example} 
\Cu^3_{K^N_p}[\Xi  g_j] = 3 a \sum_{k\in\Z} p_k \left\{ 1  - 2p_{k+j} -p_{k+2j}+ 2p_{k+j}p_{k-j}    \right\}  .
 \end{equation}
 
In the CUE case, when $p^N_k=\1_{|k|\le N}$, it is easy to check that $\Cu^3_{K^N}[\Xi  g_j] =0$ for any $ j \in \Z_+$ as we expected. However, it seems clear that for some generic choice of coefficients $0\le p_k^N \le 1$, the expression (\ref{CUE_C3_example})  will be non-zero. For instance in the case (\ref{mode_remove}), for any $j,m \ll N$,
\begin{equation*} 
\Cu^3_{K^N_p}[\Xi  g_j] =  12 a\big(1-  \1_{j\le \lfloor m/2\rfloor }\big) \ .   
\end{equation*}
Hence, if we remove a mode near the edge of the {\it spectrum} of the CUE kernel  ($m < 2j $), the linear statistics $\Xi g_j$ is not Gaussian in the limit $N\to\infty$. Moreover we can check that the variance is bounded:
\begin{equation*} \Cu^2_{K^N_p}[\Xi  g_j]= \sum_{u\in\Z} \left|\hat{g}_j(u)\right|^2 \sum_{k\in\Z} p_k(1 - p_{k+u}) 
\le j(2+a^2) \ .
\end{equation*}


This implies that, if it exists, the limit as $N\to\infty$ of the determinantal process with correlation kernel $K_p^N$ is not a Gaussian process. In particular, even if the correlation kernel is reproducing, we can get non-Gaussian behavior.
This example also shows that it is the {\it edge} of the spectrum of the correlation kernel which is influencing the distribution of the point process, see also theorem~4 in~\cite{Soshnikov_01}. Moreover, using  lemma~\ref{thm:CUE_cumulants}, we can also check that when removing $M$ different modes at the {\it edge} of the spectrum, all cumulants are of order $M$ for large $N$. Hence, if we remove sufficiently many modes, the system begins to behave like a Poisson process when~$N\to\infty$.
Finally, note that according to (\ref{Bernoulli}), removing modes is comparable to smoothing the spectrum of the correlation kernel. Hence, this example illustrates why the modified ensembles of definition~\ref{modified_GUE} are not Gaussian at the critical scale ($\delta=\alpha$).    
Actually, the strategy to obtain (\ref{G4'}) is the same as in this example but the computations are much more complicated.


\subsection{Overview of the rest of the paper}

In section~\ref{sect:Modified_CUE}, we begin by analyzing the modified CUEs of definition~\ref{modified_CUE}.  This setting is simpler than that of section~\ref{sect:intro_2} and we can focus on the combinatorial structure of the cumulants. 
In particular,  in sections~\ref{sect:CUE_CLT} and \ref{sect:CUE_critical},  we show that, if the spectrum of the kernel is given by  $p^N_k=\Psi\left(\frac{|k|-N}{\tau N^\alpha}\right)$, then the results of figure \ref{fig:phase} hold  for the modified CUEs as well. The main results of section~\ref{sect:intro_3} are proved in sections~\ref{sect:Modified_GUE} and~\ref{sect:Critical}.  The asymptotics of the variance in the Poisson regime, formula (\ref{variance_1}), is computed in section~\ref{sect:Variance}, while in the GUE regime,  theorem~\ref{thm:MNS_CLT}  is proved in section~\ref{sect:MNS_CLT}.
Both the critical modified CUEs and GUEs are analyzed, in a common framework, in section~\ref{sect:Critical}. In particular, the proof of theorem~\ref{thm:MNS_C^n} is divided in two steps. First the limits of the cumulants are established in theorem~\ref{thm:C^n} (see also proposition~\ref{thm:C^n'}). Then the weak convergence of linear statistics is established in corollary~\ref{thm:weak_convergence}. 
 In section~\ref{sect:Poisson}, we show that  the random processes $\Xi_{\Psi,\tau}$ defined by  theorem~\ref{thm:MNS_C^n} are not Gaussian and we prove the special property of the MNS ensemble, see theorem~\ref{thm:MNS_property},  by computing the generating function of the coefficient $\B^n_{\Psi}$. We also prove that, as it is expected from figure~\ref{fig:phase}, the random field $\Xi_{\Psi,\tau}$ converges to a Gaussian process in both limits $\tau\to0$ and $\tau\to\infty$; see proposition ~\ref{thm:double_limit}.
Finally, in section~\ref{sect:C3C4}, we show that the critical MNS ensemble is not Gaussian by constructing the example~\eqref{G4'}. All the asymptotics that are required to analyze the modified GUEs are gathered in section~\ref{sect:Asymptotic}. 
In appendix \ref{A:variance}, we prove theorem~\ref{thm:modified_variance} and, as an example, we compute the critical variance of linear statistics of the MNS ensemble. In appendix \ref{A:G}, we prove some technical lemmas which shows that there is no counterpart of the main combinatorial lemma for the modified ensemble and are used to prove  \ref{G4'},  
Finally, in the remainder of the paper, we will use the following conventions:\\
$\bullet$ $x_N\sim z_N$ if $\lim_{N\to\infty}  x_N / z_N =1$.\\
$\bullet$ $x_N\simeq z_N$ if $\lim_{N\to\infty}( x_N - z_N) =0$.\\
$\bullet$ $x_N=\bar{O}(z_N)$ if there exist $\kappa>0$ and $C>0$ such that 
$ |x_N| \le  C  z_N|\log N|^\kappa$.

\section{Modified Circular Unitary Ensembles} \label{sect:Modified_CUE}

We present the counterparts of the results of section~\ref{sect:intro_3} for the modified CUEs (definition~\ref{modified_CUE}). Along the way, we set up definitions and lemmas that will also be used in section~\ref{sect:Modified_GUE} and~\ref{sect:Critical}.  Circular ensembles can be thought of as  simplified models which are helpful to understand the combinatorial structure behind the cumulants of linear statistics of the MNS  model because no asymptotic estimates are required to pass to the limit. In section~\ref{sect:Soshnikov}, we review the method introduced in ~\cite{Soshnikov_00a}. In section~\ref{sect:CUE_CLT}, we show that the modified CUEs exhibit  the same transition as in figure \ref{fig:phase}. Finally, in section~\ref{sect:CUE_critical}, we provide asymptotically equivalent  kernels for the modified CUEs in the critical regime $\delta=\alpha$ and we deduce a limit theorem from the results of section~\ref{sect:Critical}.

 \subsection{Soshnikov's method: proof of lemma~\ref{thm:CUE_cumulants}} \label{sect:Soshnikov}

In~\cite[lemma 1]{Soshnikov_01}, Soshnikov proved that the cumulants of linear statistics of a determinantal process are given by
\begin{equation*} 
\Cu^n_K[\Xi f]= \sum_{|{\bf m }|=n} \M({\bf m})
  \tr[f^{m_1}K f^{m_2}\cdots f^{m_\ell} K] \ ,
  \end{equation*}
where $K$ is the correlation kernel of the process  (the sum is over all {\it compositions} ${\bf m}$ of $n$). See also  equation~(\ref{cumulant_2}) in section~\ref{sect:intro_5} for a formal derivation and an explanation of the notation.
Applying this formula to a modified CUE and using some elementary Fourier analysis we obtain lemma~\ref{thm:CUE_cumulants}.\\


{\it Proof of  lemma~\ref{thm:CUE_cumulants}.} 
The correlation kernel of the modified CUE is
$\displaystyle K^N_{p}(x,y)= \sum_{k\in\Z} p_k^N e^{i 2\pi k (y-x)} $,
and  for any composition ${\bf m}$ of $n$, by formula (\ref{cumulant_4}),
\begin{equation} \label{CUE_trace}
\tr[f^{m_1}K^N_{p}\cdots f^{m_j} K^N_{p}] = \sum_{\kappa\in\Z^{\ell({\bf m})}} \prod_{i=1}^{\ell({\bf m})} p^N_{\kappa_i} \widehat{f^{ m_i}}(\kappa_{i}-\kappa_{i-1})
 \end{equation}
 where by convention $k_0=k_{\ell({\bf m})}$.  For any indices $s, r \in \Z$ and $m\in \N$, we know that
 \begin{equation*}  \widehat{f^m}(s-r)= \sum_{k\in \Z^{m-1}} \hat{f}(k_1-r)\hat{f}(k_2-k_1)\cdots \hat{f}(s-k_{m-1}) \ .  \end{equation*}
Let  $\overline{ m}_j= m_1+\cdots+ m_j$  as in definition (\ref{N_2}).  For any $i=0,\cdots, \ell-1$, we can write
 \begin{equation*}  \widehat{f^{ m_{i+1}}}(k_{\overline{ m}_{i+1}}-k_{\overline{ m}_{i}})= \sum \hat{f}(k_{\overline{ m}_i+1}-k_{\overline{ m}_{i}})\hat{f}(k_{\overline{ m}_i+2}-k_{\overline{ m}_i+1})\cdots \hat{f}(k_{\overline{ m}_{i+1}}-k_{\overline{ m}_{i+1}-1})  \end{equation*} 
 and if we make the change of variables $\kappa_i=k_{\overline{ m}_i}$ in equation (\ref{CUE_trace}), putting everything together we get
 \begin{equation*}
\tr[f^{ m_1}K^N_{p} f^{ m_2}\cdots f^{ m_j} K^N_{p}]= \sum_{k\in\Z^{n+1}:k_0=k_n} \prod_{i=1}^n \hat{f}(k_{i}-k_{i-1})
  \prod_{i=1}^{\ell({\bf m})} p^N_{k_{\overline{ m}_i}} \ . 
 \end{equation*}
We can also make the change of variables $u_i=k_i-k_{i-1}$ for all $i=1,\dots, n$ in the previous sum. This maps $\{k\in\Z^{n}:k_0=k_n\}$ into $\{(k_0,u)\in\Z\times\Z^n_0\}$ and we obtain
 \begin{equation*}
\tr[f^{ m_1}K^N_{p} f^{ m_2}\cdots f^{ m_j} K^N_{p}]= \sum_{u\in\Z^n_0} \prod_i \hat{f}(u_i)
\sum_{k_0\in\Z}  \prod_{i=1}^{\ell({\bf m})} p^N_{k_0+u_1+\cdots+u_{\overline{ m}_i}} \ .
 \end{equation*}
Hence lemma~\ref{thm:CUE_cumulants} follows directly from formula (\ref{cumulant_2}).\qed\\


Observe that we recover lemma 1 in~\cite{Soshnikov_00a} for  Dyson's CUE by taking $p_k^N=\1_{|k| \le N}$, since then
\begin{equation} \label{H_sine}
 \sum_{k\in\Z}\prod_{i=0}^{l-1} p^N_{k+u_1+\cdots+u_{\overline{ m}_i}}
 = \left[ 2N+1- \underset{0\le i <l}{\max}\{u_1+\cdots+u_{\overline{ m}_i}\}  - \underset{0\le i <l}{\max}\{ -u_1-\cdots- u_{\overline{ m}_i}\}   \right]^+   \ .
 \end{equation}

In section~\ref{sect:counter_examples}, we used
lemma~\ref{thm:CUE_cumulants} to show that a particular modified CUE has non-Gaussian fluctuations at the macroscopic scale.  In the sequel,  we will use it to investigate fluctuations at mesoscopic scales.
Let $f\in C_0(\R)$, $0<\delta<1$, and recall that $\T=[-\frac{1}{2},\frac{1}{2}]$ with the endpoints identified. When the parameter~$N$ is sufficiently large,  the function $f(\cdot N^\delta)$ is supported in $[-\frac{1}{2},\frac{1}{2}]$ and it can be extended to some function $f_\delta\in C(\T)$.
Then, the Fourier coefficients of $f_\delta$ are given by, for any $u\in\Z$, 
\begin{equation*} \widehat{f_\delta}(u)= N^{-\delta} \hat{f}(uN^{-\delta} ) \ .    \end{equation*}
Hence,
\begin{equation}\label{CUE_cumulant_2}
C^n_{K^N_p}[\Xi  f_\delta]= N^{-n\delta}  \sum_{u\in\Z^n_0} \prod_i \hat{f}(u_i N^{-\delta})
\sum_{|{\bf m}|=n} \M({\bf m}) \sum_{k\in\Z}\ \prod_{i=1}^{\ell({\bf m})} p^N_{k+u_1+\cdots+u_{\overline{ m}_i}} \ .
 \end{equation}

\subsection{Central Limit Theorems} \label{sect:CUE_CLT}

From now on, we will assume that the spectrum of the modified CUE correlation kernel is given by
$ p^N_k =\Psi\left(\frac{|k|-N}{\tau N^\alpha}\right)$, see (\ref{modified_CUE_2}). Moreover, to keep the notation simple, we will write $p_k$ instead of $p_k^N$.\\\


We start by proving a classical CLT at the Poisson scales ($\delta<\alpha$). The proof relies on a simple variance computation. 
Observe that the asymptotic variance of theorem \ref{thm:CUE_Poisson} matches that of  theorem~\ref {thm:Poisson_CLT} only up to a multiplicative constant. The difference is due to our normalization. Namely, the scaling (\ref{modified_CUE_2}) implies that the modified CUEs have density $2N$ at the origin and
\begin{equation} \label{LS_mean_2}
 \E{K_p^N}{\Xi f_\delta}  = \left\{ 2N + \underset{N\to\infty}{O}(N^\alpha)\right\} N^{-\delta} \int  f(x)  dx \ .
 \end{equation}

\begin{theorem} \label{thm:CUE_Poisson}
Consider a modified CUE with correlation kernel $(\ref{modified_CUE_1}-\ref{modified_CUE_2})$ and let $f\in H^{1/2}_0(\R)$.
For any scale $0\le \delta <\alpha <1$, the centered and rescaled linear statistic $ N^{\frac{\delta-\alpha}{2}} \tilde\Xi f_\delta  $ converges in distribution to a Gaussian random variable with variance $2\tau \|f\|_{L^2(\R)}^2$.
\end{theorem}

\proof   
When $n=2$, equation (\ref{CUE_cumulant_2}) reads
\begin{equation} \label{CUE_C2}
C^2_{K^N_p}[\Xi f_\delta] = N^{-2\delta}\sum_{u\in\Z} \hat{f}(uN^{-\delta})\hat{f}(-uN^{-\delta}) \sum_{k\in\Z} p_k(1- p_{k+u}) \ . 
  \end{equation}
We let $\sigma_k^2=p_k(1-p_k)$ for any $k \ge 0$. Recall that $p_{-k}=p_k$,  then for any $u\in\Z$, 
\begin{align*}  \sum_{k\in\Z} p_k(1- p_{k+u})
& = p_0(1-p_u) +  \sum_{k>0} p_k \big(2- p_{k+u}-p_{k-u} \big)  \\
&= p_0(1-p_u) + 2\sum_{k>0} \sigma_k^2  +  \sum_{k>0} p_k \big(2p_k- p_{k+u}- p_{k-u} \big)
\end{align*}
Since the coefficients $p_k \in [0,1]$ and the shape $\Psi$ is non-increasing, we can check that for any $u\in\Z$,
 \begin{equation} \label{CUE_variance_estimate}
   \left| \sum_{k>0} p_k(p_k- p_{k+u}) \right| \le |u|
\hspace{.6cm}\text{so that}\hspace{.6cm}
\left| \sum_{k\in\Z} p_k(1- p_{k+u}) - 2\sum_{k>0} \sigma_k^2 \right| \le 2|u| +1 . 
\end{equation}
If we combine this estimate with formula (\ref{CUE_C2}), since the test function $f$ is real-valued,
 \begin{equation} \label{variance_split}
  \left| N^\delta\Cu^2_{K^N_p}[\Xi f_\delta] - 2 N^{-\delta}\sum_{u \in \Z} \left|\hat{f}(uN^{-\delta})\right|^2  \sum_{k>0} \sigma^2_k  \right| \le  N^{-\delta}\sum_{u \in \Z} \left|\hat{f}(uN^{-\delta})\right|^2 \big(2 |u|+1\big)    \ . 
   \end{equation}
If we assume that the test function $f\in H^{1/2}(\R)$, the r.h.s. satisfies
 \begin{equation*} N^{-\delta}\sum_{ u\in\Z} \left|\hat{f}(uN^{-\delta})\right|^2\big(2 |u|+1\big)
\le C N^{\delta} \int_0^\infty |\hat{f}(v)|^2 |v| dv \ . \end{equation*}
 Then (\ref{variance_split}) yields
 $$ N^{\delta-\alpha} \Cu^2_{K^N_p}[\Xi f_\delta]  = 2N^{-\delta}\sum_{u \in\Z} \left|\hat{f}(uN^{-\delta})\right|^2 \ N^{-\alpha}  \sum_{k>0} \sigma^2_k + \underset{N\to\infty}{O}\left(N^{\delta-\alpha} \|f\|_{H^{1/2}}^2  \right) \ .$$
 Moreover, according to formula (\ref{B^2}), a Riemann sum  approximation gives 
\begin{equation}  \label{Nvariance}
 \lim_{N\to\infty} N^{-\alpha} \sum_{k>0} \sigma^2_k  
=\int_\R \Psi(t \tau^{-1} )\big( 1-\Psi(t \tau^{-1})  \big)  dt  =  \tau \B^2_\Psi  \ .
\end{equation} 
By convention $\B^2_\Psi=1$ and we conclude that when $\delta<\alpha$, 
 \begin{equation*}  \lim_{N\to\infty} N^{\delta-\alpha} \Cu^2_{K^N_p}[\Xi f_\delta]  
 = 2 \tau  \int_{-\infty}^\infty |\hat{f}(\xi)|^2 d\xi  \ .  \end{equation*}
Since the variance of the random variable $\Xi f_\delta$ is diverging like $N^{\alpha-\delta}$ and its expected value is of order $N^{1-\delta}$ by equation (\ref{LS_mean_2}), the CLT follows from Soshnikov's theorem 1 in~\cite{Soshnikov_01}.  \qed \\


Note that, using the upper-bound (\ref{CUE_variance_estimate})  and the limit (\ref{Nvariance}), we get
$$  \sum_{k\in\Z} p_k(1- p_{k+u})  \le 2 |u|  + \underset{N\to\infty}{O}(N^\alpha) \ . $$
Hence, by formula (\ref{CUE_C2}),
\begin{equation} \ 
C^2_{K^N_p}[\Xi f_\delta] \le 2 N^{-\delta}\sum_{u\in\Z} \left|\hat{f}(uN^{-\delta})\right|^2  
\left\{  |u N^{-\delta}| + \underset{N\to\infty}{O}(N^{\alpha-\delta}) \right\} \ .
  \end{equation}
This implies that for any $f\in H^{1/2}_0$, the variance of the linear statistic $\Xi f_\delta$ remains bounded in the regime $\delta \ge \alpha$. Actually, if  (\ref{modified_CUE_2}) holds, we have in the regime $\delta>\alpha$,  
\begin{equation} \label{CUE_variance}
\lim_{N\to\infty}\Var_{K^N_p}[\Xi f_\delta] =  \|f\|_{H^{1/2}}^2  \ .
\end{equation}
This suggests that at any scale $\delta>\alpha$, we should have the same limit theorem for the modified CUEs as for the mesoscopic CUE and  sine process.  
We can prove formula (\ref{CUE_variance}) in the same way we obtained theorem~\ref{thm:CUE_Poisson} but the argument is already quite technical and becomes really sophisticated if we are interested in computing the limits of  the higher-order cumulants. 
A better approach consists in deducing the CLT from Soshnikov's theorem~\cite{Soshnikov_00a} by proving that the cumulants of a given linear statistics have the same limits regardless of the shape $\Psi$ of the modified CUE.

\begin{theorem} \label{thm:CUE_regime}
Consider a modified CUE with correlation kernel $(\ref{modified_CUE_1}-\ref{modified_CUE_2})$ and let $f\in H^{1}_0(\R)$.
For any scale $1 \ge  \delta >\alpha >0$, the linear statistics $\Xi f_\delta$ converges in distribution to a Gaussian random variable with variance $\|f\|_{H^{1/2}}^2$. \end{theorem}

\proof  
Let us decompose
\begin{equation} \label{coeff_splitting}
 p_k= \1_{|k|\le N} + \varepsilon_{k} \ . \end{equation}
By assumption, $ \varepsilon_{k}=  \Psi(\frac{|k|-N}{N^\alpha})-1$ when $| k|\le N$ and  $ \varepsilon_{k}= \Psi(\frac{|k|-N}{N^\alpha})$ when $| k|> N$.
We can write
$$\Cu^n_{K^N_p}[\Xi  f_\delta]= \Cu^n_{K^N}[\Xi  f_\delta] + \mathcal{E}_N^n(f,\delta,\alpha,\Psi) \ ,$$
where $\mathcal{E}_N^n(f,\delta,\alpha,\Psi)$ collects all the term which contains  at least one factor $\varepsilon_{k+u_1+\cdots+u_{\overline{ m}_i}}$ when we insert the decomposition (\ref{coeff_splitting}) into formula (\ref{CUE_cumulant_2}) and  expand the products $\displaystyle \prod_{i=0}^{l-1} p_{k+u_1+\cdots+u_{\overline{ m}_i}}$. Plainly, all other terms exactly add up to $C^n_{K^N}[\Xi  f_\delta] $.
Since $|p_k|,|\varepsilon_k|\le 1$ for all $k\in\Z$,  we get
\begin{equation*}\left| \mathcal{E}^N(f,\delta,\alpha,\Psi)\right|
\le C_n N^{-n\delta}  \sum_{u\in\Z^n_0} \left| \prod_i  \hat{f}(u_i N^{-\delta}) \right|  \sum_{k\ge 0}  | \varepsilon_{k} | \ ,\end{equation*}
where $\displaystyle C_n=  2\sum_{l=1}^n \frac{1}{l}\sum_{{\bf m}\in\N^l} {n \choose  {\bf m}}$.  Moreover, by the definition of $\varepsilon_k$, we have the estimates:
\begin{align*}\sum_{0\le k\le N} |\varepsilon_k|&= \sum_{-N\le k\le 0}1- \Psi(k N^{-\alpha})
\le  C N^\alpha \int_{-\infty}^0 1-\Psi(t)\ dt  \ , \\
\sum_{N< k}|\varepsilon_k| &=\sum_{0\le k} \Psi(k N^{-\alpha}) \le C N^\alpha \int_0^\infty \Psi(t)\ dt \ . \end{align*}
Both integrals are finite since $\Psi\in\F$ and there exists a positive constant $C'_n>0$  such that
\begin{equation*} \left| \mathcal{E}^N(f,\delta,\alpha,\Psi)\right| \le C'_n N^{\alpha-\delta} \int_{\R^n_0} \prod_i \left| \hat{f}(v_i) \right|  d^{n-1}v  \ . \end{equation*}

The assumption $f\in H^{1}_0(\R)$, guarantees that for any $n\in\N$, 
\begin{equation*}
\int_{\R^n_0} \prod_i \left| \hat{f}(v_i) \right|  d^{n-1}v  <\infty 
\end{equation*}
so that $\mathcal{E}^N(f,\delta,\alpha,\Psi) = O(N^{\alpha-\delta})$ as $N\to\infty$.
Therefore, all the cumulants $C^n_{K^N_p}[\Xi  f_\delta]$ and $C^n_{K^N}[\Xi  f_\delta]$ have the same limit and the CLT follows directly from Theorem~1 in~\cite{Soshnikov_00a}. \qed \\

\begin{remark} In the terminology of definition $\ref{def_equivalent_kernel}$, we have proved that  the rescaled correlation kernels 
$N^{-\delta}K_p^N(N^{-\delta}x,N^{-\delta}y)$ and $N^{-\delta}K^N(N^{-\delta}x,N^{-\delta}y)$  are asymptotically equivalent when the condition $\delta>\alpha$ is satisfied. We could also have deduced this fact from lemma $\ref{thm:lemma_perturbative}$ below by checking that  the CUE kernel $K^N$ given by $(\ref{CUE_kernel})$ satisfies the property $L^1B$ at any scale $\delta\in [0,1]$.
\end{remark}


%
%

\subsection{The critical regime} \label{sect:CUE_critical}

It remains to look at what happens at the critical scale $\delta=\alpha$. We have already seen that the variance remains bounded as $N\to \infty$.  We can compute its limit by applying a Riemann sum approximation to formula~(\ref{CUE_C2}). By symmetry,
\begin{equation*} 
\Cu^2_{K^N_p}[\Xi f_\alpha] \simeq 2 N^{-2\alpha}\sum_{u>0} \hat{f}(uN^{-\alpha})\hat{f}(-uN^{-\alpha}) \sum_{k>0} (p_{k-u} + p_{k+u})(1-p_k) \ . 
  \end{equation*}

and, since $p_{N+j} = \Psi(\frac{j}{\tau N^\alpha}) $ for any $j >-N$, we can check that for any $0<\alpha < 1$ and for any $\tau>0$,
\begin{align} 
\lim_{N\to\infty} \Var_{K^N_p}[\Xi f_\alpha] 
&\label{variance_5} = 2\int_0^\infty \hat{f}(u) \hat{f}(-u) \int_\R  \bigg(\Psi\left(\frac{t-u}{\tau}\right)+ \Psi\left(\frac{t+u}{\tau}\right) \bigg) \bigg( 1-\Psi\left(\frac{t}{\tau}\right) \bigg) dtdu \\
& \label{modified_CUE_variance} 
=   2\int_\R \big| \hat{f}(u) \big|^2 \int_\R \Psi\left(\frac{t+u}{\tau}\right)\left( 1-\Psi\left(\frac{t}{\tau}\right) \right) dtdu \ . 
  \end{align}


Because of some subtle cancellations, it is difficult to use formula (\ref{CUE_cumulant_2}) to compute  the limits of the higher-order cumulants by  Riemann sum approximations.  Another approach is to rewrite the correlation kernel of  the modified CUE before computing the cumulants. From definition \ref{modified_CUE}, a summation by parts yields
\begin{align} \notag
 K_p^N(x,y) &= \sum_{k=0}^\infty  (p^N_k -p^N_{k+1}) K^{k}(x,y) \\
 \label{CUE_sine_approximation_1}
 &= \frac{1}{\tau N^\alpha} \sum_{k=-N}^\infty  \Phi\left(\frac{k+\xi_k}{\tau N^\alpha}\right) \frac{\sin\big( (2N+2k+1) \pi (x-y)\big)}{\sin(\pi(x-y))} \ ,
\end{align}
where $\xi_k\in(0,1)$ are given by the mean-value theorem. We can use formula (\ref{CUE_sine_approximation_1}) to relate the kernel $K_p^N$ to the sine kernel and we will be able to use the ideas of~\cite{Soshnikov_00a} to compute the limits of the cumulants of linear statistics of the modified CUEs.

\begin{proposition}  \label{thm:CUE_kernel_sine}
At the critical scale $\delta=\alpha$, the modified CUE kernel $K_p^N$ and the kernel  $L^N_{\Psi,\eta}$ given by $(\ref{kernel_L})$ with $\eta(k)= (N+ k +\frac{1}{2})N^{-\alpha}$ are {\it asymptotically equivalent} in the sense of definition~$\ref{def_equivalent_kernel}$.
\end{proposition}

By  proposition~\ref{thm:GUE_kernel_sine}, a similar approximation holds for the modified GUEs. There is only a minor difference in the definition of the function $\eta$ and the limits of the cumulants of both models  will be computed in a common framework in section~\ref{sect:Critical}.
In order to prove proposition~\ref{thm:GUE_kernel_sine} and~\ref{thm:CUE_kernel_sine}, we need to provide a criterion to check whether two kernels are asymptotically equivalent. First, we need to introduce a new definition. A similar concept was introduced in~\cite{RV_07b} to control cumulants of some complex determinantal processes.

\begin{definition} \label{def_L^1B} A family of  kernels $(L^N)_{N>0}$ satisfies the property $L^1B$ if for any compact set $A \subseteq \X$,  there exists a sequence of functions $\Gamma_N:\X \to \R^+$ and $\nu>0$ such that and all  $(x,y)\in A^2$, 
$$ | L^N(x,y) | \le \Gamma_N(x-y) \ , $$
and
$$\| \Gamma_N \|_{L^1(\tilde{A})}  = \underset{N\to\infty}{O}\big(|\log N|^{\nu}\big) $$ 
where  $\tilde{A}=\{x= y-z : y, z\in A\}$. 
\end{definition}

\begin{lemma} \label{thm:lemma_perturbative} Two families of  kernels $(L^N)_{N>0}$ and $(K^N)_{N>0}$ are asymptotically equivalent if the family  $(L^N)_{N>0}$ has the property $L^1B$ and there exists $\kappa>0$ such that for any compact set $A\subseteq \X$, 
$$ \sup\left\{ | L^N(x,y)-K^N(x,y) | : (x,y) \in A^2\right\} =\underset{N\to\infty}{O}(N^{-\kappa}) \ .$$    \end{lemma}

\proof Let $\ell \in \N$ and $f_1,\dots, f_\ell \in C_0(\R) $. If we replace  $K^N= L^N+E^N$, we get 
\begin{equation*} \tr\left[ K^N f_1 \cdots K^N f_\ell\right] 
= \tr\left[ L^N f_1 \cdots L^N f_\ell\right] 
+ \sum_{J^k \in \{L^N,E^N\}} \tr[J^1 f_1\cdots J^{\ell} f_\ell ] \ .\end{equation*}
Note that all terms of the last sum contains at least one operator $E^N$. By assumption, we can suppose that all the test functions are supported in a compact set $A$ and  there exists two positive constants $C$ and $\kappa$ such that $  \sup\left\{ | E^N(x,y) | : (x,y) \in A^2\right\} \le C N^{-\kappa} $.
If we first  look at  a trace which contains a single operator $E^N$, by formula (\ref{cumulant_4}), we get the estimate
$$ \left| \tr[J^1 f_1 \cdots J^{\ell} f_\ell ]\right|
\le C N^{-\kappa} \prod_{k=2}^{\ell} \|f_k\|_\infty \int_{A^\ell} \left|  f_1(z_1) L^N(z_1,z_2)\cdots L^N(z_{\ell-1},z_\ell) \right| d\mu(x_1)\cdots d\mu(x_\ell) \ . $$
Since $L^N$ has the property $L^1B$, there exists $\Gamma_N:\X\to\R^+$ such that $|L^N(x_k,x_{k+1})|\le\Gamma_N(x_k-x_{k+1})$ and a change of variables yields
\begin{align*}\left| \tr[J^1 f_1\cdots J^{l} f_{\ell} ] \right|
&\le C N^{-\kappa}  \prod_{k=2}^{\ell}  \|f_k\|_\infty \|\Gamma_N\|_{L^1(\tilde{A})}^{\ell-1} \|f_1\|_{L^1}  \ .
\end{align*}
A similar argument shows that any trace which contains $j$ operators $E^N$ is bounded by $N^{-j\kappa}$ times a logarithmic correction coming from $\|\Gamma_N\|_{L^1(\tilde{A})}$.  Therefore, using the notation $\bar{O}$ introduced in section~\ref{sect:intro_5}, we get
\begin{equation} \label{cumulant_3}
 \tr\left[ K^N f_1 \cdots K^N f_\ell\right] 
= \tr\left[ L^N f_1 \cdots L^N f_\ell\right] +\underset{N\to\infty}{\bar{O}}\left( N^{-\kappa}\right) \ .
\end{equation}
 This completes the proof. \qed\\

{\it Proof of proposition~\ref{thm:CUE_kernel_sine}.}  
A Taylor expansion of the function $\sin(\pi(x-y) N^{-\alpha})$ in the denominator of   formula (\ref{CUE_sine_approximation_1}) shows that
\begin{equation}  \label{CUE_sine_approximation_2}
 N^{-\alpha}K_p^N(xN^{-\alpha},yN^{-\alpha})
 =  \frac{1}{\tau N^\alpha} \sum_{k=-N}^\infty  \Phi\left(\frac{k+\xi_k}{\tau N^\alpha}\right) \frac{\sin\big( 2\pi\eta(k) (x-y)\big)}{\pi(x-y)}  + \underset{N\to\infty}{O}\left( N^{-\alpha}\right) \ ,
\end{equation}
and the error term is uniform over any compact subset of $\R^2$.
Then, by lemma~\ref{thm:lemma_perturbative}, it is enough to prove that the RHS family of kernels denoted  $L^N_{\Psi,\eta}(x,y)$  satisfies the property $L^1B$. Note that the family of kernels $L^N_{\Psi,\eta}(x,y)$ is translation-invariant on $\R$ and we can choose $\Gamma_N(x)=|L^N_{\Psi,\eta}(0,x)|$.    
It is well-known that there exists a universal constant $C>0$ such that for any $s>0$ and $n>0$,
\begin{equation*}  \int_{-s}^s \left| \frac{\sin n (x-y)}{(x-y)} \right| dy \le C \log(sn) \ .   \end{equation*}
This implies that
 \begin{equation} \label{L_bound}
 \left\|  L^N_{\Psi,\eta} \right\|_{L^1[-s,s]}
 \le   \frac{C}{\tau N^\alpha} \sum_{k=-N}^\infty  \Phi\left(\frac{k+\xi_k}{\tau N^\alpha}\right)  \log( s\eta(k)) . \\
  \end{equation}

Since $\eta(k)=N^{1-\alpha} + \frac{k+1/2}{N^{\alpha}}$ and $\displaystyle \int_0^\infty \Phi(t) \log t\ dt<\infty $, we deduce from the estimate (\ref{L_bound}) that there is a constant $C'$ which only depends on the shape $\Psi$  such that $\left\|  L^N_{\Psi,\eta} \right\|_{L^1[-s,s]} \le C' \log(sN)$.    \qed\\

Proposition~\ref{thm:CUE_kernel_sine} implies that the determinantal processes with correlation kernels $K_p^N$ and  $L^N_{\Psi,\eta}$  have the same limit at the critical scale. By corollary~\ref{thm:weak_convergence}, this yields the following limit theorem for linear statistics of the critical modified CUEs. 

\begin{theorem} \label{thm:CUE_C^n}
Let $f\in H^{1}(\R)$ with compact support, $0<\alpha<1$ and $\Psi\in\F^*$. The linear statistic $\Xi f_\alpha$ of the determinantal process with correlation kernel $(\ref{modified_CUE_1}-\ref{modified_CUE_2})$ converges in distribution as $N\to \infty$ to a random variable $\Xi_{\Psi,\tau} f$ whose cumulants are given by
\begin{equation*} 
 \Cu^n\big[\Xi_{\Psi,\tau} f \big]=  
2\tau \B^n_\Psi \int_{\R} f(t)^n dt\ - 2\sum_{|{\bf m}|=n} \M({\bf m})\underset{\R^n_0} {\int du}  \underset{\R^n_<}{\int dx}\  \Re\left\{ \prod_{i=1}^n \hat{f}(u_i)\Phi(x_i) \right\} \G_\tau^{\bf m}(u,x)  \ .  
\end{equation*}
\end{theorem}

Up to a scaling, it is the same limit theorem as for the critical modified GUEs, theorem~\ref{thm:MNS_C^n}.

\section{Central Limit Theorems for the Modified GUEs} \label{sect:Modified_GUE}

We begin in section~\ref{sect:Variance} by proving some technical lemmas that are needed to compute the asymptotic variance of linear statistics of the modified GUEs. In particular we get formula~(\ref{variance_1}) for the variance at Poisson scales. 
In section~\ref{sect:MNS_CLT}, we prove theorem~\ref{thm:MNS_CLT} by comparing the rescaled correlation kernel of a modified GUE to the GUE kernel  using the perturbative method developed in section~\ref{sect:CUE_critical}. All these results are based on the asymptotics of the Hermite polynomials and the GUE kernel which are presented in section~\ref{sect:Asymptotic}.

\subsection{Proof of theorem~\ref{thm:Poisson_CLT}} \label{sect:Variance}

We start by proving a classical formula for the variance of linear statistics which is valid in a general context.

\begin{lemma}  \label{thm:variance}
Given a determinantal process with a correlation kernel $K$ of type $(\ref{kernel_0})$, 
for any test function $f\in C_0(\X)$, we have
\begin{equation*}
 \Var_K\left[ \Xi f \right]
= \sum_{k=0}^\infty \sigma_k^2 \int f(x)^2 |\varphi_k(x)|^2 \mu(dx)  
+\frac{1}{2} \iint (f(x)-f(y))^2 |K(x,y)|^2 \mu(dx) \mu(dy)  
\end{equation*}
where $\sigma_k^2=p_k(1-p_k)$.
\end{lemma}

\proof If we apply formula  (\ref{cumulant_2}) when $n=2$, 
\begin{align*}\Var_K\left[\Xi f \right]=& \int f(x)^2 K(x,x) \mu(dx)-\iint f(x)f(y)K(x,y)K(y,x)\mu(dx)\mu(dy) \\
=&\frac{1}{2}\iint (f(x)-f(y))^2K(x,y)K(y,x)\mu(dx)\mu(dy) \\
&-\iint f(x)^2K(x,y)K(y,x)\mu(dx)\mu(dy) +\int f(x)^2 K(x,x) \mu(dx) \ .
\end{align*}
Note that when the kernel $K$ is reproducing, the last two terms cancel. In general, since the function $\psi_k$ are orthonormal, we get
\begin{align*} &-\iint f(x)^2K(x,y)K(y,x)\mu(dx)\mu(dy) +\int f(x)^2 K(x,x) \mu(dx)\\
&= -\sum_{k,j} p_kp_j\int \varphi_j(y)\overline{\varphi_k}(y)\mu(dy)\int f(x)^2 \varphi_k(x)\overline{\varphi_j}(x)\mu(dx) +\sum_k p_k \int f(x)^2 \varphi_k(x)\overline{\varphi_k}(x)\mu(dx)\\
&=\sum_k \sigma_k^2 \int f(x)^2 |\varphi_k(x)|^2\mu(dx) \ . & \qed
\end{align*} 

For linear statistics of the modified GUEs, there is no counterpart of equation (\ref{CUE_cumulant_2}), but we can use lemma~\ref{thm:variance} to compute the asymptotic of the variance.
We call the {\bf reproducing variance} the quantity
\begin{equation}\label{variance_reproducing}
 V_0(f)=  \frac{1}{2} \iint |f(x)-f(y)|^2 \left|K_{\Psi,\alpha}^N(x,y)\right|^2 dx dy  \ .
\end{equation}
This definition comes from the fact that, if the correlation kernel $K$ is reproducing, then $\sigma^2_k=0$ for all $k \in \Z_+$ and $\Var_K\left[\Xi f \right]= V_0(f)$  for any linear statistic.  On the other hand, we call the {\bf Poisson variance} 
the quantity
\begin{equation}\label{variance_Poisson}
 V_\sigma(f)= \sum_{k=0}^\infty \sigma_k^2 \int f(x)^2 |\varphi_k(x)|^2 dx  \ .
 \end{equation}
This definition is motivated by the observation that considering a constant test function, we get $\Var_K[\#]= V_\sigma(1)$ and this quantity measures the {\it extra randomness} induced in the process from the fact that the correlation kernel is non-reproducing. 
In particular for a modified GUE we have
\begin{equation}\label{variance_coefficient} 
\sigma_k^2=\Psi\left(\frac{k-N}{\tau N^\alpha}\right) \left(1- \Psi\left(\frac{k-N}{\tau N^\alpha}\right)\right)  
\end{equation}
and this implies that
 \begin{equation}  \label{variance_modification_3}  
  \Var_{K^N_{\Psi,\alpha}}[\#] = \sum_{k=0}^\infty \sigma_k^2 \sim \tau N^\alpha \int_{-\infty}^{\infty} \Psi(t)\big( 1-\Psi(t) \big) dt  \ .
 \end{equation} 
This gives formula (\ref{variance_modification_1}), see also equation  (\ref{variance_modification_2}) in the introduction for a probabilistic interpretation.
We shall see that, except at the critical scale $\delta=\alpha$, only one component of the variance is asymptotically relevant. We begin by computing an asymptotic formula for the Poisson variance.

 \begin{lemma} \label{thm:variance_sigma}
 For any $0 < \alpha , \delta <1$ and for any function $f\in C_0(\R)$ we have 
 \begin{equation*}  V_\sigma(f_\delta) =  \frac{\tau}{2} N^{\alpha-\delta}  \int_\R f(x)^2 dx +\underset{N\to\infty}{o  }(N^{\alpha-\delta})  \ .   \end{equation*}
\end{lemma}

\begin{remark} It is not difficult to adapt the proof so that  lemma $\ref{thm:variance_sigma}$ holds for any function $f\in L^2(\R)$ which is uniformly continuous. In particular, by Morrey's inequality, this covers all test functions in the Sobolev space $H^1(\R)$. 
\end{remark}

\proof   Since the rescaled Hermite functions satisfy $\|\varphi_k\|_{L^2}=1$ for all $k \ge 0$ and we suppose that the test function $f$ is bounded, by formula  (\ref{variance_Poisson}), for any $0<\epsilon<1$,
 \begin{equation*} V_\sigma(f_\delta)
=  \sum_{|k|<N^{1-\epsilon}}  \sigma_{N+k}^2 \int f_\delta(x)^2 | \varphi_{N+k}(x ) |^2 dx 
+ O\bigg(\|f\|_\infty^2   \sum_{|k|>N^{1-\epsilon}}  \sigma_{N+k}^2  \bigg) \ .
 \end{equation*}
 The condition $\Psi\in\F$ guarantees that, if $\epsilon < 1-\alpha$,  the error term is converging to $0$ as $N\to\infty$. 
Actually, under the stronger assumption that $\Psi\in\F^*$, this term decays faster than any power of $N$ and it will be neglected in the following.
 Moreover the assumption that $f$ has compact support in conjunction with the condition $|k|<N^{1-\epsilon}$ implies that we can use the bulk asymptotic for the Hermite functions, formula (\ref{Hermite_asymptotic_bulk}).
Namely for any $x\in \supp(f)$,
\begin{equation}\label{Hermite_asymptotic_bulk_3} 
\varphi_{N+k}(x N^{-\delta})= \cos \left[ (N+k) \left(\frac{\pi}{2}- F(x_NN^{-\delta})\right)\right]+\underset{n\to\infty}{O}\left( N^{-\delta}\right) \ ,
\end{equation} 
where we set $x_N=x \frac{\pi}{2}\sqrt{\frac{N}{N+k}}$. Then
 \begin{equation*} V_\sigma(f_\delta)
=  N^{-\delta} \sum_{|k|<N^{1-\epsilon}}  \sigma_{N+k}^2 \left\{ \int f(x)^2 \left| \cos \left[ (N+k) \left(\frac{\pi}{2}- F(x_NN^{-\delta})\right)\right]\right|^2dx + O\left(N^{-\delta}\right) \right\}  \ .
 \end{equation*}
 Observe that according to formula (\ref{variance_modification_3}), we have  $\displaystyle  \sum_{|k|<N^{1-\epsilon}}  \sigma_{N+k}^2 = O( N^\alpha )$ and the previous estimate gives
 \begin{align*} V_\sigma(f_\delta)
&= \frac{N^{-\delta}}{2} \bigg\{ \Var_{K^N_{\Psi,\alpha}}[\#]\  \|f\|_{L^2}^2 \\
&\hspace{1.7cm} +  \sum_{|k|<N^{1-\epsilon}} (-1)^{N+k}  \sigma_{N+k}^2 \int f(x)^2 \cos \left[ 2(N+k) F(x_N N^{-\delta})\right]dx  + O(N^{\alpha-\delta} )  \bigg\}   \ .
 \end{align*}

The second term is a sum of oscillatory integrals and we will show that it converge to 0 as $N\to\infty$. 
Let us make the change of variable $z= N^\delta  F(x_N N^{-\delta})$. By definition \ref{F}, 
\begin{align*} & \int f(x)^2 \cos \left[ 2(N+k) F(x_N N^{-\delta})\right]dx \\
&\hspace{1cm}
=\frac{2}{\pi} \sqrt{1+\frac{k}{N}} \int  f^2\left(\frac{2N^\delta}{\pi} \sqrt{1+\frac{k}{N}} G( zN^{-\delta})\right) G'\left(z N^{-\delta}\right) \cos \left[2(N+k)N^{-\delta}z\right]  dz \ .\end{align*}

Since the function $f$ is uniformly continuous (we assume that $f$ has compact support in $[-\frac{L}{2},\frac{L}{2}]$), there exists a sequence $\epsilon_N \searrow 0$ such that  uniformly over all $|z|<L$ 
and all $|k|< N^{1-\epsilon}$,
$$ \left|  f^2\left(\frac{2N^\delta}{\pi} \sqrt{1+\frac{k}{N}} G(z N^{-\delta})\right)
-f^2\left(\frac{2G'(0)z}{\pi}\right) \right| \le  \epsilon_N \ .$$
Since $G'(0)=1/2$, it follows that for any $|k| < N^{1-\epsilon}$,
\begin{equation*} \int f(x)^2 \cos \left[ 2(N+k) F(x_N N^{-\delta})\right]dx
= \frac{1}{\pi} \int f^2\left(\frac{z}{\pi}\right)  \cos \left[2(N+k)N^{-\delta}z\right]  dz
+O\left(\epsilon_N\right) \ .
\end{equation*}

Since the sequence $(N-N^{1-\epsilon})N^{-\delta}\to\infty$ as $N\to\infty$, by the Riemann-Lebesgue lemma, we can also assume that 
\begin{equation*}\sup_{|k|\le N^{1-\epsilon}}
\left|\int f(x)^2 \cos \left[ 2(N+k) F(x_N N^{-\delta})\right]dx \right| \le  \epsilon_N \ .\end{equation*}
Going back to the Poisson variance, we have shown that
 \begin{equation*} V_\sigma(f_\delta)
=\frac{N^{-\delta}}{2}  \Var_{K^N_{\Psi,\alpha}}[\#]   \left\{ \|f\|_{L^2}^2
 +\underset{N\to\infty}{O}\left(\epsilon_N\right) \right\} \ .
\end{equation*}
The lemma follows after replacing $ \Var_{K^N_{\Psi,\alpha}}[\#] $ by formula (\ref{variance_modification}).   \qed\\
 
In order to prove (\ref{variance_1}) it remains to estimate the reproducing variance $V_0(f)$. Before proceeding we need to recall a few properties of the GUE correlation kernel (\ref{GUE_kernel}).  We refer to section~\ref{sect:Asymptotic} for additional details.  First,  note that according to the convention of definition \ref{modified_GUE},   the Christoffel-Darboux formula, \eqref{CD_kernel}, implies that for any $x, y \in \R$ and $M\ge 0$, 
\begin{equation}\label{CD}
 (x-y) K_0^M(x,y) =  \sqrt{M/N} \big(  \varphi_M(x)\varphi_{M-1}(y)-\varphi_{M-1}(x)\varphi _M(y) \big) 
\end{equation}

Moreover, the uniform bound for the Hermite functions,(\ref{Hermite_bound_uniform}), implies that there exists $C>0$ such that for any $n \ge 0$, 
 \begin{equation}   \label{Hermite_bound_uniform_2}
   \| \varphi_{n}\|_\infty \le CN^{1/4}  n^{-1/12} \ .
  \end{equation} 
 
 In particular, $\big| K_0^M(x,y) \big| \le C \sqrt{N} M $ and by formula \eqref{CD}, this gives us a bound for the GUE kernel. For any $x, y \in \R$, we have
\begin{equation*} 
\big| K_0^M(x,y) \big| \le C  \big( \sqrt{N} M \big)\wedge \frac{M^{1/3}}{|x-y|} \ .
\end{equation*}
 
The connection with the modified GUE kernel comes form a summation by parts:
\begin{align}
 K_{\Psi,\alpha}^N(x,y) 
 &\notag = \sum_{k=0}^\infty  (p_k -p_{k+1}) K_0^{k}(x,y)  \\
 & \label{kernel_approximation} = \frac{1}{\tau N^\alpha} \sum_{k=-N}^\infty  \Phi\left(\frac{k+\xi_k}{\tau N^\alpha}\right)K_0^{N+k}(x,y)  \ .
  \end{align}
where $\Phi=-\Psi'$ and $\xi_k \in (0,1)$ are given by the mean-value theorem.
If we further suppose that the shape $\Psi\in\F^*$, \eqref{class_F'},  there exists $\c>0$ so that for any $\Gamma>0$, 

\begin{equation}  \label{bound_1}
\bigg|  \frac{1}{\tau N^\alpha} \sum_{|k| >\Gamma N^\alpha } 
 \Phi\left(\frac{k+\xi_k}{\tau N^\alpha}\right) K_0^{N+k}(x,y)  \bigg| \le
C  e^{-\c\Gamma} \left(  N^{3/2} \wedge \frac{N^{1/3}}{|x-y|} \right) . 
\end{equation}
By formula \eqref{kernel_approximation}, this implies that for any sequence $\Gamma_N >0$, 
 \begin{equation} \label{kernel_approximation_1}
K_{\Psi,\alpha}^N(x,y)
 =  \frac{1}{\tau N^\alpha} \sum_{|k| \le \Gamma_N N^\alpha } 
 \Phi\left(\frac{k+\xi_k}{\tau N^\alpha}\right) K_0^{N+k}(x,y) + 
 \underset{N\to\infty}{O}\left(e^{-\c\Gamma_N} \left(  N^{3/2} \wedge \frac{N^{1/3}}{|x-y|} \right)\right) ,
\end{equation}
uniformly for all $x, y \in \R$. 

\begin{remark} \label{rk:exp_decay} 
The assumption $\Psi\in\F^*$ implies that, choosing  $\Gamma_N=(\log N)^2$, the error in formula $(\ref{kernel_approximation_1})$  decays faster than any power of $N$.
Analogous approximations hold for more general shapes, although with a worst error term which may not be good enough for all mesoscopic scales. Moreover, the condition $\Psi\in\F^*$ makes the proof almost trivial, otherwise we would need to take into account the speed of decay of $\Psi$ and to  produce more precise estimates.
 \end{remark}

  \begin{lemma} \label{thm:variance_zero}
 For any $0<\alpha<1$ and any scale $0<\delta \le 1$, there exists a constant $C>0$ such that for any function $f\in H^{1/2}_0(\R)$,   the reproducing variance satisfies for all sufficiently large $N$, 
 $$V_0(f_\delta) \le C \|f\|_{H^{1/2}}^2 \ .$$  
\end{lemma}
 
 \proof To simplify the notation, let us assume that the temperature $\tau=1$
and set the parameters $\xi_k=0$. 
 We will also let $x_N=xN^{-\delta}, y_N=yN^{-\delta}$ and fix  $L>0$ so that $\supp(f) \subset [-\frac{L}{2},\frac{L}{2}] $. 
 By formula (\ref{variance_reproducing}) and the approximation \eqref{kernel_approximation_1}, 
\begin{align} 
 V_0(f_\delta) 
 \notag &=  \frac{1}{2} \iint |f_\delta(x)-f_\delta(y)|^2  \left|K_{\Psi,\alpha}^N(x,y)\right|^2 dx dy \\
\label{approx_1} &\le  \iint |f_\delta(x)-f_\delta(y)|^2 \big| \tilde{K}_{\Psi,\alpha}^N(x,y)  \big|^2 dxdy
+  \underset{N\to\infty}{O}\left(\| f_\delta\|_{H^{1/2}}^2 N^{1/3}  e^{-\c\Gamma_N} \right) 
\end{align}
where 
\begin{equation}  \label{approx_2}
 \tilde{K}_{\Psi,\alpha}^N(x,y) = N^{-\alpha} \sum_{|k| \le \Gamma_N N^\alpha } \Phi\left(k N^{-\alpha}\right) K_0^{N+k}(x,y) \ . 
\end{equation}
 By formula \eqref{norm_1/2}, $\| f_\delta\|_{H^{1/2}}^2 = \| f\|_{H^{1/2}}^2$ and, if we let   
 $\Gamma_N=(\log N)^2$, the error in the previous estimate decays faster than any power of $N$ and it remains to show that  
 \begin{equation}   \label{approx_3}
 \iint |f_\delta(x)-f_\delta(y)|^2 \big| \tilde{K}_{\Psi,\alpha}^N(x,y) \big|^2 dxdy 
 =  \iint |f(x)-f(y)|^2 \big| N^{-\delta} \tilde{K}_{\Psi,\alpha}^N(x_N,y_N)\big|^2 dxdy 
 \le C \|f\|_{H^{1/2}}^2  .
 \end{equation}

 According to the sine-kernel approximation (\ref{sine_approximation}), if the density $N$ is sufficiently large compared to $L$, there exists a constant $C>0$ such that for all $|k|\le \Gamma_N N^\alpha$ $(\alpha<1)$ and  for all $x,y \in [-L,L]$,
 $$\left|N^{-\delta}K_0^{N+k}(x_N,y_N)\right| \le \frac{C}{ |x-y|}$$
  Since, for any $\Gamma>0$, 
 \begin{equation}  \label{sum_phi}
 \frac{1}{\tau N^\alpha}  \sum_{|k| \le \Gamma N^\alpha } 
 \Phi\left(\frac{k+\xi_k}{\tau N^\alpha}\right) \le  1 \ ,
 \end{equation}
this implies that  for all $x,y \in [-L,L]$
\begin{equation*}
\left|N^{-\delta}\tilde{K}_{\Psi,\alpha}^N(x_N,y_N)\right| 
\le \frac{C}{|x-y|} \  .
\end{equation*}
 
Hence, by a change of variables,
\begin{align}
 \iint |f_\delta(x)-f_\delta(y)|^2 \big| \tilde{K}_{\Psi,\alpha}^N(x,y) \big|^2 dxdy  
 &  \label{variance_decay_0}
 \le C  \iint_{[-L,L]^2}  \left|\frac{f(x)-f(y)}{x- y} \right|^2 dxdy \\
 &\notag +\iint_{\R^2\backslash[-L,L]^2}\left|f(x)-f(y)\right|^2 \left|  N^{-\delta} \tilde{K}_{\Psi,\alpha}^N(x_N,y_N)\right|^2 dxdy \ .
  \end{align} 
For any $L>0$, the integral (\ref{variance_decay_0}) is bounded by $\|f\|_{H^{1/2}}^2$ and to obtain the upper-bound \eqref{approx_3}, it suffices to show that  there exists a constant $C(f)\ge0$ which only depends on the test function $f$ such that 
\begin{equation}  \label{variance_decay}
 \iint_{\R^2\backslash[-L,L]^2}\left|f(x)-f(y)\right|^2 \left| N^{-\delta} \tilde{K}_{\Psi,\alpha}^N(x_N,y_N)\right|^2 dxdy\le \frac{C(f)}{L} \ .
\end{equation}
Thus, choosing the parameter $L$ sufficiently large, by formula (\ref{approx_1}), this implies that the variance $V_0(f_\delta) \le C\|f\|_{H^{1/2}}^2$. 
The rest of the proof is rather technical and is devoted to prove the estimate (\ref{variance_decay}). Since the function $f$ is supported in $ [-\frac{L}{2},\frac{L}{2}]$, by symmetry of the GUE kernel, we obtain
\begin{align} \label{variance_decay_1} 
& \iint_{\R^2\backslash[-L,L]^2}\left|f(x)-f(y)\right|^2 \left| N^{-\delta} \tilde{K}_{\Psi,\alpha}^N(x_N,y_N)\right|^2 dxdy \\
&\notag \le 4 N^{-2\alpha} \underset{\begin{subarray}{c}   |k| \le \Gamma N^\alpha \\  |j| \le \Gamma N^\alpha \end{subarray}}{\sum}  \Phi(kN^{-\alpha})  \Phi(jN^{-\alpha})  
\underset{\begin{subarray}{c}  |y| < L/2 \\ x>L \end{subarray}}{\iint} 
 \left| \frac{f(x)-f(y)}{x-y}\right|^2 \left|(x_N-y_N)^2 K_0^{N+k}(x_N,y_N) K_0^{N+j}(x_N,y_N) \right| dxdy \ .
\end{align}
 Since  $ \varphi_k(x)=\sqrt{\frac{\pi\sqrt{N}}{\sqrt{2}}}h_k\left( x\frac{\pi\sqrt{N}}{\sqrt{2}} \right) $, we deduce from the bulk asymptotic (\ref{Hermite_asymptotic_bulk}) with $\gamma=1/6$ that there exists a universal constant $C>0$ such that for all $|k| \le \Gamma_N N^\alpha$ and for all $|x|<\frac{2(\sqrt{N+k}-1)}{\pi \sqrt{N}}$,
  \begin{equation} \label{Hermite_asymptotic_uniform_bulk}
   | \varphi_{N+k}(x) |  \le \frac{C}{(4\frac{N+k}{N}- \pi^2 x^2)^{1/4}} . 
 \end{equation}
In particular, for any $|y|< L/2$,  we have $ | \varphi_{N+k}(y_N)| \le C$, and by formula (\ref{CD}),  
 \begin{align*} 
  \big| (x_N-y_N)^2   K_0^{N+k}(x_N,y_N) K_0^{N+j}(x_N,y_N)  \big| 
\le &C^2 \big(  | \varphi_{N+k}(x_N)\varphi_{N+j}(x_N) | +
 | \varphi_{N+k}(x_N)\varphi_{N+j-1}(x_N) |  \\
&+  | \varphi_{N+k-1}(x_N)\varphi_{N+j}(x_N)| + 
 | \varphi_{N+k-1}(x_N)\varphi_{N+j-1}(x_N) | \big) \ .
  \end{align*}
Let
  \begin{equation*}
 J_{k,j}= 
\underset{\begin{subarray}{c}  |y| < L/2 \\ x> L \end{subarray}}{\iint} 
\left| \frac{f(x)-f(y)}{x-y}\right|^2 | \varphi_{N+k}(x_N) \varphi_{N+j}(x_N)|\ dxdy \ .
\end{equation*}

By \eqref{variance_decay_1}, we see that there exists $C>0$ such that
\begin{align} 
&\notag \iint_{\R^2\backslash[-L,L]^2}\left|f(x)-f(y)\right|^2 \left| N^{-\delta} \tilde{K}_{\Psi,\alpha}^N(x_N,y_N)\right|^2 dxdy  \\
 &\label{variance_decay_4}\hspace{2cm} \le C N^{-2\alpha} \underset{\begin{subarray}{c}   |k| \le \Gamma N^\alpha \\  |j| \le \Gamma N^\alpha \end{subarray}}{\sum}  \Phi(kN^{-\alpha})  \Phi(jN^{-\alpha})  
\big\{ J_{k,j}+  J_{k,j-1} +  J_{k-1,j} +  J_{k-1,j-1}  \big\} .
\end{align}

So, in order to prove the estimate \eqref{variance_decay}, by \eqref{sum_phi}, it remains to show that  
$ J_{k,j} \le C(f) /L $ for all $|k|,|j| \le \Gamma N^\alpha$. To do so, we shall use the asymptotics  from section~\ref{sect:Asymptotic}. First of all, since $\supp(f) \subset  [-\frac{L}{2},\frac{L}{2}] $, we have for any $|x|>L$,
 \begin{equation*}  
  \left| \frac{f(x)-f(y)}{x-y}\right| \le \1_{y\in \supp(f)} 
  \frac{2\|f\|_\infty}{|x-L/2|} \le  \1_{y\in \supp(f)}
 \frac{4\|f\|_\infty}{|x|} 
  \end{equation*} 
 and, if $C(f)= 4 \|f\|_\infty^2|\supp f| $, we get
\begin{equation*}
\int_{-L/2}^{L/2}  \left| \frac{f(x)-f(y)}{x-y}\right|^2  dy \le \frac{C(f)}{|x|^2}  \ .
\end{equation*}
 Hence, by a change of variables,
 \begin{equation} \label{variance_decay_2}
   J_{k,j} \le C(f) N^{-\delta}  \int_{LN^{-\delta}}^\infty | \varphi_{N+k}(x) \varphi_{N+j}(x)| \frac{dx}{x^2} \ .
\end{equation}

Suppose that $ j\ge k$ and let $a_\pm = \frac{2(\sqrt{N+k}\pm 1)}{\pi \sqrt{N}}$. We split the integral:
 \begin{equation} \label{variance_decay_3}
\int_{LN^{-\delta}}^\infty | \varphi_{N+k}(x) \varphi_{N+j}(x)| \frac{dx}{x^2} = 
\bigg\{ \int_{LN^{-\delta}}^{a_-} + \int_{a_-}^{a_+} + \int_{a_+}^\infty \bigg\} | \varphi_{N+k}(x) \varphi_{N+j}(x)| \frac{dx}{x^2} \ .
\end{equation}

Using the upper-bound \eqref{Hermite_asymptotic_uniform_bulk}, the first integral gives
\begin{equation*}
\int_{LN^{-\delta}}^{a_-}   | \varphi_{N+k}(x) \varphi_{N+j}(x)| \frac{dx}{x^2} 
\le C^2 \int_{LN^{-\delta}}^{a_-} \frac{dx}{x^2 (4\frac{N+k}{N}- \pi^2 x^2)^{1/2}}  \ .
\end{equation*}
Since $a_-^2 \le  4\frac{N+k}{\pi N}$, when $N$ is sufficiently large, we obtain
\begin{equation} \label{approx_4}
\int_{LN^{-\delta}}^{a_-}   | \varphi_{N+k}(x) \varphi_{N+j}(x)| \frac{dx}{x^2} 
\le \frac{C N^{\delta}}{L} \ . 
\end{equation}
Using the uniform bound  (\ref{Hermite_bound_uniform_2}), since $a_+ -a_- = \frac{4}{\pi} N^{-1/2}$, the contribution from the edge gives
\begin{equation}\label{approx_5}
\int_{a_-}^{a_+} | \varphi_{N+k}(x) \varphi_{N+j}(x)| \frac{dx}{x^2} \le C N^{-1/6} \ .
\end{equation}
Finally, we have
$$
\int_{a_+}^\infty | \varphi_{N+k}(x) \varphi_{N+j}(x)| \frac{dx}{x^2}
 \le C  N^{1/6} \int_{a_+}^\infty | \varphi_{N+k}(x) | \frac{dx}{x^2} 
$$
and the estimate \eqref{Hermite_asymptotic_exterior} implies that for all $x > a_+$, 
\begin{equation*}
\big| \varphi_{N+k}(x)  \big| \le C N^{1/4} e^{- \frac{ 2\sqrt{2} \pi  N^{3/4}}{ 3\sqrt{N+k}}}  .  
\end{equation*}
Thus, we obtain for all $|k|,|j| \le \Gamma N^\alpha$, 
\begin{equation} \label{approx_6}
\int_{a_+}^\infty | \varphi_{N+k}(x) \varphi_{N+j}(x)| \frac{dx}{x^2}
 \le C  N^{5/12} e^{- 2\sqrt{2} N^{1/4} } . 
\end{equation}

If we put together  the estimates (\ref{variance_decay_3} - \ref{approx_6}), we have proved that, when $N$ is sufficiently large, 
  \begin{equation} \label{variance_decay_3}
\int_{LN^{-\delta}}^\infty | \varphi_{N+k}(x) \varphi_{N+j}(x)| \frac{dx}{x^2} \le  \frac{C N^{\delta}}{L} \ . 
\end{equation}

Hence, it follows from formula (\ref{variance_decay_2}) that for any $|j|, |k| \le \Gamma_N N^\alpha$, the integral  $J_{k,j}  \le C(f)/ L $.  By  (\ref{variance_decay_4}) and \eqref{sum_phi}, we conclude that the estimate \eqref{variance_decay} holds and this completes the proof. \qed\\ 

We are now ready to finish the proof of formula  (\ref{variance_1}) and hence of theorem~\ref{thm:Poisson_CLT}. It
 follows immediately from  lemmas~\ref{thm:variance},~\ref{thm:variance_sigma} and~\ref{thm:variance_zero}
that in the regime $\delta<\alpha$, for any test function $f\in H^{1/2}_0 \cap L^\infty(\R)$,
\begin{equation*} 
 \Var_{K^N_{\Psi,\alpha}}[\Xi f_\delta] = \frac{\tau}{2} N^{\alpha-\delta} \int_\R f(x)^2 dx +\underset{N\to\infty}{o  }(N^{\alpha-\delta}) \ .
 \end{equation*}

The same argument shows that, in the regime $\delta\ge \alpha$,
\begin{equation}\label{variance_estimate}
  \Var_{K^N_{\Psi,\alpha}}[\Xi f_\delta] \le C\left( \|f\|_{L^2}^2 +  \|f\|_{H^{1/2}}^2 \right)
\end{equation}
 At the GUE scales ($\delta>\alpha$), the limit of the variance is given by theorem~\ref{thm:MNS_CLT} which is proved in the next section. At the critical scale $\delta=\alpha$, by lemma~\ref{thm:variance_sigma}, the Poisson variance $ V_\sigma(f_\alpha)$ converges to $\frac{\tau}{2}  \|f\|^2_{L^2}$ and the limit of the reproducing variance $ V_0(f_\alpha)$ is computed in appendix \ref{A:variance} by a Riemann sum approximation.

  \subsection{Proof of theorem~\ref{thm:MNS_CLT}} \label{sect:MNS_CLT} 

 \begin{theorem}  \label{thm:Szego_2} 
Let $\Xi$ be the GUE eigenvalue process with correlation kernel $K_0^N$ given by $(\ref{GUE_kernel})$. For   any $0<\delta<1$ and any function $f\in H^{1/2}_0\cap L^\infty(\R)$, as the number of eigenvalues $N\to\infty$, 
$$\tilde\Xi f_\delta\ \Rightarrow\ \No\left(0, \|f\|^2_{H^{1/2}}\right) \ . $$
\end{theorem}


Theorem~\ref{thm:Szego_2} was first established in~\cite{BK_99a, BK_99b} for the resolvent function $ x\mapsto (x-z)^{-1}$ where $\Im z>0$. A general proof was given only recently in~\cite{FKS_13}, 
their argument exploits  a nice connection between the characteristic polynomial of a GUE matrix and a log-correlated Gaussian process. In ~\cite{BEYY_14}, a generalization of theorem~\ref{thm:Szego_2} is given for  Gaussian $\beta$-Ensembles.
Yet another generalization to certain classes of Orthogonal Polynomial ensembles is made in~\cite{ BD_14, L_15a}. In particular, the proof of theorem~\ref{thm:Szego_2}  in~\cite{L_15a} is based on the cumulant computations presented in section~\ref{sect:Critical} and the sine-kernel asymptotics of theorem~\ref{thm:meso_sine_kernel}. 
%
%
%
We now turn to the approximation of the modified GUE correlation kernels at the so-called  GUE scales ($\delta>\alpha$). By definition \ref{def_equivalent_kernel}, proposition~\ref{thm:GUE_sine} below combined with theorem~\ref{thm:Szego_2} implies the central limit theorem~\ref{thm:MNS_CLT}.

\begin{proposition} \label{thm:GUE_sine}
For any shape $\Psi\in\F^*$, the modified GUE correlation kernel $K_{\Psi,\alpha}^N$ and the GUE kernel $K_0^N$ are asymptotically equivalent at any scale $\delta>\alpha$, 
$$ N^{-\delta}K_{\Psi,\alpha}^N( N^{-\delta}x,N^{-\delta}y) \cong N^{-\delta}K_0^N( N^{-\delta}x,N^{-\delta}y) \ . $$
\end{proposition}

\proof To simplify the notation, let us assume that temperature $\tau=1$
and set the parameters $\xi_k=0$. 
The condition $\Psi\in\F^*$ implies that for any $\Gamma>0$, 
\begin{equation*}
 \sum_{k>\Gamma N^{\alpha}} \Psi\left(kN^{-\alpha}\right)  
+ \sum_{k<- \Gamma N^{\alpha}}\left( 1- \Psi\left(k N^{-\alpha}\right)\right)  \le 
C  N^{\alpha} e^{-\Gamma} \ .
\end{equation*}

So that, if we let $\Gamma_N=(\log N)^2$, these sums decay faster than any power of $N$ and combined with the uniform bound (\ref{Hermite_bound_uniform_2}), this implies that
\begin{align*}
&K_{\Psi,\alpha}^{N}(x,y)\simeq \sum_{k=0}^{N-1} \varphi_{k}(x) \varphi_{k}(y) 
+\sum_{k=0}^{\Gamma_N N^\alpha} \Psi\left(k N^{-\alpha}\right) \varphi_{N+k}(x) \varphi_{N+k}(y)\\
&\hspace{4cm} -\sum_{k=1}^{\Gamma_N N^\alpha} \left(1- \Psi\left(-k N^{-\alpha}\right)\right) \varphi_{N-k}(x) \varphi_{N-k}(y)
\end{align*}
with a uniform error of order $N^{1/2+\alpha} e^{-(\log N)^2}$ as $N\to\infty$. 
Moreover, for any $L>0$, the bulk estimate (\ref{Hermite_asymptotic_uniform_bulk}) implies that  that all $x,y\in [-L,L]$, 
\begin{equation*}
\sum_{|k|\le \Gamma_N N^\alpha} \left| \varphi_{N+ k}(N^{-\delta}x) \varphi_{N+ k}(N^{-\delta}y) \right|
\le 2C^2 \Gamma_N N^\alpha \ .
\end{equation*}
Since the function  $\Psi\in [0,1] $, using the notation $\bar{O}$ introduced in section~\ref{sect:intro_5}, this yields
 \begin{equation*}
N^{-\delta}K_{\Psi,\alpha}^N( N^{-\delta}x,N^{-\delta}y)= N^{-\delta} \sum_{k=0}^{N-1}
\varphi_{k}( N^{-\delta}x) \varphi_k(N^{-\delta}y) + \underset{N\to\infty}{\bar O}\left( N^{\alpha-\delta} \right) 
\end{equation*}
uniformly for all $x,y\in [-L,L]$.
The sum in the RHS is the rescaled GUE correlation kernel.  By lemma~\ref{thm:lemma_perturbative}, to prove that the kernels $K_{\Psi,\alpha}^N$ and $K_0^N$ are asymptotically equivalent, it remains to show the latter satisfies the property $L^1B$. 
Taking $M=N$ in the  approximation (\ref{sine_approximation}) implies that there exists a positive constant $C_L$  such that for any $0<\delta \le 1$,

\begin{equation*}
 \left| N^{-\delta}K_0^N(xN^{-\delta},yN^{-\delta}) \right|
 \le \Gamma_N(x-y) :=  C_L \begin{cases} N &\text{if } |x-y| < \frac{\log N}{N} \\
\frac{1}{|x-y|} &\text{if } |x-y| \ge \frac{\log N}{N}
 \end{cases}
 \ ,
\end{equation*}
and we  immediately check that 
$\displaystyle \int_{-L}^L \Gamma_N(z)dz  \le 4C_L \log N $. \qed

\section{Cumulants of the Critical models} \label{sect:Critical}

In this section, we prove theorem~\ref{thm:MNS_C^n} and~\ref{thm:CUE_C^n}, then we analyze the random processes $\Xi_{\Psi,\tau}$ which arise from  the critical modified ensembles. Actually we will not investigate directly the modified ensembles but  the processes with kernel $L^N_{\Psi,\eta}$ given by (\ref{kernel_L}). By propositions ~\ref{thm:GUE_kernel_sine} and~\ref{thm:CUE_kernel_sine}, there are two choices of the function $\eta$  which correspond to the modified GUEs and CUEs respectively. However, our analysis works as long as $\eta$ satisfies the conditions (\ref{eta_condition_1} - \ref{eta_condition_2}) below.
In section~\ref{sect:MNS_Critical}, we show that $L^N_{\Psi,\eta}$ is the correlation kernel of a determinantal process and we prove proposition~\ref{thm:GUE_kernel_sine}. The convergence of smooth linear statistics of these processes is established  in section~\ref{sect:C^n}; see corollary~\ref{thm:weak_convergence}.
The main result in section~\ref{sect:Poisson} is that, for all $\Psi\in\F$ such that $\Psi \neq \psi $ and for all  $\tau>0$, the random variables $\Xi_{\Psi,\tau} f$ are not Gaussian; see propositions~\ref{thm:Gaussian} and~\ref{thm:MNS_property}. 
In section~\ref{sect:C3C4}, we show that, despite the special property of the MNS shape $\psi$, the MNS ensemble at the critical scale also converges to a random process which is not Gaussian.

\subsection{Asymptotically equivalent kernels for the critical modified GUEs} \label{sect:MNS_Critical}

\begin{lemma}\label{thm:kernel_L}
Let  $N,\tau,\Gamma>0$, $\alpha \in (0,1)$, $\Psi\in\F$ and $\eta$ be a non-decreasing function. The kernel $L^N_{\Psi,\eta}$ given by $(\ref{kernel_L})$ defines a translation-invariant determinantal process on $\R$. 
\end{lemma}

\proof  The fundamental property of the kernel $L^N_{\Psi,\eta}$ is that it is translation-invariant. Hence we can define its Fourier transform
\begin{equation} \label{kernel_Fourier_2} 
\hat{L}^N_{\Psi,\eta}(v)=   \frac{1}{\tau N^\alpha} \sum_{|k|\le \Gamma N^\alpha}  \Phi\left(\frac{k +\xi_k}{\tau N^\alpha}\right) \1_{ [-\eta(k),\eta(k)]}(v) \ .
\end{equation}
Plainly, the function $\hat{L}^N_{\Phi,\eta} \in L^1(\R)$ and by (\ref{kernel_L}),
\begin{equation} \label{kernel_Fourier_1}
 L^N_{\Psi,\eta}(x,y) = \int_\R  \hat{L}^N_{\Psi,\eta}(v)  e^{i 2\pi v (x-y) } dv \ .
 \end{equation}
 This definition comes from the article~\cite{Soshnikov_01} and it is also established that for any translation-invariant kernel $L$, the condition  $0 \le \hat{L}\le 1$, guarantees that it defines a determinantal point process.  
The parameters $\xi_k$ have been chosen so that 
$\frac{1}{\tau N^\alpha} \Phi\left(\frac{k +\xi_k}{\tau N^\alpha}\right)=\Psi\left(\frac{k}{\tau N^\alpha}\right)- \Psi\left(\frac{k+1}{\tau N^\alpha}\right)$
 and it follows that for any $v\in\R$, 
\begin{equation*}
\hat{L}^N_{\Psi,\eta}(v) \le \Psi(-\Gamma \tau^{-1}) \le 1 \ .
\end{equation*} 
Moreover, since $\Phi \ge 0$ by assumption, $\hat{L}^N_{\Psi,\eta} \ge 0$ and we conclude that  $L^N_{\Psi,\eta}$ is the correlation kernel of some determinantal process. \qed\\

{\it Proof of proposition~\ref{thm:GUE_kernel_sine}.} Let $1/3<\alpha<1$, $\Psi \in\F^*$, and $\Gamma_N=(\log N)^2$. We also assume that $\tau=1$.
We combine the approximation (\ref{kernel_approximation_1}) of the modified GUE kernel $K_{\Psi,\alpha}^N$ with the asymptotic formula of theorem~\ref{thm:micro_sine_kernel} with $x_0=0$. 
Namely, taking $\delta=\alpha$ in formula (\ref{sine_approximation_2}), we get
\begin{align*}
N^{-\alpha}K_{\Psi,\alpha}^N(xN^{-\alpha},yN^{-\alpha})
=N^{-\alpha} \sum_{|k| \le \Gamma_N N^\alpha } 
 \Phi\left(\frac{k+\xi_k}{ N^\alpha}\right) 
 \frac{\sin\big[ \pi N^{1-\alpha}\sqrt{1+k/N} (x-y)\big]}{\pi(x-y)} +  O \left( N^{1-3\alpha} \right) \ .
\end{align*}

The estimates (\ref{L_bound}) shows that, with $\eta(k)=\frac{1}{2}N^{1-\alpha}\sqrt{1+k/N}$, the kernel $L^N_{\Psi,\eta}$ given by (\ref{kernel_L}) has the property $L^1B$ and  it follows from lemma~\ref{thm:lemma_perturbative} that 
\begin{equation*}
N^{-\alpha}K_{\Psi,\alpha}^N(xN^{-\alpha},yN^{-\alpha}) \cong  L^N_{\Psi,\eta}(x,y) \ .
\end{equation*}
\qed

Proposition~\ref{thm:GUE_kernel_sine} and ~\ref{thm:CUE_kernel_sine} imply  that each of the modified ensembles have an asymptotically equivalent kernel at the critical scale of the form  $L^N_{\Psi,\eta}$  with 
\begin{equation} \label{eta} 
\begin{array}{l} i)\ \  \eta(k)=(N +k)N^{-\alpha} \text{ for the modified CUEs} \ ,  \\
ii)\ \ \eta(k)=\frac{1}{2}N^{1-\alpha}\sqrt{1+k/N} \text{ for the modified GUEs} \ .
\end{array}\end{equation}

In the sequel, we will compute the limits of the cumulants  for any determinantal process with kernel $L^N_{\Psi,\eta}$ which satisfies the following conditions.
The function $\eta$ is non-decreasing and it satisfies uniformly for all $|k| \le \Gamma_N N^\alpha $,
\begin{equation} \label{eta_condition_1}
 \eta(k) =N^\nu +\beta k N^{-\alpha} +\underset{N\to\infty}{O}(N^{-\epsilon}) \ ,
\end{equation}
where $\nu,\beta, \epsilon>0$  such that $N^\nu \gg \Gamma_N $ and 
\begin{equation} \label{eta_condition_2}
 \lim_{N\to\infty} N^\nu \max\{\Psi(\Gamma_N), 1-\Psi(-\Gamma_N)\} =0 \ .
\end{equation}
In particular, for the modified GUEs (resp$.$ CUEs), the asymptotics (\ref{eta_condition_1}) holds with $\nu=1-\alpha$ and $\beta=1/4$ (resp$.$ $\beta=1$)
and, if the shape $\Psi \in\F^*$, the condition (\ref{eta_condition_2}) holds for any $\alpha \in (0,1)$ with 
$\Gamma_N=(\log N)^2$.

\subsection{Proof of theorem~\ref{thm:MNS_C^n}} \label{sect:C^n}


Given the expression (\ref{kernel_Fourier_1}) of the correlation kernel $L_{\Psi,\eta}^N$, we can repeat the proof of lemma ~\ref{thm:CUE_cumulants} replacing sums by integrals and we get the following formula

\begin{equation} \label{C^n_1}  
C^n_{L^N_{\Psi,\eta}}[\Xi f]= 
\int_{\R^n_0} \prod_i \hat{f}(u_i)\sum_{|{\bf m}|=n} \M({\bf m})  \left( \int_\R\  \prod_{i=1}^{\ell({\bf m})}\hat{L}^N_{\Phi,\eta}(v+u_1+\cdots+u_{\overline{ m}_i}) dv \right)d^{n-1}u \ ,
\end{equation}
where the sum is over all compositions ${\bf m}$ of the number $n$,  (\ref{M}). Combining this formula with (\ref{kernel_Fourier_2}), we get an  expression for the cumulants  that is appropriate to pass to the limit as $N\to\infty$. In this section, to simplify the notation, we will assume that $\xi_k=0$ and (unless stated otherwise) all sums run over $|k_i| \le \Gamma_N N^\alpha$.
Define for any composition ${\bf m}$ of $n\ge 2$, the function 
\begin{equation} \label{H_1}
\h_{\bf  m}(u,k) = \int_\R\  \prod_{i=1}^{\ell({\bf m})}  \1_{|v+u_1+\cdots+u_{\overline{ m}_i}|\le \eta(k_i)}\ dv \ .
\end{equation}

\begin{lemma} \label{thm:C^n_Fourier}
For  any function $f\in C_0(\R)$, the cumulants of linear statistics of the determinantal process with correlation kernel $(\ref{kernel_L})$ are given by
\begin{equation} \label{C^n_2}  
C^n_{L^N_{\Psi,\eta}}[\Xi f] \simeq \left(\frac{1}{\tau N^\alpha}\right)^n \sum_{k_1\le \cdots\le k_n}\ \prod_{i=1}^{n} \Phi\left(\frac{k_i}{\tau N^\alpha}\right) \int_{\R^n_0} \prod_i \hat{f}(u_i)  \sum_{|\bf m|=n} \M({\bf m}) \sum_{\sigma \in \Sy(n)} \h_{\bf  m}(u, \sigma k)\ d^{n-1}u \ ,
\end{equation}
where the function $\h_{\bf m}$ is given by $(\ref{H_1})$.
\end{lemma}

\proof To simplify the notation, let us assume that $\tau=1$. 
By formula (\ref{kernel_Fourier_2}), for any composition ${\bf m}$ of $n$ of length $\ell$ and any $v \in \R^\ell$, we have
\begin{equation}  \label{H_4}
\prod_{i=1}^{\ell} \hat{L}^N_{\Phi,\eta}(v_i)  =  N^{-\alpha\ell} \sum_{k_1,\cdots, k_\ell}\ \prod_{i=1}^{\ell}  \Phi\left(\frac{k_i}{ N^\alpha}\right) \1_{|v_i|\le \eta(k_i)} \ .
\end{equation}

If we let  $\epsilon_N^n=0$ and for any $1\le \ell <n$,
\begin{equation*} \label{H_error}
\epsilon^\ell_N = 1- N^{- \alpha (n-\ell)} \sum_{k_{\ell+1},\cdots, k_n}\ \prod_{i=\ell+1}^{n} \Phi\left(\frac{k_i}{ N^\alpha}\right) \ ,
\end{equation*}
we get 
\begin{equation*}
\prod_{i=1}^{\ell} \hat{L}^N_{\Phi,\eta}(v_i) \big\{1- \epsilon_N^\ell\big\}=N^{-\alpha n} \sum_{k_1,\cdots, k_n}\ \prod_{i=1}^{n} \Phi\left(\frac{k_i}{ N^\alpha}\right) \prod_{i=1}^{\ell} \1_{|v_i|\le \eta(k_i)} \ .
\end{equation*}
By (\ref{H_1}), this implies that
\begin{equation}  \label{H_5}
 \big\{ 1-\epsilon_N^\ell  \big\} \int_\R \prod_{i=1}^{\ell} \hat{L}^N_{\Phi,\eta}(v+u_1+\cdots+u_{\overline{ m}_i}) dv
 = N^{-\alpha n}  \sum_{k_1,\cdots,k_n}\ \prod_{i=1}^{n} \Phi\left(\frac{k_i}{ N^\alpha}\right) \h_{\bf  m}(u,k)  \ . 
\end{equation}
Observe that by (\ref{H_4}), since 
$\displaystyle N^{-\alpha}\sum_{\kappa} \Phi(\kappa N^{-\alpha}) \le 1$ and we assume that the condition  (\ref{eta_condition_1}) holds, there exists $C>0$ such that
$$ \int_\R \prod_{i=1}^{\ell} \hat{L}^N_{\Phi,\eta}(v_i) dv_1  \le  C  N^\nu \ . $$  
Moreover, by definition $0 \le \epsilon^\ell_N \le \Psi(\Gamma)$, 
 so that according to the condition (\ref{eta_condition_2}),
$$  \lim_{N\to\infty} \epsilon_N^\ell \int_\R \prod_{i=1}^{\ell} \hat{L}^N_{\Phi,\eta}(v+u_1+\cdots+u_{\overline{ m}_i}) dv = 0 \ . $$

Thus, by formula (\ref{H_5}),
\begin{equation*}
\int_\R \prod_{i=1}^{\ell} \hat{L}^N_{\Phi,\eta}(v+u_1+\cdots+u_{\overline{ m}_i}) dv \simeq N^{-\alpha n}  \sum_{k_1\le\cdots\le k_n}\ \prod_{i=1}^{n} \Phi\left(\frac{k_i}{ N^\alpha}\right)  \sum_{\sigma \in \Sy(n)} \h_{\bf  m}(u, \sigma k)   \ .
\end{equation*}
We conclude by using formula (\ref{C^n_1}). \qed\\

We can use the notation  (\ref{Lambda} - \ref{N_3}) to compute the function $\h_{\bf m}$ given by (\ref{H_1}). The computation is analogous to the proof of formula (\ref{H_sine}).  
  In the sequel, we will always use the conventions $\ell=\ell({\bf m})$, $\Lambda^{\bf m}_{i, s}=\Lambda^{\bf m}_{i, s}(u)$  and $\s=\s_{\ell({\bf m})}(\sigma) $.

\begin{lemma} \label{thm:H_3}  Let  $k\in\Z^n_\le$ and $u\in\R^n_0$. For any $\sigma \in \Sy(n)$ and any composition ${\bf m}$ of $n \ge 2$,
\begin{equation*}
\h_{\bf m}(u,\sigma k)
= \left[ 2 \eta(k_{\sigma(\s)})-\max_{i \le \ell}\left\{\Lambda^{\bf m}_{i,\s} -\eta(k_{\sigma(i)}) +  \eta(k_{\sigma(\s)})
\right\}
-\max_{i \le \ell}\left\{-\Lambda^{\bf m}_{i,\s} -\eta(k_{\sigma(i)}) +  \eta(k_{\sigma(\s)})
\right\} \right]^+ \ .  
\end{equation*}
\end{lemma}

\proof Let
$\displaystyle v_i=  \sum_{j=1}^{{\overline{ m}_i}} u_j $. 
The change of variable $ w=v- v_{\s}$ in (\ref{H_1}) gives
\begin{equation}
 \h_{\bf  m}(u,\sigma k) 
 =  \int_\R\  \prod_{i=1}^{\ell}   \1_{|w+\Lambda_{i,\s}^{\bf m}|\le \eta(k_{\sigma(i)})} \ dw  
\label{H_2}= \bigg| \bigcap_{i=1}^{\ell}\bigg\{ w : | w+\Lambda_{i,\s}^{\bf m}|\le \eta(k_{\sigma(i)})\bigg\}  \bigg|  \ . 
\end{equation} 

By definition $\min_{i\le\ell} \{\sigma(i) \}=\sigma(\s)$, and since the function $\eta$ is non-decreasing, for any $k\in \Z^n_\le$,
 \begin{equation} \label{min_eta}
\eta(k_{\sigma(\s)}) =\min_{i\le \ell}\eta(k_{\sigma(i)}) \ .\end{equation}

Then,  since $\Lambda_{\s, \s}^{\bf m}=0$ by (\ref{Lambda}), we get
\begin{align*} &\bigcap_{i=1}^{\ell}\bigg\{|w+\Lambda_{i,\s}^{\bf m}|\le \eta(k_{\sigma(i)})\bigg\}  \\
&= \left[-\eta(k_{\sigma(\s)}) + \max_{i\le l}\big\{-\Lambda_{i,\s}^{\bf m} -\eta(k_{\sigma(i)}) + \eta(k_{\sigma(\s)}) \big\},\ \eta(k_{\sigma(\s)}) - \max_{i\le l}\big\{\Lambda_{i,\s}^{\bf m} -\eta(k_{\sigma(i)}) + \eta(k_{\sigma(\s)}) \big\}\right]   \ .
\end{align*}
This interval is non-empty if the condition
$$ 2 \eta(k_{\sigma(\s)}) > \max_{i\le \ell}\big\{\Lambda_{i,\s}^{\bf m} -\eta(k_{\sigma(i)}) + \eta(k_{\sigma(\s)}) \big\}
+ \max_{i\le \ell}\big\{-\Lambda_{i,\s}^{\bf m} -\eta(k_{\sigma(i)}) + \eta(k_{\sigma(\s)}) \big\}$$
is satisfied, in which case the lemma follows from equation (\ref{H_2}). \qed\\

We are now ready to  prove our main result, i.e.~to compute the limits of the cumulants of linear statistics of the modified ensembles by applying a Riemann sum approximation to formula (\ref{C^n_2}). The argument is similar to the proof of Lemma 2 in~\cite{Soshnikov_00a} but it is more complicated. To keep the proof as transparent as possible, it relies on three lemmas which will be proved afterwards.

\begin{theorem} \label{thm:C^n}  Assume that the conditions $(\ref{eta_condition_1})$ and $(\ref{eta_condition_2})$ are satisfied and let $f\in H^1_0(\R)$.
For any $n\ge 2$, 
\begin{equation*} \label{C^n} 
\lim_{N\to\infty}\Cu^n_{L^N_{\Psi,\eta}}[\Xi f] =
2\beta\tau \B^n_\Psi \int_{\R} f(t)^n dt\
 -2\underset{\R^n_0} {\int du}  \underset{\R^n_<}{\int dx}\ \Re\left\{ \prod_{i=1}^n  \hat{f}(u_i)\Phi(x_i) \right\}\sum_{|{\bf m}|=n} \M({\bf m}) \G_{\beta\tau}^{\bf m}(u,x)\ ,  
\end{equation*}
where the function $\G^{\bf m}_{\tau}(u,x)$ and constant $\B^n_\Psi$ are defined  by formulae $(\ref{G})$ and  $(\ref{B})$.
\end{theorem}

\proof
Throughout the proof, we will use the familiar inequality (\ref{sum_phi}) without any reference.
Let $|u|_1=|u_1|+\cdots +|u_n|$ and
\begin{align}\label{Upsilon_0}
& \Upsilon^n_N(u)
=  \left(\frac{1}{\tau N^\alpha}\right)^n \sum_{k_1\le\cdots\le k_n}\ \prod_{i=1}^{n} \Phi\left(\frac{k_i}{\tau N^\alpha}\right) \times  \\
&\hspace{3cm}\notag\sum_{|{\bf m}|=n}  \M({\bf m})\sum_{\sigma \in \Sy(n)}   \left( \eta(k_{\sigma(\s)}) -\max_{i \le \ell}\left\{\Lambda^{\bf m}_{i,\s} -\eta(k_{\sigma(i)}) + \eta(k_{\sigma(\s)}) 
\right\} \right) \ .
\end{align}

If the parameter $N$ is sufficiently large, we claim that  for any $u \in \R^n_0$, any $\sigma \in \Sy(n)$, and for all
 $k\in\Z^n_\le$ such that $|k|_\infty< \Gamma_N N^{\alpha}$,
\begin{align} \notag 
&\left| \h_{\bf m}(u,\sigma k)- 2 \eta(k_{\sigma(\s)})+\max_{i \le \ell}\left\{\Lambda^{\bf m}_{i,\s} -\eta(k_{\sigma(i)}) +  \eta(k_{\sigma(\s)})\right\}
+\max_{i \le \ell}\left\{-\Lambda^{\bf m}_{i,\s} -\eta(k_{\sigma(i)}) +  \eta(k_{\sigma(\s)})\right\}\right|  \\
&\hspace{1cm}\label{Sos_bound}
 \le \begin{cases}  0 &\text{if } |u|_1 \le \frac{N^\nu}{2}   \\
		        18 |u|_1  &\text{else }
\end{cases} \ . \end{align}
First note that, since $ \Lambda^{\bf m}_{\s, \s}=0$, 

\begin{equation*} 
0\le \max_{i \le \ell}\left\{ \pm \Lambda^{\bf m}_{i,\s} -\eta(k_{\sigma(i)}) + \eta(k_{\sigma(\s)}) 
\right\} \ .
\end{equation*}
 Moreover, by (\ref{min_eta}),
\begin{equation*} 
 \max_{i \le \ell}\left\{ \pm \Lambda^{\bf m}_{i,\s} -\eta(k_{\sigma(i)}) + \eta(k_{\sigma(\s)}) 
\right\} \le  \max_{i \le \ell}\left\{ \pm \Lambda^{\bf m}_{i,\s} \right\} \ .
\end{equation*}
By formula (\ref{Lambda}), for any composition ${\bf m}$ of $n$, we have $|  \Lambda^{\bf m}_{i,\s} | \le |u|_1$ for all $i \le \ell$. Hence, we conclude that
\begin{equation} \label{Lambda_bound}
0\le \max_{i \le \ell}\left\{ \pm \Lambda^{\bf m}_{i,\s} -\eta(k_{\sigma(i)}) + \eta(k_{\sigma(\s)}) 
\right\}
  \le |u|_1 \ .
  \end{equation}
When the parameter $N$ is  large, condition (\ref{eta_condition_1}) implies that for any $|\kappa|\le \Gamma N^\alpha$, 
\begin{equation}\label{eta_condition_3}
\frac{N^\nu}{2} < \eta(\kappa) < 2 N^\nu \ .
\end{equation}
Thus, if we also suppose that  $|u|_1\le \frac{N^\nu}{2}$, by  (\ref{Lambda_bound}), 
$$  \eta(k_{\sigma(\s)})> \max_{i \le \ell}\left\{ \pm \Lambda^{\bf m}_{i,\s} -\eta(k_{\sigma(i)}) +  \eta(k_{\sigma(\s)})\right\}  \ .$$
By lemma~\ref{thm:H_3}, we conclude that when $|u|_1\le \frac{N^\nu}{2}$,
$$ \h_{\bf m}(u,\sigma k)- 2 \eta(k_{\sigma(\s)})+\max_{i \le \ell}\left\{\Lambda^{\bf m}_{i,\s} -\eta(k_{\sigma(i)}) +  \eta(k_{\sigma(\s)})\right\}+\max_{i \le \ell}\left\{-\Lambda^{\bf m}_{i,\s} -\eta(k_{\sigma(i)}) +  \eta(k_{\sigma(\s)})\right\}=0 \ . $$

For the second estimate, we observe that the estimate  (\ref{eta_condition_3}) implies that
$$ 0 \le \h_{\bf m}(u,\sigma k) \le   2 \eta(k_{\sigma(\s)}) \le 4 N^\nu \ . $$ 
Then, by the triangle inequality and (\ref{Lambda_bound}), the l.h.s$.$ of  (\ref{Sos_bound}) is bounded by $8N^\nu + 2|u|_1$. 
Thus, we have also proved (\ref{Sos_bound}) in the case when $|u|_1>\frac{N^{\nu}}{2}$.  If we combine this estimate with formula (\ref{C^n_2}) for the cumulants of the random variable $\Xi f$, there exists a positive constant $C_n$ which only depends on $n$ such that if the parameter $N$ is sufficiently large,
\begin{equation} \label{C^n_3}
\left|\Cu^n_{L^N_{\Psi,\eta}}[\Xi f] - \int_{\R^n_0} \prod_i \hat{f}(u_i)\big\{ \Upsilon^n_N(u) +  \Upsilon^n_N(-u) \big\} d^{n-1} u \right|
\le C_n  \int_{\R^n_0} \1_{\left\{|u|_1 >\frac{N^{\nu}}{2}\right\}}  |u|_1  \prod_i \left| \hat{f}(u_i) \right| d^{n-1}u \ ,
\end{equation}
where the function  $\Upsilon^n_N(u)$ is given by (\ref{Upsilon_0}). 
Taking $f_1=\cdots=f_n=f$ in lemma~\ref{thm:integrability_condition} below implies that the RHS of (\ref{C^n_3}) converges to 0 as $N\to\infty$. Thus, the limits of the cumulants are given by
\begin{align} \notag
\lim_{N\to\infty}\Cu^n_{L^N_{\Psi,\eta}}[\Xi f] 
&= \lim_{N\to\infty} \int_{\R^n_0} \prod_i \hat{f}(u_i) \big\{ \Upsilon^n_N(u) +  \Upsilon^n_N(-u) \big\} d^{n-1} u \\   
&\label{C^n_4} 
=2\lim_{N\to\infty} \int_{\R^n_0} \Re\left\{ \prod_i \hat{f}(u_i) \right\} \Upsilon^n_N(u)\ d^{n-1} u \ .
\end{align}

The next step is to compute the limit of $\Upsilon^n_N(u)$ as $N\to\infty$; see equation (\ref{Upsilon_3}).
Observe that, according to condition (\ref{eta_condition_1}) and since the $\max$ function is Lipschitz continuous, we get uniformly for all  $u\in \R^n_0$,
\begin{align} \label{Upsilon_1}
 \Upsilon^n_N(u)
=\left(\frac{1}{\tau N^\alpha}\right)^n & \sum_{k_1\le\cdots\le k_n}\ \prod_{i=1}^{n} \Phi\left(\frac{k_i}{\tau N^\alpha}\right) \sum_{|{\bf m}|=n}  \M({\bf m})  \times \\
&\notag  \sum_{\sigma \in \Sy(n)} 
 \left(  \eta(k_1) + \beta \frac{k_{\sigma(\s)}- k_1}{N^\alpha}   -\max_{i \le \ell}\left\{\Lambda^{\bf m}_{i,\s} -\beta\frac{k_{\sigma(i)} - k_{\sigma(\s)}}{N^\alpha} \right\}  \right) 
 + \underset{N\to\infty}{O}\left(N^{-\epsilon}\right) \ .
\end{align}
By lemma~\ref{thm:combinatorics_1}  below, $ \sum  \M({\bf m}) =0$ and we can remove the two terms $\eta(k_1)$ and $k_1N^{-\alpha}$ from formula (\ref{Upsilon_1})  since they do not depend on ${\bf m}$,
Hence we have proved that
\begin{equation} \label{Upsilon_2}
\Upsilon^n_N(u)=\left(\frac{1}{\tau N^\alpha}\right)^n  \sum_{k_1\le\cdots\le k_n}\ \prod_{i=1}^{n} \Phi\left(\frac{k_i}{\tau N^\alpha}\right) \sum_{|{\bf m}|=n}  \M({\bf m}) 
 \left(\beta  \sum_{\sigma \in \Sy(n)} \frac{k_{\sigma(\s)}}{N^\alpha}  - \G^{\bf m}_{\beta N^{-\alpha}}\left(u,k\right)  \right)  + \underset{N\to\infty}{O}\left(N^{-\epsilon}\right) \ ,
\end{equation}
where $\G^{\bf m}_{\beta N^{-\alpha}}$ is given by  (\ref{G}). 
Then, a Riemann sum approximation implies that 
\begin{equation*} 
\lim_{N\to\infty} \Upsilon^n_N(u) :=  \Upsilon^n_\infty(u) 
 = \int_{\R^n_<}  \prod_{i=1}^{n} \Phi(z_i /\tau)  \sum_{|{\bf m}|=n}  \M({\bf m}) 
 \left( \beta  \sum_{\sigma \in \Sy(n)} z_{\sigma(\s)} -\G^{\bf m}_{\beta }(u,z) \right) d^nz \ .
\end{equation*}
The first sum is independent of the Fourier variable $u\in\R^n_0$ and it can be computed explicitly; see lemma~\ref{thm:combinatorics_2} below.  Furthermore, making  the change of variables $x_i = z_i \tau $, we obtain
\begin{equation} \label{Upsilon_3}
 \Upsilon^n_\infty(u) = \beta\tau \B^n_\Psi- \int_{\R^n_<}  \prod_{i=1}^{n} \Phi(x_i) 
 \sum_{|{\bf m}|=n}  \M({\bf m}) \G^{\bf m}_{\beta\tau}(u,x) d^nx \ .
\end{equation}
Now, we can deduce the limits of the cumulants of the random variable $\Xi f$ from equations (\ref{C^n_4}) and (\ref{Upsilon_3}). 
By (\ref{G}) and the estimate $|  \Lambda^{\bf m}_{i,\s} | \le |u|_1$, we get
\begin{equation} \label{G_bound}
\sup\big\{  \G^{\bf m}_\tau(x,u)  : x\in\R^n_>, \tau>0 \big\} \le n! |u|_1 \ .
\end{equation}
Moreover, since $|k_{\sigma(\s)}| \le |k_1|+|k_n|$ for any $k \in \Z_\le^n$, by formula (\ref{Upsilon_2}), there exists a constant $C$ which only depends on $n$ such that for any $N>0$,   
\begin{equation} \label{Upsilon_4}
\left| \Upsilon^n_N(u) \right| 
\le C \int_{\R^n_<}  \prod_{i=1}^{n} \Phi(x_i /\tau) \big(1+ |x_1|+ |x_n| + |u|_1 \big) d^nx \ .
\end{equation}
The assumption $\Psi\in\F$ guarantees that the RHS of (\ref{Upsilon_4}) is finite. From lemma~\ref{thm:integrability_condition} below,  (\ref{Upsilon_3}) and the dominated convergence theorem, we conclude that
\begin{equation*} \label{C^n_5}   
\lim_{N\to\infty}\Cu^n_{L^N_{\Psi,\eta}}[\Xi f]
=2  \int_{\R^n_0} \Re\left\{ \prod_i \hat{f}(u_i) \right\}\left(\beta\tau \B^n_\Psi- \int_{\R^n_<}  \prod_{i=1}^{n} \Phi(x_i) 
 \sum_{|{\bf m}|=n}  \M({\bf m}) \G^{\bf m}_{\beta\tau}(u,x) d^nx \right) d^{n-1}u  \ .
\end{equation*}

The final observation is that the integral over $\R^n_0$ can be  written as a convolution. Namely a change of variables gives
\begin{equation*}
\int_{\R^n_0} \prod_{i=1}^n \hat{f}(u_i)   d^{n-1}u
= \int_{\R^{n-1}}  \hat{f}(v_1) \prod_{i=2}^{n-1} \hat{f}(v_i-v_{i-1}) \hat{f}(-v_{n-1})   d^{n-1}v  
= \underbrace{\hat{f} *\cdots * \hat{f}}_n(0) \ .
\end{equation*}
If we replace $ \hat{f} *\cdots * \hat{f}= \widehat{f^n}$ and evaluate at 0, we get
$\displaystyle \int_{\R^n_0} \prod_{i=1}^n \hat{f}(u_i)   d^{n-1}u =  \int_{\R} f(t)^n dt $
and the proof of theorem~\ref{thm:C^n} is completed. \qed \\

Now we prove the lemmas that we used to get theorem~\ref{thm:C^n}. 
The first lemma is classical, it was already used  in~\cite{Soshnikov_00a}, as well as in the context of other invariant ensembles, \cite{AHM_11, RV_07a, RV_07b}.

\begin{lemma} \label{thm:combinatorics_1} For any $n \ge 1$,
 \begin{equation*}
\sum_{|{\bf m}|=n} \M({\bf m})= \begin{cases} 1 &\text{if } n=1 \\ 0 &\text{if } n\ge 2\end{cases} \ ,
\hspace{.6cm} \text{and} \hspace{.6cm}
 \sum_{|{\bf m}|=n} \left| \M({\bf m}) \right| \le n! 2^{n-1}  \ .
 \end{equation*}
 \end{lemma}

We have assumed that our test function $f$ has compact support since the original problem is to study mesoscopic linear statistics. However this assumption is not necessary to prove theorem~\ref{thm:C^n}. 
We shall certainly require that $f\in L^1(\R)$ and, according to the estimate  (\ref{C^n_3}), the  regularity condition needed to prove  theorem~\ref{thm:C^n}  is that, for any $n \ge 2$, the integral 
\begin{equation*}
\int_{\R^n_0}\left|\hat{f}(u_1)\cdots \hat{f}(u_n)\right| \big( 1+ |u|_1 \big) d^{n-1}u <\infty \ .
\end{equation*}
A sufficient condition is provided by the next lemma since by assumption:
 $\|\hat{f}\|_\infty \le \| f \|_{L^1}<\infty$.

\begin{lemma} \label{thm:integrability_condition} 
Let $n \ge 2$. For any functions $f_1,\dots, f_n \in H^1(\R)$,
\begin{equation}  \label{f_estimate}
 \int_{\R^n_0}\left| \hat{f_1}(u_1)\cdots \hat{f_n}(u_n) \right| \big(1+ |u_1|+\cdots +|u_n| \big) d^{n-1}u \le n2^n \prod_{j=1}^n \left( \|\hat{f_j}\|_\infty+ \|f_j\|_{H^1}   \right)  \ .    \end{equation}
\end{lemma}

\proof 
By the Cauchy-Schwartz inequality,
\begin{align} \int_\R \bigg| \hat{f_1}(u_1)\hat{f_2}\bigg(- \sum_{j\neq 2 }u_j\bigg) \bigg| \big(1+ |u_1| \big) du_1
& \notag
 \le 2 \left(2\|\hat{f_1}\|_{\infty}\|\hat{f_2}\|_{\infty} + \int_{|u|>1} \bigg| \hat{f_1}(u)\hat{f_2}\bigg(-u- \sum_{j> 2 }u_j\bigg) \bigg| |u| du     \right) \\
& \label{f_estimate_1}
 \le 2 \left(2\|\hat{f_1}\|_{\infty}\|\hat{f_2}\|_{\infty} +\|f_1\|_{H^1}\|f_2\|_{L^2} \right)  
\end{align}
A similar argument shows that for any $f \in L^1 \cap H^1(\R)$, 
\begin{equation} \label{f_estimate_2}
 \|\hat{f}\|_{L^1} \vee \|\hat{f}\|_{L^2} \le 2 \left(  \|\hat{f}\|_{\infty}+ \|f\|_{H^1}  \right) \ .
 \end{equation}
Hence, it follows from (\ref{f_estimate_1}) that
$$  \int_\R \bigg| \hat{f_1}(u_1)\hat{f_2}\bigg(- \sum_{j\neq 2 }u_j\bigg) \bigg| \big(1+ |u_1| \big) du_1 
\le 4\big(\|\hat{f_1}\|_\infty +\|f_1\|_{H^1} \big)(\|\hat{f_2}\|_\infty +\|f_2\|_{H^1} \big) \ ,$$ 
and, if we combine this estimate with (\ref{f_estimate_2}),
\begin{align*} 
\int_{\R^n_0}\left| \hat{f_1}(u_1)\cdots \hat{f_n}(u_n) \right| \big(1+ |u_1| \big) d^{n-1}u   
&\le 4\big(\|\hat{f_1}\|_\infty +\|f_1\|_{H^1} \big)(\|\hat{f_2}\|_\infty +\|f_2\|_{H^1} \big) \prod_{j>2} \|\hat{f_j}\|_{L^1} \\
& \le 2^n \prod_{j=1}^n \left( \|\hat{f_j}\|_\infty+ \|f_j\|_{H^1}   \right) \ .
\end{align*}
The upper-bound (\ref{f_estimate}) follows by symmetry. \qed \\



The next lemma shows how the shape-dependent constant $ \B^n_\Psi$ defined by (\ref{B}) arises in (\ref{Upsilon_3}).

\begin{lemma} \label{thm:combinatorics_2} For any $n \ge 1$, 
\begin{align} \label{B_1}
\B^n_\Psi & = \int_{\R^n_<}  \prod_{i=1}^{n} \Phi(x_i) \sum_{\sigma\in\Sy(n)}  \sum_{|{\bf m}|=n}  \M({\bf m}) x_{\sigma(\s)} d^nx \\
\notag &=  \sum_{k=0}^{n-1} b_k^n \int_\R z \Phi(z)\Psi(z)^k\big(1-\Psi(z)\big)^{n-1-k}  dz 
\end{align}
where,  according to formula  $(\ref{N_3})$,
 $\displaystyle x_{\sigma(\s)}=\min_{i\le \ell({\bf m})}\{ x_{\sigma(i)}\}$ for any $x\in \R^n_<$, 
  and the coefficients $b_k^n$ are given by formula $(\ref{b})$.  
\end{lemma}

\proof 
Let $\mathbb{P}_n$ be  the uniform measure over the group $\Sy(n)$, so that we can view  $\sigma(\s_l)= \min_{i \le l}\{\sigma(i)\}$ as a random variable. Then, we can rewrite equation (\ref{B_1}) as 
 \begin{equation} \label{B_2}
 \B^n_\Psi = n! \int_{\R^n_<}  \prod_{i=1}^{n} \Phi(x_i)  \sum_{|{\bf m}|=n}  \M({\bf m}) \E{n}{x_{\sigma(\s_\ell)}} d^nx \ .
\end{equation}

We claim that for any $l=1,\dots,n$  and for any $k=0,\dots,n-1$,
\begin{equation} \label{Combinatorics_4}\P{\min_{i\le l}\sigma(i) \ge n-k}={ k+1\choose l}\frac{l! (n-l)!}{n!} \ .\end{equation}
To see this, observe that if $k+1<l$, since there are only $k$ elements in $\{1,\dots,n\}$ which are greater than $n-k$, one of the $l$ first elements of $\sigma$ has to be less than $n-k$ and therefore the probability in question is 0.\\
On the other hand, if $l\le k+1$, then $n-k$ is smaller than the minimum of the $l$ first entries of $\sigma$ if and only if these entries are drawn from the set $\{n-k,\dots, n\}$. Since the order of these entries and that of the $(n-l)$ last entries is irrelevant, the number of such permutations is $\displaystyle {k+1 \choose l } l!(n-l)!$.\\
Hence, by definition and equation (\ref{Combinatorics_4}), 
 the distribution of  $\sigma(\s_l)$ is given by
\begin{equation*}\mathbb{P}_n\left[\sigma(\s_l) = n-k\right]= {n \choose l}^{-1} { k\choose l-1} \ .\end{equation*}
Then, by definition of $\M$, (\ref{M}),
\begin{equation} \label{Combinatorics_5}
\sum_{|{\bf m}|=n}   \M({\bf m}) \E{n}{x_{\sigma(\s_\ell)}} 
 = \sum_{l=1}^n \frac{(-1)^{l+1}}{l}{n \choose l}^{-1}
  \sum_{\begin{subarray}{c} |{\bf m}|=n \\ \ell({\bf m})=l\end{subarray}}  {n \choose {\bf m}} 
\sum_{k=0}^{n-1} { k\choose l-1} x_{n-k} \ .
\end{equation}
If we integrate successively over $x_1,\cdots, x_{n-k-1}$ and over $x_n,\cdots, x_{n-k+1}$, and use the relationship $\Phi=-\Psi'$, we find that for any $k=0,\dots,n-1$,
\begin{equation}  \label{Combinatorics_3}
\int_{\R^n_<} \prod_{i=1}^n \Phi(x_i) x_{n-k} d^nx 
= \frac{1}{k! (n-k-1)!} \int_\R \Psi(x)^k (1-\Psi(x))^{n-k-1}  \Phi(x) x  dx \ . 
\end{equation} 
 Then, if we combine formulae (\ref{B_2}), (\ref{Combinatorics_5}) and (\ref{Combinatorics_3}), we get
 \begin{equation*}
 \B^n_\Psi = \sum_{l=1}^n (-1)^{l+1}\frac{n!}{l!}{n \choose l}^{-1} 
   \sum_{\begin{subarray}{c} |{\bf m}|=n \\ \ell({\bf m})=l\end{subarray}}  {n \choose {\bf m}} 
\sum_{k=0}^{n-1} \frac{1}{(k-l+1)!(n-k-1)!} \int_\R \Psi(x)^k (1-\Psi(x))^{n-k-1}  \Phi(x) x  dx \ .
 \end{equation*}
We see that we can simplify $\frac{n!}{l!}{n \choose l}^{-1}$ from the previous formula and exchange the sums over $l$ and $k$. In the end, we obtain
\begin{equation*}
\B^n_\Psi = \sum_{k=0}^{n-1} \bigg(  \sum_{l=1}^n (-1)^{l+1}{n- l \choose k+1-l}
  \sum_{\begin{subarray}{c} |{\bf m}|=n \\ \ell({\bf m})=l\end{subarray}}  {n \choose {\bf m}}  \bigg) \int_\R \Psi(x)^k (1-\Psi(x))^{n-k-1}  \Phi(x) x  dx \ .
 \end{equation*}
 If we define the array $b^n_k$ according to (\ref{b}), the lemma is proved. \qed\\

The {\it Cumulant problem} is generally not discussed directly in the literature, so we provide a simple criterion which guarantees uniqueness of the law of a random variable given its cumulants.  

\begin{lemma} \label{thm:Cumulant_Problem}
Given a sequence of random variables $X_N$ whose Laplace transform is well-defined and
such that for any $n\ge1$, the cumulant $\Cu^n[X_N] \to \Cu^n_\infty$ as $N\to\infty$. If there exists constants $c,v>0$ such that 
$$  | \Cu^n_\infty | \le c n! v^n \,, $$
then there exists a random variable $X_\infty$ whose cumulants satisfy  $\Cu^n[X_\infty]= \Cu^n_\infty$ and the sequence $X_N\Rightarrow X$.
\end{lemma}

The condition of lemma~\ref{thm:Cumulant_Problem}  is very natural and its proof follows from a straightforward repetition of the argument that is used when dealing with the Hamburger moment problem (see e.g$.$ section 3.3.3 in~\cite{Durrett_10}).
Next, we use this criterion to deduce from theorem ~\ref{thm:C^n} the weak convergence of smooth linear statistics $\Xi f$ for any determinantal process with correlation kernel $L^N_{\Psi,\eta}$.

\begin{definition} \label{C_Poisson} 
In the sequel, the quantity $\displaystyle 2\tau \B^n_\Psi \int_{\R} f(t)^n dt $ will be called the 
{\bf Poisson component} of the $n^{\text{th}}$ cumulant 
 and we will use the decomposition
$\displaystyle \lim_{N\to\infty}\Cu^n_{L^N_{\Psi,\eta}}[\Xi f] = 2\tau \B^n_\Psi \int_{\R} f(t)^n dt +  \mathfrak G ^n_{\Psi,\tau}[f] $
where 
\begin{equation} \label{GG}
 \mathfrak G ^n_{\Psi,\tau}[f] =- 2\int_{\R^n_0} du\ \Re\left\{\prod_{i=1}^n \hat{f}(u_i) \right\} \int_{\R^n_<}\prod_{i=1}^n \Phi(x_i) \sum_{|{\bf m}|=n} \M({\bf m}) \G_\tau^{\bf m}(u,x)\ dx  \ . \\
\end{equation}
\end{definition}
This name is motivated by linear statistics of the Poisson point process whose cumulants are equal to~$\displaystyle \tau' \int_{\R} f(t)^n dt  $ where $\tau'$ is the intensity.

\begin{corollary} \label{thm:weak_convergence}
 Consider the determinantal process with correlation kernel $L^N_{\Psi,\eta}$ and let $f\in H^1_0(\R)$. If the conditions $(\ref{eta_condition_1}-\ref{eta_condition_2})$  hold, then the random variable $\Xi f$ converges in distribution as $N\to\infty$ to a random variable $\Xi_{\Psi,\tau'} f$ where $\tau'=\beta \tau$ and whose cumulants are  given by   
\begin{equation} \label{C_split}
\Cu^n\big[\Xi_{\Psi,\tau'} f\big]  =  2\tau' \B^n_\Psi \int_{\R} f(t)^n dt + \mathfrak G ^n_{\Psi,\tau'} [f] \ .
\end{equation} 
\end{corollary}

\proof  We can estimate the growth of $\mathfrak G ^n_{\Psi,\tau'}[f]$ and the Poisson component separately.
We start by giving an upper-bound for the constant $\B^n_\Psi$ .
By formula (\ref{B_2}),
 \begin{equation*} 
\B^n_\Psi=  n!\int_{\R^n_<} \prod_{i=1}^n \Phi(x_i) \sum_{|{\bf m}|=n}   \M({\bf m}) \E{n}{x_{\sigma(\s_\ell)}} d^nx \ .
\end{equation*} 
Obviously for any $x\in\R^n_<$,  $\E{n}{x_{\sigma(\s_\ell)}} \le  x_n$ and if we use formula (\ref{Combinatorics_3}),
\begin{equation*} \int_{\R^n_<}  \prod_{i=1}^{n} \Phi(x_i) \E{n}{x_{\s_\ell}} d^nx  
\le  \int_{\R^n_<}  \prod_{i=1}^{n} \Phi(x_i) x_n d^nx
 = \frac{1}{(n-1)!} \int_\R (1-\Psi(x))^{n-1} \Phi(x) x dx \ .
\end{equation*}
Moreover, since $0\le \Phi=-\Psi'$ and $0\le \Psi\le 1$, we have for any $n\ge1$, 
\begin{equation*}  \int_\R (1-\Psi(x))^{n-1} \Phi(x) x dx 
\le \int_0^\infty \Phi(x)x dx
= \int_0^\infty \Psi(x) dx \ .
\end{equation*}
On the other hand, if we use that  $x_1 \le \E{n}{x_{\sigma(\s_\ell})} $ and apply the same method, we can show that
\begin{equation*}  \int_{\R^n_<}  \prod_{i=1}^{n} \Phi(x_i) \E{n}{x_{\s_\ell}} d^nx  
\ge  \frac{1}{(n-1)!} \int_\R \Psi(x)^{n-1} \Phi(x) x dx \ge  \frac{-\Psi(0)^{n-1}}{(n-1)!}  \int_{-\infty}^0 \big(1- \Psi(x)\big) dx \ .
\end{equation*}
These estimates show that there exists a positive constant $C$ which only depends on the shape $\Psi$ such that for any $l=1,\dots,n$,
 \begin{equation*} 
n! \left| \int_{\R^n_<}  \prod_{i=1}^{n} \Phi(x_i) \E{n}{x_{\sigma(\s_l)}} d^nx   \right| \le C n \ .
 \end{equation*} 
Using the estimate of lemma~\ref{thm:combinatorics_1} and (\ref{B_2}) this implies that
 $|\B^n_\Psi| \le C (n+1)! 2^{n-1}  $.
 Hence, for any $n \ge 2$, the Poisson component is bounded by
 \begin{equation} \label{growth_1} 
\left|  \B^n_\Psi \int_{\R} f(t)^n dt \right| \le C' (n+1)! \big(2\|f\|_\infty\big)^{n-1} \|f\|_{L^1} \ .
 \end{equation}

In the second half of the proof, we estimate the growth of $\mathfrak G ^n_{\Psi,\tau'}[f]$, (\ref{GG}). Applying the upper-bound (\ref{G_bound}), we see that
\begin{equation} \label{growth_2}
\left|\mathfrak G ^n_{\Psi,\tau'}[f] \right| \le 2 n! \sum_{|{\bf m}|=n} | \M({\bf m})|  \int_{R^n_0} \prod_{i=1}^n \big|\hat{f}(u_i)\big| \left( |u_1|+\cdots +|u_n| \right) d^{n-1}u  \int_{\R^n_<}\prod_{i=1}^n \Phi(x_i) d^nx    \ .
 \end{equation}
 By symmetry 
 $\displaystyle  \int_{\R^n_<}\prod_{i=1}^n \Phi(x_i) d^nx  =\frac{1}{n!} \left( \int_\R \Phi(x) dx \right)^n =\frac{1}{n!} $ and lemmas~\ref{thm:combinatorics_1} and ~\ref{thm:integrability_condition} provide bounds for the other factors of the RHS of (\ref{growth_2}). 
We obtain
\begin{equation}\label{growth_3}
 \left|\mathfrak G ^n_{\Psi,\tau'}[f] \right|
 \le (n+1)! 4^n \left( \|f\|_\infty+ \|f\|_{H^1}\right)^n \ .
 \end{equation}

The estimates (\ref{growth_1}) and (\ref{growth_3}) show that the limits of theorem~\ref{thm:C^n} satisfy the criterion of lemma~\ref{thm:Cumulant_Problem} for any choice of parameters $\tau'>0$, $\Psi\in\F$, and $f\in H^1_0(\R)$. Hence they corresponds to the cumulants  of some random variable which  we denote by $\Xi_{\Psi,\tau'} f$, and $\Xi f \Rightarrow \Xi_{\Psi,\tau'} f$ as $N\to\infty$. \qed\\

For the critical modified CUEs, according to formulae (\ref{eta}-$i$) and (\ref{eta_condition_1}), the parameter $\beta=1$. Hence, theorem~\ref{thm:CUE_C^n}  follows directly from proposition~\ref{thm:CUE_kernel_sine} and corollary~\ref{thm:weak_convergence}.
Likewise, in the GUE setting, the parameter $\beta=1/4$  by formula (\ref{eta}-$ii$). Provided that $1/3<\alpha<1$, by   proposition~\ref{thm:GUE_kernel_sine}, we conclude that  at the critical scale, a modified GUE  with shape $\Psi\in\F^*$ converges in distribution to the random field $\Xi_{\Psi,\tau/4}$.
In order to deal with all mesoscopic scales, we can use the asymptotic expansion of theorem~\ref{thm:meso_sine_kernel} instead of theorem~\ref{thm:micro_sine_kernel}. Namely, if we combine formula (\ref{kernel_approximation_1}) with the sine-kernel approximation (\ref{sine_approximation}), we obtain for any scales $0<\alpha,\delta <1$,
 \begin{align}\label{kernel_approximation_2}
&N^{-\delta}K_{\Psi,\alpha}^N(xN^{-\delta},yN^{-\delta}) \\
&\notag =\frac{1}{\tau N^{\alpha}} \sum_{|k| \le \Gamma_N N^\alpha } 
 \Phi\left(\frac{k}{\tau N^\alpha}\right) 
 \frac{ \sin\left[ (N+k) \left( F\left(\frac{\pi}{2} \frac{x N^{1/2-\delta}}{\sqrt{N+k}} \right)- F\left(\frac{\pi}{2} \frac{y N^{1/2-\delta}}{\sqrt{N+k}} \right) \right)\right] }{\pi(x-y)} 
+\underset{N\to\infty}{O}(N^{-\delta})  \ .
\end{align}
The RHS of (\ref{kernel_approximation_2}) is not a translation-invariant, so we cannot defined its Fourier transform. However, it is related to the kernel $L^N_{\Psi,\eta}$, (\ref{kernel_L}), by a change of variables and we can exploit this fact to compute the limits  of critical linear statistics of the modified GUE at any scale, including the regime $0<\alpha\le 1/3$.

\begin{proposition}  \label{thm:MNS_trace}
Let $\Psi\in\F^*$, $f\in C_0(\R)$, and $0<\alpha < 1$. For any $n\ge 2$
\begin{equation*}
\lim_{N\to \infty}\Cu^n_{K_{\Psi,\alpha}^N}[\Xi f_\alpha] =  \lim_{N\to \infty}\Cu^n_{L^N_{\Psi,\eta}}[\Xi g_N] \ ,
\end{equation*}
where  
\begin{equation} \label{g_function}
g_N(x)= f\left( \frac{2}{\pi} N^{\alpha} G\left(\frac{\pi x }{N^{\alpha}} \right) \right) 
\end{equation}
 and the function $G$ is given by definition $\ref{F}$. 
\end{proposition}

\proof 

Let $\supp(f) \subset [-L,L]$.
Observe that for any $|k| \le \Gamma_N N^\alpha $ and any  $x,y \in [-2L,2L]$, a Taylor expansion give
\begin{equation*}
F\left(x \frac{ N^{1/2-\alpha}}{\sqrt{N+k}} \right)- F\left(y \frac{ N^{1/2-\alpha}}{\sqrt{N+k}} \right)
= \sqrt{\frac{N}{N+k}} \left\{ F\left( \frac{x}{N^\alpha} \right) - F\left( \frac{y}{N^\alpha}\right)\right\}
+\underset{N\to\infty}{O}\left((x-y)\frac{\Gamma_N N^{-\alpha}}{N+k} \right)  \ .
\end{equation*}
Thus, 
taking $\delta=\alpha$ in equation (\ref{kernel_approximation_2}), we get for any $0<\alpha<1$,
 \begin{align}  \label{kernel_approximation_0}
&N^{-\alpha}K_{\Psi,\alpha}^N(xN^{-\alpha},yN^{-\alpha}) \\
&\notag= \frac{1}{\tau N^\alpha} \sum_{|k| \le \Gamma N^\alpha } 
 \Phi\left(\frac{k}{\tau N^\alpha}\right) 
 \frac{ \sin\left[ \sqrt{N(N+k)} \left( F\left(\frac{\pi}{2} \frac{x}{N^\alpha} \right)- F\left(\frac{\pi}{2} \frac{y}{N^\alpha} \right) \right) \right] }{\pi(x-y)} +\underset{N\to\infty}{O}(N^{-\alpha}) \ ,   
\end{align}
 where the error term is uniform for all $x,y \in [-2L,2L] $.
Following the proof of  lemma~\ref{thm:lemma_perturbative}, this approximation implies that for any any composition ${\bf m}$,
\begin{align} \label{kernel_approximation_3}
\tr[f^{m_1}_\alpha K_{\Psi,\alpha}^N &\cdots f^{m_\ell}_\alpha K_{\Psi,\alpha}^N]  
=\sum_{\begin{subarray}{c} k\in\Z^\ell \\ |k_j| \le  \Gamma_N N^\alpha \end{subarray}}  
\prod_{j=1}^\ell \Phi\left(\frac{k_j}{\tau N^\alpha}\right) \times \\
&\notag
 \int_{[-L,L]^\ell} \prod_{j=1}^\ell  f(x_j)^{m_j} 
\frac{ \sin\left[\sqrt{N(N+k)}\left( F\left(\frac{\pi}{2} \frac{x_j}{N^{\alpha}} \right)- F\left(\frac{\pi}{2} \frac{x_{j+1}}{N^{\alpha}} \right) \right)\right] }{\pi(x_j-x_{j+1})}  d^\ell x\
 +\underset{N\to\infty}{\bar{O}}(N^{-\alpha}) \ ,
\end{align}
where $x_{\ell+1}=x_1$. There exists $N_L \in\N$ such that  for all $N \ge N_L$, we can make the change of variables $y_j = \pi^{-1} N^{\delta} F\left(\frac{\pi}{2} \frac{x_j}{N^{\delta}} \right)$ in the integral (\ref{kernel_approximation_3}).
Since $0<F'(x) \le 2$ for any $|x| <1$, this change of variables maps the interval $|x_j|<L$ to some subset of $|y_j|<L$ for any $N \ge N_L$. Hence, if we let $g_N(y)= f\left( \frac{2}{\pi} N^{\alpha} G\left(\frac{\pi y }{N^{\alpha}} \right) \right)$ and 
 $\eta(k)$ is given by  (\ref{eta}-$ii$), we obtain
\begin{align} \label{kernel_approximation_4}
&\int_{[-L,L]^\ell} \prod_{j=1}^\ell  f(x_j)^{m_j} 
\frac{ \sin\left[\sqrt{N(N+k)}\left( F\left(\frac{\pi}{2} \frac{x_j}{N^{\alpha}} \right)- F\left(\frac{\pi}{2} \frac{x_{j+1}}{N^{\alpha}} \right) \right)\right] }{\pi(x_j-x_{j+1})}  d^\ell x  \\
&\notag\hspace{4cm} 
 = \int_{[-L,L]^\ell} \prod_{j=1}^\ell  g_N(y_j)^{m_j} 
\frac{ G'\left(\frac{\pi y_j }{N^{\alpha}} \right) \sin\left[2 \pi \eta(k)\left( y_j-y_{j+1} \right)\right] }{N^\alpha\left(G\left(\frac{\pi y_j }{N^{\alpha}} \right) - G\left(\frac{\pi y_{j+1} }{N^{\alpha}} \right)\right)}  d^\ell x \ .
\end{align}

A Taylor expansion gives for any $y,z\in [-4L,4L]$,
\begin{equation*}  G'\left(yN^{-\alpha} \right)^{-1} N^\alpha\left\{G(y N^{-\alpha})-G( zN^{-\alpha})\right\}
= (y-z) \left\{ 1+ \underset{N\to\infty}{O}\left( (y-z) N^{-\alpha}\right) \right\}  \ .
\end{equation*} 
This implies that for any $|k| \le \Gamma N^\alpha $ ,
\begin{equation} \label{G_approximation}
\frac{ G'\left(\frac{\pi y }{N^{\alpha}} \right) \sin\left[2 \pi \eta(k)\left( y-z \right)\right] }{N^\alpha\left(G\left(\frac{\pi y }{N^{\alpha}} \right) - G\left(\frac{\pi z }{N^{\alpha}} \right)\right)} 
=\frac{ \sin\left[2 \pi \eta(k)\left(  y-z \right)\right] }{\pi(y-z)} 
+ \underset{N\to\infty}{O}\left( N^{-\alpha}\right) \ .
\end{equation}

Hence, if we combine formulae (\ref{kernel_approximation_3}), (\ref{kernel_approximation_4}) and (\ref{G_approximation}), we have proved that

\begin{align*} 
\tr[f^{m_1}_\alpha K_{\Psi,\alpha}^N &\cdots f^{m_\ell}_\alpha K_{\Psi,\alpha}^N]  \\
&=\sum_{\begin{subarray}{c} k\in\Z^\ell \\ |k_j| \le  \Gamma N^\alpha \end{subarray}}  
\prod_{j=1}^\ell \Phi\left(\frac{k_j}{\tau N^\alpha}\right) 
 \int_{[-L,L]^\ell} \prod_{j=1}^\ell  g_N(y_j)^{m_j}\frac{ \sin\left[2 \pi \eta(k)\left(  y-z \right)\right] }{\pi(y-z)} 
+ \underset{N\to\infty}{\bar O}\left( N^{-\alpha}\right) \ . 
\end{align*}

By (\ref{kernel_L}), we can write this equation as
\begin{equation*}
\tr[f^{m_1}_\alpha K_{\Psi,\alpha}^N \cdots f^{m_\ell}_\alpha K_{\Psi,\alpha}^N] 
= \tr[g_N^{m_1} L^N_{\Psi,\eta} \cdots g_N^{m_\ell} L^N_{\Psi,\eta}] +\underset{N\to\infty}{\bar{O}}(N^{-\alpha})  \ ,
\end{equation*}
and the proposition follows from  formula (\ref{cumulant_2}). \qed\\

By proposition~\ref{thm:MNS_trace}, to complete the proof of theorem~\ref{thm:MNS_C^n}, it remains to extend the argument of theorem~\ref{thm:C^n} to deal with test functions of the form (\ref{g_function}).
As we can see from the next lemma, such functions depend mildly on the density $N$  and it is not difficult to finish the proof.

 \begin{lemma}\label{thm:g_N} Let $f\in C^1(\R)$ with support in $[-L,L]$ and $0<\delta \le 1$. For any $N > (2 L)^{1/\delta}$, the function $g_N(x)= f\left( \frac{2}{\pi} N^{\delta} G\left(\frac{\pi x }{N^{\delta}} \right) \right)$ has compact support in $[-L,L]$. Moreover, we have 
\begin{equation*}
   \|\hat{g}_N-\hat{f}\|_{\infty} =\underset{N\to\infty}{O}(N^{-\delta}) 
   \hspace{.5cm} \text{and} \hspace{.6cm}
 \|g_N-f\|_{H^1}=\underset{N\to\infty}{o(1)} \ .
\end{equation*}
\end{lemma}

\proof  By definition \ref{F},  $0\le F' \le 2$ and the map $x\mapsto \frac{2}{\pi} N^{\delta} G\left(\frac{\pi x }{N^{\delta}} \right)$ is a dilation. Therefore, when $N > (2 L)^{1/\delta}$, the map $g_N$ is well-defined on $[-L,L]$ and $g_N(x)=0$ for all $x\in [-L,L]\backslash \supp(f)$. By continuity, we can assume that $g_N(x)=0$ for all $|x|>L$.
Hence $g_N\in C^1_0(\R)$ with $\supp(g_N) \subset \supp(f)$. 
Then, by Lipschitz continuity of $f$, for any $|x|<L$,
$$| g_N(x)- f(x)| \le CN^\delta \left| 2G\left(\frac{\pi x }{N^{\delta}}\right)- \frac{\pi x }{N^{\delta}} \right|  
\le C L^2 N^{-\delta} \ ,   $$
where we used that $G$ is smooth with $G(0)=0$ and $G'(0)=1/2$. This implies that
$$  \|\hat{g}_N-\hat{f}\|_{\infty} \le \int_{[-L,L]} | g_N(x)- f(x)|  dx = O(N^{-\delta}) \ .$$ 
Similarly, by the triangle inequality,
$$| g_N'(x)- f'(x)| \le \|f'\|_{\infty} \left| 2G'\left(\frac{\pi x }{N^{\delta}}\right)-1\right| 
+  \left|f'\left( \frac{2}{\pi} N^{\delta} G\left(\frac{\pi x }{N^{\delta}} \right) \right)- f'(x)\right| \ .$$
Since $f'$ is continuous, 
$\displaystyle \lim_{N\to\infty}| g_N'(x)- f'(x)| = 0 $ for all $x\in [-L,L]$. By (\ref{norm_1}) and  the dominated convergence theorem, we conclude  that 
$$  \|g_N-f\|_{H^1}^2 =\frac{1}{4\pi^2} \int_{[-L,L]} | g_N'(x)- f'(x)|^2 dx  \to 0 $$
as $N\to\infty$. \qed \\

 \begin{proposition} \label{thm:C^n'}
Let $f \in C^1_0(\R)$, $0<\alpha<1$, and $\Psi \in\F^*$. If $\eta$ is given by $(\ref{eta}$-$ii)$ and $g_N$ is given by $(\ref{g_function})$, then for any $n\ge2$,
\begin{equation*}
\lim_{N\to\infty}\Cu^n_{L^N_{\Psi,\eta}}[\Xi g_N] = \Cu^n\big[\Xi_{\Psi,\tau/4} f\big] \ .
\end{equation*}
\end{proposition}

\proof We can take $f=g_N$ in formula (\ref{C^n_2}). In particular the estimate (\ref{C^n_3}) is still valid  for the test function $g_N$. Let $|u|_1=|u_1|+\cdots+|u_n|$ and the function $ \Upsilon_N^n(u)$ be given by (\ref{Upsilon_0}).  We want to prove that 
\begin{equation}  \label{g_estimate_1}
\lim_{N\to\infty} \int_{\R^n_0} \1_{\left\{|u|_1 >\frac{N^{1-\alpha}}{2}\right\}}   |u|_1  \prod_i \left| \hat{g}_N(u_i) \right| d^{n-1}u =0 \ ,
\end{equation}
and 
\begin{equation} \label{g_estimate_2}
\lim_{N\to\infty} \int_{\R^n_0} \left|  \prod_i \hat{g}_N(u_i) - \prod_i \hat{f}(u_i)\right| 
\left| \Upsilon_N^n(u) \right| d^{n-1}u =0  \ .
\end{equation}

Indeed  these limits imply that
\begin{align*}
 \lim_{N\to\infty}\Cu^n_{L^N_{\Psi,\eta}}[\Xi g_N]  
&=2\lim_{N\to\infty} \int_{\R^n_0} \Re\left\{ \prod_i \hat{g}_N(u_i) \right\} \Upsilon^n_N(u)\ d^{n-1} u \\
&=2\lim_{N\to\infty} \int_{\R^n_0} \Re\left\{ \prod_i \hat{f}(u_i) \right\} \Upsilon^n_N(u)\ d^{n-1} u \ , 
\end{align*}
and the rest of the proof is identical to the proof of theorem~\ref{thm:C^n}. 
To complete our argument, it remains to show (\ref{g_estimate_1}) and (\ref{g_estimate_2}).
First observe that
\begin{equation} \label{C^n_6}
  \int_{\R^n_0} \1_{\left\{|u|_1 >\frac{N^{1-\alpha}}{2}\right\}}   |u|_1  \prod_i \left| \hat{g}_N(u_i) \right| d^{n-1}u 
\le \sum_{k=1}^n \int_{\R^n_0} \1_{\left\{|u_k| >\frac{N^{1-\alpha}}{2n}\right\}}   |u|_1  \prod_i \left| \hat{g}_N(u_i) \right| d^{n-1}u  \ .
\end{equation}
Let $\mathcal A_N =\left\{ v \in \R : |v| >\frac{N^{1-\alpha}}{2n}\right\}$ and define the function $q_N$ by its Fourier transform $\hat{q}_N= \1_{\mathcal A_N}\hat{g}_N$.
Then, by  lemma~\ref{thm:integrability_condition}, we have
\begin{equation} \label{g_estimate_6}
\int_{\R^n_0} \1_{\left\{|u|_1 >\frac{N^{1-\alpha}}{2}\right\}}   |u|_1  \prod_i \left| \hat{g}_N(u_i) \right| d^{n-1}u 
\le n^22^{n-1} \left( \|\hat{q}_N \|_\infty+ \|q_N\|_{H^1}\right) \left( \|\hat{g}_N\|_\infty+ \|g_N\|_{H^1}\right)^{n-1}  \ . 
\end{equation}

Since, $\|F'\|_\infty\le 2$, a change of variables yields
\begin{equation} \label{g_estimate_3}
 \|\hat{g}_N\|_{\infty} \le \int \left|  f\left( \frac{2}{\pi} N^{\alpha} G\left(\frac{\pi x }{N^{\alpha}}  \right) \right) \right| dx  \le 
\int \left|f(z)\right| dz = \| f \|_{L^1} \ .  \end{equation}
Moreover, 
$$
\|g_N\|_{H^1}^2 =  \int  \big| g_N'(x)\big|^2 dx = 2  \int_{-L}^L  \big| f'(z)\big|^2 G'\left(F \left( \frac{\pi z }{2N^{\alpha}} \right)\right) dz ,  
$$
and, since $G'(F(t)) = F'(t) = 2 \sqrt{1-t^2}$ for any $|t|<1$, we obtain that $
\|g_N\|_{H^1}^2 \le 4 \| f\|_{H^1}^2$  for all $N$ sufficiently large. 
By (\ref{g_estimate_6}), this implies that
\begin{equation*}
\int_{\R^n_0} \1_{\left\{|u|_1>\frac{N^{1-\alpha}}{2}\right\}}   |u|  \prod_i \left| \hat{g}_N(u_i) \right| d^{n-1}u 
\le n^22^{n-1} \left( \|\hat{q}_N \|_\infty+ \|q_N\|_{H^1}\right) \left( \|\hat{f}\|_{L^1}+ 2 \|f\|_{H^1}\right)^{n-1} \ . 
\end{equation*}
Obviously $ \|\hat{q}_N \|_\infty \to 0$ and to conclude that (\ref{g_estimate_1}) holds it remains to estimate $\|q_N\|_{H^1}$. The main observation is that
$$\int \left| \hat{q}_N(v)- \1_{\mathcal A_N}\hat{f}(v)\right|^2 |v|^2 dv \le \|g_N- f \|^2_{H^1} \ .  $$ 
Then, by the triangle inequality,
$$ \|q_N\|_{H^1}^2 \le 2 \left(\|g_N- f \|^2_{H^1} + \int_{\mathcal A_N}\left| \hat{f}(v)\right|^2 |v|^2 dv\right) \ . $$
The first term converges to 0 by lemma~\ref{thm:g_N}, and so does the second term  by the dominated convergence theorem. We conclude that $\|q_N\|_{H^1} \to 0$ and we have established (\ref{g_estimate_1}).
The proof of the estimate (\ref{g_estimate_2}) is very similar.
We observe that
$$ \prod_i \hat{g}_N(u_i) - \prod_i \hat{f}(u_i) = \sum_{j=1}^n \left( \hat{g}_N(u_j) - \hat{f}(u_j)\right) 
\prod_{i>j}\hat{g}_N(u_i) \prod_{i<j} \hat{f}(u_i) \ . $$ 
By (\ref{Upsilon_4}), there exists a constant $C>0$ which  depends only on $n$ and the shape $\Psi$ such that $\left| \Upsilon_N^n(u) \right| \le C \left\{1+|u|_1\right\} $ and  we obtain
\begin{align*}
 &\int_{\R^n_0} \left|  \prod_i \hat{g}_N(u_i) - \prod_i \hat{f}(u_i)\right| 
\left| \Upsilon_N^n(u) \right| d^{n-1}u \\
&\hspace{3cm}\le C \sum_{j=1}^n  \int_{\R^n_0}    \left| \hat{g}_N(u_j) - \hat{f}(u_j)\right| \left| \prod_{i>j}\hat{g}_N(u_i) \prod_{i<j} \hat{f}(u_i) \right| \big\{1+ |u|_1 \big\} d^{n-1}u  \ .
\end{align*}
Thus, by lemma~\ref{thm:integrability_condition} and the facts $\|\hat{g}_N\|_{\infty} \le \|f\|_{L^1}$ and $\|g_N\|_{H^1} \le 2 \|f\|_{H^1}$,   we get
\begin{align} \label{g_estimate_7}
& \int_{\R^n_0} \left|  \prod_i \hat{g}_N(u_i) - \prod_i \hat{f}(u_i)\right| 
\left| \Upsilon_N^n(u) \right| d^{n-1}u\\
&\notag\hspace{4cm}\le  Cn^22^{n-1}  \left( \| \hat{g}_N-\hat{f} \|_{\infty} + \|  g_N-f \|_{H^1} \right) \left( \| f \|_{L^1} + 2 \|f \|_{H^1} \right)^{n-1} \ .
\end{align}
 Lemma~\ref{thm:g_N} implies that the RHS of equation (\ref{g_estimate_7}) converges to 0 as $N\to\infty$ and the limit (\ref{g_estimate_2}) follows.  \qed\\

By proposition~\ref{thm:MNS_trace} and~\ref{thm:C^n'}, if $f\in C^1_0(\R)$, for any $0<\alpha < 1$ and any $n\ge 2$,
\begin{equation} \label{cumulant_convergence}
\lim_{N\to \infty} \Cu^n_{K_{\Psi,\alpha}^N}[\Xi f_\alpha] 
= \Cu^n\big[\Xi_{\Psi,\tau/4} f\big] \ .
\end{equation}

In the proof of corollary~\ref{thm:weak_convergence}, we have shown that the sequence $\left(  \Cu^n\big[\Xi_{\Psi,\tau/4} f\big] \right)_{n\ge2}$ given by theorem~\ref{thm:C^n}  satisfies the condition of lemma~\ref{thm:Cumulant_Problem}. This 
 implies that,  considering the determinantal process with correlation kernel $K_{\Psi,\alpha}^N$, the random variable  
 \begin{equation} \label{weak_convergence}
 \Xi f_\alpha \Rightarrow \Xi_{\Psi,\tau/4} f
 \end{equation} as $N\to\infty$. In order to complete the proof of theorem~\ref{thm:MNS_C^n}, we use a density argument to extend (\ref{weak_convergence}) to all test functions in $H^{1/2}_0 \cap L^\infty(\R)$.\\

{\it Proof of theorem~\ref{thm:MNS_C^n}.}
Let $\tau'=\tau/4$. First observe that for any $x,y \in \R$,
$$ \left| e^{ix} -e^{iy} \right|^2 \le 4 |x-y|^2  \ .$$  
By Chebychev's inequality, this implies that, if $X$ and $Y$ are mean-zero random variables defined on the same probability space, for any $\xi\in\R$,
\begin{equation*} 
 \left| \E{}{e^{i\xi X } -e^{i\xi Y}} \right| \le 4 |\xi| \sqrt{ \Var\left[ X-Y \right] } \ .
\end{equation*}
 For critical linear statistics of the modified GUE,
using the estimate (\ref{variance_estimate}) which is valid when $\delta=\alpha$,  we obtain for any test functions $f,h\in H^{1/2}_0(\R)$,     
\begin{equation} \label{characteristic_1}
 \left| \E{K^N_{\Psi,\alpha}}{e^{i\xi \Xi f_\alpha } -e^{i\xi \Xi h_\alpha}} \right| 
 \le  C |\xi| \sqrt{ \|f-h\|_{L^2}^2 +  \|f-h\|_{H^{1/2}}^2} \ .
\end{equation}

Moreover, formula (\ref{C2}) implies that under the same assumptions,
\begin{equation} \label{characteristic_2}
 \left| \mathbb{E}\left[e^{i\xi \Xi_{\Psi,\tau'} h}-e^{i\xi \Xi_{\Psi,\tau'} f}\right] \right| 
 \le C |\xi|  \sqrt{  \|f-h\|_{L^2}^2 +    \|f-h\|_{H^{1/2}}^2 } \ .
\end{equation}

By the triangle inequality, 
\begin{align} \notag 
& \big| \E{K^N_{\Psi,\alpha}}{e^{i\xi \Xi f_\alpha }} - \E{}{e^{i\xi \Xi_{\Psi,\tau'} f}} \big| \\
& \label{characteristic_3} \hspace{.7cm}\le  \left| \E{K^N_{\Psi,\alpha}}{e^{i\xi \Xi f_\alpha } -e^{i\xi \Xi h_\alpha}} \right| 
 +  \left| \mathbb{E}\left[e^{i\xi \Xi_{\Psi,\tau'} h}-e^{i\xi \Xi_{\Psi,\tau'} f}\right] \right|   
  +  \left| \E{K^N_{\Psi,\alpha}}{e^{i\xi \Xi h_\alpha }} - \E{}{e^{i\xi \Xi_{\Psi,\tau'} h}} \right| \ . 
\end{align}
If we suppose that $h\in C^1_0(\R)$, by (\ref{weak_convergence}), the last term in the RHS of (\ref{characteristic_3}) converges to 0 as $N\to\infty$. 
Thus, using the upper-bound (\ref{characteristic_1}) and (\ref{characteristic_2}), for any $f\in H^{1/2}_0\cap L^\infty(\R)$ and $\xi\in\R$, 
\begin{equation}  \label{characteristic_4}
\limsup{N\to\infty}  \left| \E{K^N_{\Psi,\alpha}}{e^{i\xi \Xi f_\alpha }} - \E{}{e^{i\xi \Xi_{\Psi,\tau'} f}} \right| 
\le 2 C  |\xi|  \sqrt{  \|f-h\|_{L^2}^2 +    \|f-h\|_{H^{1/2}}^2 } \ . 
\end{equation}

Since, the space $C^1_0$ is dense in the Sobolev space $H^{1/2}_0$ with respect to the norm $\sqrt{\|\cdot\|_{L^2}^2 +  \|\cdot\|_{H^{1/2}}^2}$, \cite[Theorem~7.14]{Lieb_Loss},  the RHS of the inequality  (\ref{characteristic_4}) is arbitrary small by choosing $h\in C^1_0(\R)$ appropriately, and we conclude that $\Xi f_\alpha \Rightarrow \Xi_{\Psi,\tau'} f$ as $N\to\infty$. \qed

\subsection{Properties of the random process $\Xi_{\Psi,\tau}$} \label{sect:Poisson}


In this section, we study the random variables $\Xi_{\Psi,\tau} f$ which arise from the limit of linear statistics of the critical modified ensembles. 
Because of the complicated structure of the cumulants in theorem~\ref{thm:MNS_C^n}, we cannot get much information about the random fields $\Xi_{\Psi,\tau}$ except that they are not Gaussian.
However, as we expect from figure~\ref{fig:phase}, we recover Gaussian fluctuations in both limits  $\tau\to\infty$ or~$\tau=0$; see proposition~\ref{thm:double_limit}. 
Proposition~\ref{thm:Gaussian} provides a sufficient condition under which the field $\Xi_{\Psi,\tau}$ is not Gaussian and it leads us to compute the {\it Laplace transform} of the Poisson component of the field  $\Xi_{\Psi,\tau}$. 
In particular, we establish proposition \ref{thm:MNS_property}, i.e$.$ we show that the Poisson component of the random field  $\Xi_{\Psi,\tau}$ is Gaussian if and only if $\Psi$ is the MNS shape $\psi(t)=(1+e^t)^{-1}$.
Then, a natural problem that remains unanswered is whether the sequence $\mathfrak G_{\psi,\tau}$ given by (\ref{GG}) also corresponds to the cumulants of some random variables, so that the field $\Xi_{\psi,\tau}$ would be the superposition of a {\it white noise} and an independent non-Gaussian process.


\begin{proposition} \label{thm:Gaussian}
 If the shape $\Psi\in\F$ satisfies the condition $\B^n_\Psi\neq 0$ for some $n>2$. Then, for any $\tau>0$, the random process $\Xi_{\Psi,\tau} $ of corollary $\ref{thm:weak_convergence}$ is not Gaussian.
\end{proposition}

\proof

It is clear from the definition (\ref{G}) that
$\displaystyle \lim_{\tau\to\infty}\G_\tau^{\bf m}(u,x)=0 $ 
and it follows from  (\ref{GG}) that for any $n \ge 2$ and $f\in H^1_0(\R)$,
$\displaystyle \lim_{\tau\to\infty}\mathfrak G^n_{\Psi,\tau}(f)=0 $. Hence, by formula (\ref{C_split})

\begin{equation} \label{temperature_infty}
\Cu^n\big[\Xi_{\Psi,\tau} f\big] 
= 2\tau \B^n_\Psi \int_{\R} f(t)^n dt + \underset{\tau\to\infty}{o(1)} \ .
\end{equation}
Thus, the Poisson component dominates at large temperature and the random field of $\Xi_{\Psi,\tau} f$  is not Gaussian since there are test functions such that $\Cu^n[\Xi_{\Psi,\tau} f] \neq 0$ whenever $\B^n_\Psi\neq 0$.
This observation is actually valid at any temperature $\tau>0$ because of the scaling property of the cumulants. By definition,
$ \G_\tau^{\bf m}(u,x) =  \tau\G_1^{\bf m}(u/\tau,x) $ 
and the change of variables $u_i= \tau v_i$ leads to
\begin{equation*}
 \mathfrak G ^n_{\Psi,\tau}[f] =- 2\int_{\R^n_0} dv\ \Re\left\{\prod_{i=1}^n \tau \hat{f}(\tau v_i) \right\} \int_{\R^n_<}\prod_{i=1}^n \Phi(x_i) \sum_{|{\bf m}|=n} \M({\bf m}) \G_1^{\bf m}(v,x) dx  \ . \\
\end{equation*}
Hence, by (\ref{C_split}), the random variables 
$\Xi_{\Psi,\tau} f$ and $\Xi_{\Psi,1} f(\frac{\cdot}{\tau})$ have the same distribution. \qed

\begin{proposition} \label{thm:double_limit}
 For any function $f \in H^1_0(\R)$, the rescaled random variable $\tau^{-1/2} \Xi_{\Psi,\tau} f$  converges in  distribution as $\tau\to\infty$ to a Gaussian random variable with variance $\|f\|_{L^2}^2$.
On the other hand,  $ \Xi_{\Psi,\tau} f$  converges in  distribution as $\tau\to0$ to a Gaussian random variable with variance $\|f\|_{H^{1/2}}^2$. 
\end{proposition}

\proof When $\tau \to \infty$, the asymptotic of the cumulants  of the random variables $\Xi_{\Psi,\tau} f$ are given by formula~(\ref{temperature_infty}) and 
\begin{equation*}
\Cu^n\big[\Xi_{\Psi,\tau}(\tau^{-1/2} f )\big] = 
 \tau^{ 1-n/2} \B^n_\Psi \int_{\R} f(t)^n dt + \underset{\tau\to\infty}{o}( \tau^{-n/2}) \ .
\end{equation*}
Hence,
\begin{equation*}
\lim_{\tau\to\infty}  \Cu^n\big[ \tau^{-1/2}\Xi_{\Psi,\tau}( f )\big] =  \begin{cases}
 \B^2_\Psi\|f\|_{L^2}^2 &\text{if }n=2 \\ 0 &\text{if } n \ge 3 \end{cases} \ .
\end{equation*}
Taking the limit as $\tau\to0$ is more subtled. We shall see that we recover the  cumulants of the sine process given by lemma 2 in~\cite{Soshnikov_00a} and the cancellation follows from the main combinatorial lemma. We fix some composition ${\bf m}$ of $n$ and some vector $u\in\R^n_0$ and we will look at the symmetries of the function $\G_0^{\bf m}$. By definition (\ref{G}),
\begin{equation} \label{G_0}
\G_0^{\bf m}(u, x)=  \sum_{\sigma\in \Sy(n)} \max_{i \le \ell}\left\{\Lambda^{\bf m}_{i,\s}(u) \right\}
\end{equation}
where $\s=\s_{\ell({\bf m})}(\sigma)$, (\ref{N_3}). The important fact is that this expression becomes independent of the variable~$x$.
In the sequel, we will denote $\G_0^{\bf m}(u)$ instead of $\G_0^{\bf m}(u, x)$ and we define
\begin{equation} \label{triangle}
 \triangle^{\bf m}(u)= (u_1+\cdots+u_{\overline{ m}_1},\ u_1+\cdots+u_{\overline{ m}_2},\ \dots,\ u_1+\cdots+u_{\overline{ m}_{\ell-1}},\ 0)  \ .
\end{equation}
By definition (\ref{Lambda}), for any $u\in\R_0^n$, we have 
\begin{equation*} 
\Lambda^{\bf m}_{i, s}(u) 
=\begin{cases} u_{\overline{ m}_{s}+1}+\cdots+ u_{\overline{ m}_i} &\text{if } s< i \\
u_{\overline{ m}_{s}+1}+\cdots+ u_n+ u_1+ \cdots+ u_{\overline{ m}_i} &\text{if } i<s \\
0 &\text{if } i=s
\end{cases}
 \end{equation*}
For any $s=1,\dots,  \ell({\bf m})$ we let $\pi_s\in \Sy(n)$ be the cyclic permutation given by
\begin{equation*} \pi_s(i) = \overline{ m}_{s} + i\mod n \ .\end{equation*}
Then, we see that $\big\{ \Lambda^{\bf m}_{i, s}(u) : i=1,\cdots,\ell\big\} = \big\{ \triangle_i^{\bf m}(\pi_s u): i=1,\cdots,\ell\big\}  $ and, by (\ref{G_0}), we obtain
\begin{align} \notag
\G_0^{\bf m}(u) &=  \sum_{\sigma\in \Sy(n)} \max_{i \le \ell} \left\{ \triangle_i^{\bf m}(\pi_{\s} u) \right\} \ , \\
\label{zero_averaging}
\sum_{\pi \in \Sy(n)} \G_0^{\bf m}(\pi u) &=  n! \sum_{\pi \in \Sy(n)} \max_{i \le \ell}\left\{ \triangle_i^{\bf m}(\pi u) \right\}\ .
\end{align}

By dominated convergence, we can pass to the limit $\tau\to0$ in formulae (\ref{GG}) and  (\ref{C_split}). By (\ref{G_0}), the two integrals decouple and, since $\int_{\R^n_<}\prod_{i=1}^n \Phi(x_i) d^nx = 1/ n!$ , we obtain
\begin{equation} \label{C^n_0}
\lim_{\tau\searrow 0} \Cu^n\big[\Xi_{\Psi,\tau'} f\big] 
= - 2 \int_{\R^n_0} \Re\left\{ \prod_{i=1}^n \hat{f}(u_i) \right\}  \frac{1}{n!} \sum_{|{\bf m}|=n} \M({\bf m}) \G_0^{\bf m}(u)d^{n-1}u \ .\end{equation}
This limit is independent of the shape $\Psi$ and it will be denoted by $ \Cu^n\big[\Xi_0 f\big] $.
If we use the notation (\ref{triangle}), Soshnikov's main combinatorial lemma reads for any $u \in \R^n_0$,
\begin{equation} \label{MCL}
\sum_{\pi \in \Sy(n)} \sum_{|{\bf m}|=n} \M({\bf m})  \max\left\{ \triangle^{\bf m}(\pi u) \right\}
=  \begin{cases} -|u_1| &\text{if } n=2 \\ 0 &\text{if } n\ge3 \end{cases} \ .
\end{equation}
Next we symmetrize formula (\ref{C^n_0}) over all permutations of $u$, by equations (\ref{zero_averaging}) and (\ref{MCL}), we conclude that for any $n\ge 3$,
\begin{equation*}
 \Cu^n\big[\Xi_0 f\big] = 0
 \hspace{.5cm}\text{and}\hspace{.5cm}
 \Cu^2\big[\Xi_0 f\big] =  \int_\R   \hat{f}(u)  \hat{f}(-u) |u| du  \ .
\end{equation*}
This shows that the random field $\Xi_0$ is Gaussian with covariance structure $\langle f, g \rangle_{H^{1/2}}$. \qed\\

In the proof of proposition~\ref{thm:Gaussian}, we have seen that  $\tau$ is just a scaling parameter. 
Therefore, in the sequel,  we will assume that $\tau=1$ 
and write $\Xi_{\Psi}=\Xi_{\Psi,1}$, $ \G^{\bf m}=\G_{1}^{\bf m}$, etc.\\

By definition \ref{C_Poisson},  the behavior of the Poisson component of the field $\Xi_\Psi$ is encoded by the coefficients $\B^n_\Psi$, (\ref{B}). In the remainder of this section, we will compute the generating function of the sequence $\B^n_\Psi$ and prove proposition~\ref{thm:MNS_property}. We start by a combinatorial lemma.

\begin{proposition} \label{thm:b_GF}
 For any $z,w\in\C$ such that $\left|w(e^{(1+w)z}-1)\right|< |1+w|,$
 \begin{equation} \label{b_GF}
\sum_{n=1}^\infty \sum_{k=0}^{n-1} b_{k}^n \frac{w^{k+1} z^n}{n!}=\frac{w\big(e^{(1+w)z}-1\big)}{1+we^{(1+w)z}} \ .
\end{equation} \end{proposition}

\proof By equation (\ref{b}),
\begin{equation*} b^n_k=  \sum_{l=1}^{k+1} (-1)^{l+1} {n- l \choose k+1-l} 
\sum_{\begin{subarray}{c}n_1,\dots,n_l\ge 1 \\ n_1+\cdots+n_l=n \end{subarray} } \frac{n!}{n_1!\cdots n_l !} \ . \end{equation*} 
So that if we exchange the order of summation between $k$ and $l$,
\begin{align*} \sum_{n=1}^\infty  \sum_{k=0}^{n-1} b_{k}^n \frac{w^{k+1} z^n}{n!}
&=  \sum_{n=1}^\infty z^n \sum_{l=1}^{n}
\sum_{\begin{subarray}{c}n_1,\dots,n_l\ge 1 \\ n_1+\cdots+n_l=n \end{subarray} } \frac{(-1)^{l+1}}{n_1!\cdots n_l !}   \sum_{k=l-1}^{n-1}  {n- l \choose k+1-l} w^{k+1}\\
&=  \sum_{n=1}^\infty z^n \sum_{l=1}^{n} \sum_{\begin{subarray}{c}n_1,\dots,n_l\ge 1 \\ n_1+\cdots+n_l=n \end{subarray} }\frac{(-1)^{l+1}}{n_1!\cdots n_l !}  w^l(1+w)^{n-l}
\ .\end{align*}
Then, since
$ \displaystyle \sum_{n=l}^\infty a^n \sum_{\begin{subarray}{c}n_1,\dots,n_l\ge 1 \\ n_1+\cdots+n_l=n \end{subarray} } \frac{1}{n_1!\cdots n_l !}   = (e^a-1)^l $ for any $a\in\C$, if we exchange the order of summation between $l$ and $n$, we obtain
\begin{equation*} \sum_{n=1}^\infty  \sum_{k=0}^{n-1} b_{k}^n \frac{w^{k+1} z^n}{n!}
= \sum_{l=1}^\infty (-1)^{l+1}\left(\frac{w}{1+w}\right)^l \left( e^{(1+w)z}-1 \right)^l \ .
\end{equation*}
This proves (\ref{b_GF}) using the identity $\displaystyle\sum_{1\le l}  (-1)^{l+1} \xi^l = \frac{\xi}{1+\xi} $, if $|\xi|<1$. \qed

If we substitute $w=\frac{\Psi(x)}{1-\Psi(x)}$ and $z=\xi(1-\Psi(x))$, for any $|\xi|<e^{-1}$, into (\ref{b_GF}), we get
\begin{equation*} \sum_{n=1}^\infty \frac{\xi^n}{n!} \sum_{k=0}^{n-1} b_{k}^n \Psi(x)^{k+1}\big(1-\Psi(x)\big)^{n-k-1}=\frac{\Psi(x)\big(e^\xi-1\big)}{1+\Psi(x)\big(e^\xi-1\big)} \ .\end{equation*}

Integrating both sides, by definition (\ref{B}), this implies that 
\begin{equation} \label{B_GF_1}  
\sum_{n=1}^\infty \frac{\xi^n}{n!} \B^n_\Psi = \int_\R \frac{x(e^\xi-1)}{1+\Psi(x)\big(e^\xi-1\big)}  \Phi(x) dx \ .\end{equation}

 Motivated by proposition \ref{thm:Gaussian}, it is meaningful to raise the question:  which shape $\Psi$ satisfies the conditions $\B^n_\Psi =0$ for all $n>2$?
Since we assume that $\B^1_\Psi=0$ and $\B^2_\Psi =1$,  by formula (\ref{B_GF_1}), this amounts to solving the integral equation
\begin{equation} \label{B_GF_2}  
 \int_\R \frac{x(e^\xi-1)}{1+\Psi(x)\big(e^\xi-1\big)}  \Phi(x) dx
 =\frac{ \xi^2}{2} .
 \end{equation} 
 The proof of proposition~\ref{thm:MNS_property} is divided into two part, first we show that the equation (\ref{B_GF_2}) has a unique solution by transforming it into a {\it Cauchy integral}. Second, we verify that the  MNS shape $\psi(x)=1/(1+e^x)$ is a solution. \\

 {\it Proof of proposition~\ref{thm:MNS_property}.} 
 In general, since $(1-\Psi)$ is cumulative distribution function, we can make the change of  variable $s=\Psi(x)$ in the RHS of formula (\ref{B_GF_1}) and this leads to the identity 
\begin{equation} \label{B_GF}  
\sum_{n=1}^\infty \frac{\xi^n}{n!} \B^n_\Psi = \int_0^1  \frac{\Psi^{-1}(s)}{(e^\xi-1)^{-1}+ s } ds \ ,
\end{equation}
 where $\Psi^{-1}$ is the generalized inverse of $\Psi$,
 $$ \Psi^{-1}(s)= \inf \big\{ t \in \R : \Psi(t) \le t \big\} \ . $$

 If we also make the  change of variable $w=(1-e^\xi)^{-1}$, we see that equation (\ref{B_GF_2}) gives
 \begin{equation} \label{B_GF_3}  
\int_0^1  \frac{\Psi^{-1}(s)}{s-w} ds   =\frac{\sigma^2}{2} \left( \log(1 - w^{-1}) \right)^2 .
\end{equation}
 Note that, since $\Psi$ is continuous,  the function $\Psi^{-1}$  is continuous almost everywhere.
Moreover, the conditions $\Psi \in L^1(0,\infty)$ and  $(1-\Psi)\in L^1(-\infty,0)$ guarantees that the  RHS of formula (\ref{B_GF_3}) defines an analytic function in $\C\backslash[0,1]$. Thus, if we use the principal branch of the logarithm, then the  l.h.s$.$ of formula (\ref{B_GF_3}) is also analytic in the same domain. 
In particular, by analytic continuation, equation  (\ref{B_GF_3}) holds for any $w \in \C\backslash[0,1]$.
 It is well-known that,  for any $t\in(0,1)$ where the function $\Psi^{-1}$ is continuous,
 $$ \lim_{\eta\searrow 0} \Im \left( \int_0^1  \frac{\Psi^{-1}(s)}{s- (t+i \eta)} ds \right)  = \pi \Psi^{-1}(t) \ .
 $$
For the principal branch, for any $t\in(0,1)$,
$$ \lim_{\eta\searrow 0}  \log\big(1 -  (t+i \eta)^{-1}\big) =  \log\big( t^{-1} -1\big) + i \pi \  ,$$
 Hence, by equation (\ref{B_GF_3}), these limits implies that if $\Psi\in\F$ is solution of  (\ref{B_GF_1}), then  for almost all $t\in(0,1)$,
 \begin{equation} \label{MNS_property}
 \Psi^{-1}(t)=  \log\big( t^{-1} -1\big) \ . 
 \end{equation} 
 Since, by assumption, $ \Psi^{-1}$ is non-increasing, equation (\ref{MNS_property}) holds for all $t\in(0,1)$ and it is straightforward to see that this amounts to $\Psi(x)=1/(1+ e^{x})$.  
 
It remains to check that the solution of  equation (\ref{B_GF_2}) is indeed  $\psi(x)=1/(1+ e^{x})$. If we let $\rho=e^\xi-1$, by formula (\ref{B_GF_1}),
\begin{align} 
\sum_{n=1}^\infty \frac{\xi^n}{n!} \B^n_\psi
&\notag= \rho \int_\R \frac{x e^{x}}{(1+e^{x}) (1+e^x +\rho)} dx\\
&\label{PolyLog_integral_1}= \rho \int_0^\infty \frac{\log t }{(1+t) (1+t +\rho)} dt \ .    \end{align}

Thus, we want to prove that for any $\rho> -1$,
\begin{equation} \label{PolyLog_integral_2}
 \rho \int_0^\infty \frac{\log t }{(1+t) (1+t +\rho)} dt  =\frac{1}{2} \log^2(1+\rho) \ .
 \end{equation}
Indeed, since $\rho=e^\xi-1$, formulae (\ref{PolyLog_integral_1}) and (\ref{PolyLog_integral_2}) imply that
\begin{equation} \label{B_GF_MNS} \sum_{n=1}^\infty \frac{\xi^n}{n!} \B^n_\psi=\frac{1}{2}\log^2(e^\xi)= \frac{\xi^2}{2} \  . \end{equation}
To prove  (\ref{PolyLog_integral_2}), we can differentiate both sides with respect to the parameter $\rho$ and we see that it is enough to show that for any $\rho>-1$,
\begin{equation} \label{PolyLog_integral_3}
 \int_0^\infty \frac{\log t }{(1+t +\rho)^2} dt  =\frac{\log(1+\rho)}{1+\rho} \ .
 \end{equation}
For any $L>\epsilon> 0$, 
\begin{equation*}
 \int_\epsilon^L \frac{\log t }{(1+t +\rho)^2} dt  = 
 - \left\{  \frac{ \log\big(1 + (1+\rho)/t \big) }{1+\rho} + \frac{\log t}{1+ \rho +t}  \right\}_{t=\epsilon}^L \ .
\end{equation*} 

Hence, taking $L\to\infty$, 
\begin{equation*}
 \int_\epsilon^\infty \frac{\log t }{(1+t +\rho)^2} dt  = 
 \frac{(1+\rho) \log(1 +\rho + \epsilon) +  \epsilon\log\big(1+ (1+\rho)/\epsilon \big) }{(1+\rho)(1+ \rho +\epsilon)} 
\end{equation*} 
and, taking $\epsilon \to 0$, this gives formula (\ref{PolyLog_integral_3}) and then (\ref{PolyLog_integral_2}) follows.
  \qed\\

 By definition \ref{C_Poisson}, proposition~\ref{thm:MNS_property} means that the Fermi statistics $\psi(t)=1/(1+e^t)$ is the only shape in $\F$  for which the Poisson component of the field $\Xi_\psi$ is a Gaussian process.
In general, formula (\ref{B_GF}) implies that {\it Laplace transform} of the Poisson component of the random process $\Xi_\Psi$  is given by
\begin{equation} \label{Laplace_B}  
\exp\left(\sum_{n=1}^\infty \frac{\xi^n}{n!}  \B^n_\Psi \int_{\R} f(t)^n dt \right)
 = \exp\left( \int_{\R} \int_0^1  \frac{\Psi^{-1}(t)}{(e^{\xi f(x)}-1)^{-1}+ t } dt dx \right) \ ,
 \end{equation}
 for all $f\in C_0(\R)$ and all $|\xi|< e^{-1}/\|f\|_\infty$.
For a given shape $\Psi\in\F$, it seems very difficult to check whether the RHS of formula (\ref{Laplace_B}) defines a positive definite function in the variable $\xi$ so that the Poisson component of the field $\Xi_\Psi$ comes from a random process.

\subsection{The third and fourth cumulants} \label{sect:C3C4}

Propositions~\ref{thm:MNS_property} and~\ref{thm:Gaussian} imply that for any modified Ensemble whose shape $\Psi \neq \psi$, the limiting fluctuations at the critical scale are not Gaussian. 
The goal of this section is to prove that this is also the case for the  MNS Ensemble.  
Our first attempt is to compute the third cumulant of  the random variable $\Xi_\psi f $, but it turns out that it vanishes for any test function; see proposition~\ref{thm:C3}. Consequently, we construct a test function $y\in \mathcal S(\R)$ such that $\Cu^4\big[\Xi_\psi y\big]\neq 0$.
The strategy to simplify formula  (\ref{GG}) for the cumulants of the random variables $\Xi_\Psi f$ is to symmetrize the functions $\G^{\bf m}$ with respect to all permutations of the variables $u_i$ and to look for cancellations. There are even more simplifications available using the constraints $u_1+\cdots+u_n=0$ and the DHK formulae; see remark \ref{rk:DHK} below. 
 However, as we emphasized in the introduction, it turns out there is no counterpart of Soshnikov's main combinatorial lemma for the modified ensembles and already for the $4^{\text{th}}$ cumulant, it becomes quite technical to rewrite formula (\ref{GG}) in a simple way. In order to get even more simplifications, we shall only consider the following subclass of shapes which includes the MNS shape  $\psi(t)=1/(1+e^t)$.

\begin{definition} \label{shape_symmetry}  
A {\bf shape} $\Psi\in\F$ is called {\bf symmetric} if its derivative $-\Phi$ is even. 
In other words, if it satisfies for all $t\in\R$,
$$1-\Psi(-t)=\Psi(t) \ .$$
\end{definition}

We can deduce from proposition~\ref{thm:b_GF} that the triangular array $b_k^n$ satisfies for any $k=0,\dots,n-1$,
\begin{equation}\label{b_symmetry} b^n_k= (-1)^{n+1} b^n_{n-1-k} \ . \end{equation}

This implies that for any symmetric shape  the map
$$ x\mapsto  \Phi(x)\sum_{k=0}^{n-1}   b_k^n  \Psi(x)^k\big(1-\Psi(x)\big)^{n-1-k} $$
is even when the index $n$ is odd. Thus, for any $m\ge 1$,
\begin{equation}\label{B_symmetry}  
 \B^{2m+1}_\Psi= \int_\R x \Phi(x)  \sum_{k=0}^{n-1}   b_k^n  \Psi(x)^k\big(1-\Psi(x)\big)^{n-1-k} dx =0 \ .
\end{equation}

Using symmetries, we can obtain a simple formula for the $3^{\text{rd}}$ cumulant  of  any random variable~$\Xi_\Psi f $.


\begin{lemma} \label{thm:G3} 
We define the  function $\varpi:\R^2 \to \R^2$ by
\begin{equation*} \varpi(v_1,v_2)=[v_1]^++[v_2]^+ +[v_1+v_2]^+-2\max\{0,v_1,v_1+v_2\} \ .  \end{equation*}
By formula $(\ref{C_split})$, for any $\Psi\in\F$ and any function $f\in H^1_0(\R)$,
\begin{align} \label{C3}
\Cu^3\big[\Xi_\Psi f \big] &=  2\B^3_\Psi \int f(x)^3 dx \\
&\notag
+   24\int_{\R^3_0}d^2u\ \Re\left\{ \prod_i \hat{f}(u_i) \right\} \int_\R ds\ \Phi(s) \underset{(0,\infty)^2}{\iint d^2z}\ \Phi(s+z_1)\Phi(s+z_1+z_2) \varpi(u_1-z_1,u_2-z_2) \ ,
\end{align}
where $d^2u=du_1du_2$ and it is understood that $u_3=-u_1-u_2$.

\end{lemma}
\proof Appendix \ref{A:G}.\qed\\

The important feature of formula (\ref{C3}) is that the functions $\hat f$ and $\Phi$ are coupled by a function $\varpi$ which only depends on the variables  $u_i-z_i$. Moreover, it follows from the DHK formula,  (\ref{DHK_2}), that this function is anti-symmetric. Consequently, if the shape is symmetric, the  $3^{\text{rd}}$ cumulant vanishes for any test function.

\begin{proposition}\label{thm:C3}
If  the shape $\Psi \in \F$ is symmetric, then for any function $f\in H^1_0(\R)$ we have $$\Cu^3\big[\Xi_\Psi f \big]=0 \ .$$
\end{proposition}

\proof First, observe that by  (\ref{B_symmetry}), the constant $\B^3_\Psi=0$  and it remains to show that the second term in formula (\ref{C3}) vanishes as well.
For any $w\in\R^2$, we easily check that
\begin{equation}\label{DHK_2}\max\{0,v_1,v_1+v_2\}+\max\{0,v_2,v_1+v_2\}= [v_1]^++[v_2]^+ +[v_1+v_2]^+  \ .\end{equation}
This implies that another expression for the function $\varpi$ is given by
\begin{equation*} \varpi(v_1,v_2)= \max\{0,v_2,v_1+v_2\} -\max\{0,v_1,v_1+v_2\} \end{equation*}
In particular $\varpi(v_1,v_2)=-\varpi(v_2,v_1)$ and it follows that

\begin{align}  \label{C3'}
&\iiint\limits_{\R\times(0,\infty)^2} \Phi(x)\Phi(x+z_1)\Phi(x+z_1+z_2) \varpi(u_1-z_1,u_2-z_2)\  dxd^2z \\
&\notag\hspace{1cm}
=\iiint\limits_{\R\times(0,\infty)^2} \Phi(-y-z_1-z_2)\Phi(-y-z_2)\Phi(-y) \varpi(u_1-z_1,u_2-z_2)\ dyd^2z \\
&\notag\hspace{1cm}
=\iiint\limits_{\R\times(0,\infty)^2} \Phi(y)\Phi(y+z_1)\Phi(y+z_1+z_2)\varpi(u_1-z_2,u_2-z_1)\ dyd^2z\\
&\label{C3''}\hspace{1cm} 
=-\iiint\limits_{\R\times(0,\infty)^2} \Phi(y)\Phi(y+z_1)\Phi(y+z_1+z_2)\varpi(u_2-z_1,u_1-z_2)\ dyd^2z  \ .
\end{align}
At first, we made the  change of variable $y=-x-z_1-z_2$. In the second equality we used the assumption that $\Phi$ is symmetric and permuted the variables $z_1$ and $z_2$. For the last equality,  we used the anti-symmetry of the function $\varpi$.
Equation~(\ref{C3''}) shows that the integral (\ref{C3'}) changes sign under permutation of the variables $u_1$ and $u_2$. Because of this fact and the symmetry of formula (\ref{C3}), the $3^{\text{rd}}$ cumulant of the random variable $\Xi_\Psi f$ vanishes. \qed\\

\begin{remark}  \label{rk:DHK}
The Dyson, Hunt, Kac $(DHK)$ formulae are the following remarkable identities. For any $n \ge 2$,
\begin{equation} \label{DHK} 
\sum_{\pi \in \Sy(n)} \max\big\{ u_{\pi(1)}, u_{\pi(1)}+u_{\pi(2)},\dots,u_{\pi(1)}+\cdots+u_{\pi(n-1)},0\big\}
 = \sum_{\pi \in \Sy(n)} \sum_{l=1}^n \frac{1}{l} \big[ u_{\pi(1)}+\cdots+u_{\pi(l)} \big]^+ \ .
\end{equation}
When $n=2$, this gives formula $(\ref{DHK_2})$. We refer to  Simon's book $\cite{Simon_04a}$ section $6.5$ for a proof of $(\ref{DHK})$ and an application to the Strong Szeg\H{o} theorem. Actually, the proof of Sosnhikov's main Combinatorial lemma $(\ref{MCL})$ is also based on these formulae, see \cite[Appendix~A]{L_15b}.
We can also apply the formulae $(\ref{DHK})$ to the cumulants of linear statistics of the modified ensembles but, except for the third cumulant,  this only leads to partial simplifications. \end{remark}

To compute the $4^{\text{th}}$ cumulant, we  make the change of variables $x\in \R^4_< \mapsto (s, z)\in\R\times\R_+^3$ given by 
\begin{equation*}
s=x_1
\hspace{1.5cm} z_1=x_2- x_1
\hspace{1.5cm} z_2=x_3- x_2
\hspace{1.5cm} z_3=x_4- x_3 
\end{equation*}
in equation (\ref{GG}). We get
\begin{equation} \label{C4}
\mathfrak G^4_\psi (f)=- 2\underset{(0,\infty)^3}{\int d^3z}\ \Theta(z)  \underset{\R^4_0}{\int }d^3u\ \Re\left\{ \prod_{i=1}^4 \hat{f}(u_i) \right\}  \sum_{|{\bf m}|=4} \M({\bf m}) \tilde \G^{\bf m}(u,z)   \ ,
\end{equation}
where $\tilde \G^{\bf m}(u,z)$ is the image of $x\mapsto \G^{\bf m}(u,x)$ under the change of variables  (this function does not depend on the variable $s$) and 
 \begin{equation} \label{Theta}
  \Theta(z)=  \int_\R \phi(s)\phi(s+z_1)\phi(s+z_1+z_2)\phi(s+z_1+z_2+z_3) ds \ . \\
  \end{equation}

 It is worth noting that, since the function $\phi$ is even, we have $\Theta(z_3,z_2,z_1)=\Theta(z_1,z_2,z_3) $ but no further symmetry. 
 We are not able to obtain a compact formula for the $4^{\text{th}}$ cumulant and it turns out to be simpler  to compute the value of the functions $\tilde\G^{\bf m}(u,z)$ at some well-chosen points and deduce from formula (\ref{C4}) that $\mathfrak {G}^4_\psi(y) \neq 0$ when the test function $y$ which is sufficiently concentrated around these points. The technical result that we need is given in the following lemma.  
 
  \begin{lemma} \label{thm:G4} For any $z\in\R_+^3$, up to the permutation of $z_1$ and $z_3$,  we have
 \begin{align*} &\sum_{\begin{subarray}{c} v_1+\cdots+v_4=0\\ v_i\in\{-1,1\} \end{subarray} } 
  \sum_{|{\bf m}|=4} \M({\bf m})\tilde \G^{\bf m}(v,z)
= 24\bigg(4[1-z_2]^+ +2[1-z_1-z_2]^+      + \frac{[2-z_2 ]^+}{2} + \frac{ [2-z_1-z_2 -z_3]^+}{2}\\
& \hspace{.5cm}+ [2-z_1-z_2 ]^+ -[2-z_1]^+  -2 \max\{0,1-z_1,2-z_1-z_2\}  -2\max\{0,1-z_1, 2-z_1-z_2-z_3\} \bigg) \ . \end{align*}
  \end{lemma}
 \proof Appendix \ref{A:G}.\qed\\
 
 It follows from lemma~\ref{thm:G4} that
 
 \begin{equation} \label{G4_Mathematica}
 \sum_{\begin{subarray}{c} v_1+\cdots+v_4=0\\ v_i\in\{-1,1\} \end{subarray} } \int_{(0,\infty)^3} \Theta(z)  \sum_{|{\bf m}|=4} \M({\bf m})\tilde \G^{\bf m}(v,z)\  d^3z =  0.29... \ .
\end{equation} 
The integral can be performed analytically or numerically using {\it Mathematica}. To complete our argument, we also need the following approximation lemma. Its proof is rather straightforward and for completeness it will be given after our example.

\begin{lemma}\label{approximation_lemma} Let $g(x)=e^{- \pi x^2}$ and  
$\displaystyle \A(u)= \int_{(0,\infty)^3} \Theta(z)  \sum_{|{\bf m}|=4} \M({\bf m})\tilde \G^{\bf m}(u,z) d^3z$.
For any $ v\in\R^4$, we have
\begin{equation*}\label{Gaussian_approximation}\lim_{\epsilon\to0}\ \epsilon^{-4} \int_{\R^4_0} g\left(\frac{u_1-v_1}{\epsilon}\right)\cdots g\left(\frac{u_4-v_4}{\epsilon}\right) \A(u)\ d^3u =\frac{1}{2} \A(v)\ \delta_0(v_1+\cdots+v_4) \ .
\end{equation*}\\
 \end{lemma}

 We let $y(x)= 2e^{-\epsilon \pi x^2}\cos(2\pi x)$ for some $\epsilon>0$. Then
\begin{equation*}
\hat{y}(u)=\epsilon^{-1}  g\left(\frac{u-1}{\epsilon}\right)+\epsilon^{-1} g\left(\frac{u+1}{\epsilon}\right) \ ,
\end{equation*}
and by lemma \ref{approximation_lemma}, we have
\begin{equation} \label{counterexample}
  \lim_{\epsilon\to0}  \underset{\R^4_0}{\int }d^3u\ \Re\left\{ \prod_{i=1}^4 \hat{y}(u_i) \right\}\underset{(0,\infty)^3}{\int d^3z}\ \Theta(z)  \sum_{|{\bf m}|=4} \M({\bf m}) \tilde \G^{\bf m}(u,z) 
= \frac{1}{2} \sum_{\begin{subarray}{c} v_1+\cdots+v_4=0\\ v_i\in\{-1,1\} \end{subarray} }  \A(v) \ .
\end{equation}
The RHS of equation (\ref{counterexample}) is given by  (\ref{G4_Mathematica}) and by equation (\ref{C4})
$$ \lim_{\epsilon\to0} \mathfrak G_\psi (y) = -0.29... \ . $$

Since the constant $\B^4_\psi=0$ by proposition~\ref{thm:MNS_property}, 
we conclude that, if the parameter $\epsilon$ is sufficiently small, $\Cu^4\big[\Xi_\psi y\big] \neq 0$ and the linear statistics $\Xi_\psi y$ is not Gaussian. \\

{\it Proof of lemma \ref{approximation_lemma}.} Let us fix $v\in \R^4$ and let $r(v)= v_1+\cdots+v_4$.
It is easy to see that  the functions $u\mapsto \tilde{G}^{\bf m}(u,z)$ are Lipchitz continuous with respect to $|u|_\infty$ with some constant which can be chosen independently of $z\in\R_+^3$. Then, the function $\A(u)$ is also Lipschitz continuous on $\R^3$. A change of variables yields 
\begin{align*} \epsilon^{-4}\int_{\R^4_0} g\left(\frac{u_1-v_1}{\epsilon}\right)\cdots g\left(\frac{u_4-v_4}{\epsilon}\right) \A(u)\  d^3u 
&= \int_{\R^4_{-r(v)/\epsilon}} g(w_1)\cdots g(w_4) \A(v+\epsilon w)\  d^3w\\
&= \A(v) \int_{\R^4_{-r(v)/\epsilon}} g(w_1)\cdots g(w_4)\  d^3w + \underset{\epsilon\to0}{O}(\epsilon).
 \end{align*}

If we let $X_{\epsilon} \in \No(\epsilon^{-1}r(v), \frac{3}{2\pi})$, it is easy to see that
\begin{equation} \label{Gaussian_integral_4} \int_{\R^4_{-r(v)/\epsilon}} g(w_1)\cdots g(w_4)\  d^3w =   \mathbb{E}\left[e^{-\pi X_{\epsilon}^2}\right]   . \end{equation}

So that if  $v_1+\cdots+v_4=0$, i.e. $r(v)=0$, then $ \mathbb{E}\left[e^{-\pi X_{\epsilon}^2}\right]= \frac{1}{2}  $ for any $\epsilon>0$ and it follows that
\begin{equation*} \lim_{\epsilon\to0}\epsilon^{-4}\int_{\R^4_0} g\left(\frac{u_1-v_1}{\epsilon}\right)\cdots g\left(\frac{u_4-v_4}{\epsilon}\right) \A(u)\  d^4u = \frac{1}{2}\A(v).  \end{equation*}
 On the other hand if $r(v) \neq 0 $, by equation $(\ref{Gaussian_integral_4})$,
 \begin{equation*}\lim_{\epsilon\to0} \int_{\R^4_{r(v)/\epsilon}} g(w_1)\cdots g(w_4)\  d^4w=0.  \end{equation*}
 \qed 

\section{Appendices} \appendix

\section{Asymptotics of the Hermite polynomials} \label{sect:Asymptotic}

In this section, we provide some background on the asymptotics of the Hermite polynomials and the GUE kernel, see for instance~\cite[Proposition~5.1.3]{Pastur_Shcherbina}
or \cite[Section~2.2]{Deift_al_99_a}. 
These asymptotics are named after Plancherel-Rotach, \cite{PR_29},  and can be derived using the classical saddle point method.
To investigate the statistics of Hermitian invariant ensembles, one is usually interested in uniform asymptotics of the Christoffel-Darboux kernels. Based on the Riemann-Hilbert problem for orthogonal polynomials, the sine-kernel asymptotics have been established  in \cite{Deift_al_99_b} at the microscopic scale for a large class of potentials. These results have recently been extended to mesoscopic scales  in  \cite{KSSV_14}, see also \cite{L_15a}, and we will present the results for the GUE kernel. 
For all $n\in\N$, let $h_n$ be the normalized Hermite polynomial of degree $n$ with respect to $e^{-x^2}$ on $\R$ and define the Hermite functions:
\begin{equation}  \label{Hermite_wave_function}  
\vartheta_n(x)=  h_n(x)e^{-\frac{x^2}{2}} .\end{equation}

\begin{definition}  \label{F} 
 We define on $[-1,1]$ the functions $\varrho(t)=2\sqrt{1-t^2}$ and 
\begin{equation*}
F(x)=\int_0^x \varrho(t) dt = \arcsin(x)+x\sqrt{1-x^2} \ .
 \end{equation*}  
 The map $F$ is a  diffeomorphism  from $|x|<1$ to $|x|<\frac{\pi}{2}$ and we let $G$ be its inverse.
\end{definition}


The Hermite polynomials have the following bulk asymptotics:
\begin{equation} \label{Hermite_asymptotic_bulk}
 \vartheta_n(\sqrt{2n}x)=\left(\frac{2}{n(1-x^2)}\right)^{1/4}\frac{1}{\sqrt{\pi}}\left\{\cos \left(n\frac{\pi}{2}- nF(x)-\frac{1}{2}\arcsin(x)\right)+\underset{n\to\infty}{O}\left(\beta^{-3/2}n^{-3\gamma/2}\right)\right\} \ ,    
\end{equation}
for all  $|x|\le 1-\beta n^{-2/3+\gamma}$ where $0< \gamma\le 2/3$ and $0<\beta<1$.  Observe that, in the case $\gamma=2/3$, we obtain an asymptotics valid in any fixed interval $[-1+\beta,1-\beta]$  with a uniform error term of order~$N^{-1}$.   
While the Hermite functions have oscillatory behavior inside the bulk, they have exponential decay outside:
\begin{equation}  
\vartheta_n(\sqrt{2n}x)= \frac{e^{- n H(x)}\big\{ 1  +O(n^{-3\gamma/2}) \big\}}{\sqrt{\pi\sqrt{2n} \sqrt{x^2-1}}}  \ , 
\end{equation}
for all $|x| \ge 1+ n^{-2/3+\gamma}$ where $0<\gamma<2/3$ and the $H$ is even  and defined for any $x>1$ by
\begin{equation*} H(x)=\int_1^x 2\sqrt{t^2-1} dt =   x\sqrt{x^2-1}-\log(x+\sqrt{x^2-1})   . 
 \end{equation*}
In particular, since $\frac{4\sqrt{2}}{3} (x-1)^{3/2} \le H(x) $, there exist a constants $C>0$ such that  for any $x > 1 + n^{-1/2} $, 
\begin{equation} \label{Hermite_asymptotic_exterior}
\big| \vartheta_n(\sqrt{2n}x)  \big| \le C e^{- n  \frac{4\sqrt{2}}{3} (x-1)^{3/2}}  . 
\end{equation}

 At the edge, the asymptotics is also well-known but, since we are only interested in bulk  linear statistics in this paper,  we won't need any precise estimates  and we will use instead the uniform bound 
\begin{equation}  \label{Hermite_bound_uniform}  \| \vartheta_n\|_\infty \le c_H\  n^{-1/12} \ ,    \end{equation}
where $c_H$ is a universal constant and the exponent is sharp. 
We define the Christoffel-Darboux  kernel
\begin{equation} \label{CD_kernel}
K_{\CD}^N(x,y)=\sum_{n=0}^{N-1} \vartheta_n(x)\vartheta_n(y)
= \sqrt{\frac{N}{2}}\frac{\vartheta_N(x)\vartheta_{N-1}(y)-\vartheta_{N-1}(x)\vartheta_N(y)}{x-y} \ ,
  \end{equation}
  and the Wigner semicircle law, for any $|t| \le \sqrt{2}$,
  \begin{equation} \label{SC}
  \varrho_{sc}(t) =  \frac{1}{\pi} \sqrt{2-t^2}  = \frac{1}{\pi\sqrt{2}}\varrho\left(\frac{t}{\sqrt{2}}\right) \ .
  \end{equation}

 At the microscopic scale, it is well-known that we get the sine kernel in the limit for any $|x_0|<\sqrt{2}$, 
\begin{equation*} \lim_{N\to\infty} \frac{\pi}{\sqrt{2N}}K^N_{\CD}\left(\sqrt{N}x_0+ \frac{\pi \xi}{\sqrt{2N}},\sqrt{N}x_0+\frac{\pi \zeta}{\sqrt{2N}} \right)= \frac{\sin{[\pi \varrho_{sc}(x_0)(\xi-\zeta)]}}{\pi(\xi-\zeta)} \ .
\end{equation*}
Using the results of \cite{Deift_al_99_b}, this asymptotics can be extended to all mesoscopic scales.  In the sequel, $L$ is some arbitrary large  positive constant and $|x_0|<\sqrt{2}$.

\begin{theorem}\label{thm:meso_sine_kernel} For any $-1/2\le \lambda <1/2$, we have the asymptotic  formula
\begin{align*} 
N^\lambda K_{\CD}^N\left(\sqrt{N}x_0+\xi N^{\lambda},\sqrt{N}x_0 +\zeta N^{\lambda}\right)   
&= \frac{\sin \left[ N \left( F(\frac{x_0+ N^{-1/2+\lambda}\xi}{\sqrt{2}})-F(\frac{x_0+N^{-1/2+\lambda}\zeta}{\sqrt{2}})\right)\right]}{\pi(\xi-\zeta)}\\
&\hspace{.3cm}+ \underset{N\to\infty}{O}\left(N^{- 1/2+\lambda}\right) 
\end{align*}
uniformly over all $\xi,\zeta \in[-L,L]$. 
\end{theorem}

Theorem~\ref{thm:meso_sine_kernel} was first proved in~\cite{KSSV_14} using the Riemann-Hilbert formulation of~\cite{Deift_al_99_b}. In~\cite{L_15a}, we produce an elementary proof which is based on the classical steepest descent method performed in the seminal paper~\cite{PR_29}. 
Note that this approximation  takes into account the density of the Wigner semicircle law, i.e.~the fact that the GUE eigenvalues are not uniformly   distributed at the global scale. 
Namely, by definition \ref{F} and (\ref{SC}), theorem~\ref{thm:meso_sine_kernel} can be rephrased as
\begin{align} \label{meso_sine_kernel}
N^\lambda K_{\CD}^N\left(\sqrt{N}x_0+\xi N^{\lambda},\sqrt{N}x_0 +\zeta N^{\lambda}\right)    
&= \frac{\sin\left[ \pi  N^{1/2+\lambda}  \displaystyle\int_{\zeta}^{\xi} \varrho_{sc}\big(x_0+tN^{-1/2+\lambda}\big) dt   \right]}{\pi (\xi - \zeta)} \\
&\hspace{.3cm}\notag+ \underset{N\to\infty}{O}\left(N^{- 1/2+\lambda}\right) \ . 
\end{align}
Note that unlike the sine-kernel, the kernel (\ref{meso_sine_kernel}) is not translation-invariant. This raises complications to compute the limits of the cumulants for large scale linear statistics of the modified GUEs, cf.~proposition~\ref{thm:MNS_trace}.
However, at sufficiently small scales, we recover the sine-kernel as a special case of theorem~\ref{thm:meso_sine_kernel}.

\begin{theorem} \label{thm:micro_sine_kernel}  For any $-1/2\le \lambda <0$, we have the asymptotic formula
\begin{equation*} \label{micro_sine_kernel} 
N^\lambda K_{\CD}^N\left(\sqrt{N}x_0+\xi N^{\lambda},\sqrt{N}x_0 +\zeta N^{\lambda}\right)  
= \frac{\sin\left(N^{1/2+\lambda}\pi \varrho_{sc}(x_0)(\xi-\zeta)\right)}{\pi (\xi-\zeta)} 
+\underset{N\to\infty}{O}\left(N^{2\lambda}\right) 
\end{equation*}
uniformly over all $\xi,\zeta \in[-L,L]$.\\
 \end{theorem}
 
 \proof By definition \ref{F}, 
 \begin{equation*} 
 F(x)-F(y) = \varrho\left(\frac{x+y}{2} \right) (x-y) + O\left(|x-y|^3\right) \ .
\end{equation*}
If we let $x_N=\frac{x_0+ N^{-1/2+\lambda}\xi}{\sqrt{2}}$ and $y_N=\frac{x_0+ N^{-1/2+\lambda}\zeta}{\sqrt{2}}$, a Taylor expansion gives
\begin{equation} \label{SC_expansion}
\varrho\left(\frac{x_N+y_N}{2} \right) = \varrho\left(\frac{x_0}{\sqrt{2}}\right) + O\big(N^{-1/2+\lambda}\big) \ , 
\end{equation}
and it follows that
 \begin{equation*} 
 F(x_N)-F(y_N) = \frac{1}{\sqrt{2}}\varrho\left(\frac{x_0}{\sqrt{2}}\right) (\xi-\zeta)N^{-1/2+\lambda} + O\big(|\xi-\zeta|N^{-1+2\lambda}\big) \ .
\end{equation*}
Hence, by (\ref{SC}), we have proved that
$$ \sin\left[ N \left( F(\frac{x_0+ N^{-1/2+\lambda}\xi}{\sqrt{2}})-F(\frac{x_0+N^{-1/2+\lambda}\zeta}{\sqrt{2}}) \right) \right]
=  \sin\left[N^{1/2+\lambda}\pi\varrho_{sc}(x_0)(\xi-\zeta) \right] + \underset{N\to\infty}{O}\big(|\xi-\zeta|N^{2\lambda}\big) \ .
$$ 
When $\lambda<0$, the error term is converging to 0 for any $x,y \in [-L,L]$, and if we plug this approximation in the formula of theorem~\ref{thm:meso_sine_kernel}, we obtain the asymptotics of theorem \ref{thm:micro_sine_kernel}. \qed

\begin{remark} \label{rk:x0} 
In the special case $x_0=0$, since $\varrho'(0)=0$, the error term in $(\ref{SC_expansion})$ is of order $N^{-1+2\lambda}$
and the sine kernel approximation of theorem~$\ref{thm:micro_sine_kernel}$ is valid in the whole range $-1/2 \le \lambda<1/6$.
\end{remark}

The Christoffel-Darboux kernel $K_{\CD}^N$ is the same, up to a scaling, as the GUE kernel $K_0^N$ defined by (\ref{GUE_kernel}) in the introduction. Namely, if we let $\lambda=\frac{1}{2}-\delta$ for $\delta\in(0,1]$, by (\ref{CD_kernel}), we can rewrite
\begin{equation*} \label{GUE_CD} N^{-\delta}K_0^M(xN^{-\delta},yN^{-\delta}) 
= \frac{\pi}{\sqrt{2}} N^{\lambda} K_{\CD}^M\left( \frac{\pi}{\sqrt{2}}x N^{\lambda},\frac{\pi}{\sqrt{2}}y N^{\lambda}\right)
\ .\end{equation*}
Thus, for any $\epsilon>0$, if $x_0=0$ and $M=N+k$ for some $|k| \le N^{1-\epsilon}$, by theorem ~\ref{thm:meso_sine_kernel}, for all $x, y\in[-L,L]$, we have 
\begin{equation}  \label{sine_approximation}
N^{-\delta}K_0^M(xN^{-\delta},yN^{-\delta})
 = \frac{ \sin\left[ M \left( F\left(\frac{\pi}{2} \frac{x}{\sqrt{M}N^{\delta-1/2}} \right)- F\left(\frac{\pi}{2} \frac{y}{\sqrt{M}N^{\delta-1/2}} \right) \right)\right] }{\pi(x-y)} +\underset{N\to\infty}{O}(N^{-\delta}) \ .
\end{equation}

This formula holds at any mesoscopic scales. On the other hand, if we assume that $1/3<\delta\le 1$ (cf$.$ remark~\ref{rk:x0}),  since $\varrho_{sc}(0)=\frac{\sqrt 2}{\pi}$, by theorem~\ref{thm:micro_sine_kernel}, we get the following asymptotics 

\begin{equation}  \label{sine_approximation_2}
N^{-\delta}K_0^M(xN^{-\delta},yN^{-\delta}) =  \frac{\sin\big[ \pi N^{1/2-\delta}\sqrt{M} (x-y)\big]}{\pi(x-y)} + \underset{N\to\infty}{O} \left( N^{1-3\delta} \right) \ .
\end{equation}

\section{Proof of theorem~\ref{thm:modified_variance}} \label{A:variance}
 
We give two different proofs of  theorem~\ref{thm:modified_variance}. First, we can use the ideas of section~\ref{sect:Variance} to compute the limit of the reproducing variance $V_0(f_\delta)$ at the critical scale $\delta=\alpha$, see proposition \ref{thm:critical_Variance}.  Since this result is not used in the rest of the paper, we will only sketch the argument. Second, we use  (\ref{cumulant_convergence}) and compute  $\Cu^2\big[\Xi_{\Psi,\tau/4} f\big]$ using formula (\ref{MNS_C^n}). Subsequently, we check that the different formulae for the critical variance are consistent and we apply them to the MNS ensemble.

\begin{proposition}  \label{thm:critical_Variance}
For any shape $\Psi\in\F^*$, any test function $f\in H^{1/2}_0(\R)$ and any scale $0<\alpha<1$, we have
\begin{equation*}\label{critical_variance_2} 
\lim_{N\to\infty} V_0(f_\alpha)  =\frac{1}{4\pi^2}  \iint \left| \frac{f(x)-f(y)}{x-y}\right|^2
\left|\hat\Phi\left(\frac{x-y}{4 \tau^{-1}} \right)\right|^2  dxdy \ .
\end{equation*} 
\end{proposition}

\proof By formula (\ref{variance_reproducing}) combined to the estimate (\ref{variance_decay}) and the approximation (\ref{kernel_approximation_0}), we see that 
\begin{align*} \label{critical_Variance_1}
&V_0(f_\alpha) 
= \frac{1}{2\pi^2\tau^2 N^{2\alpha}} \underset{\begin{subarray}{c}   |k| \le \Gamma_N N^\alpha \\  |j| \le \Gamma_N N^\alpha \end{subarray}}{\sum}  \Phi\left(\frac{k}{\tau N^\alpha}\right)  \Phi\left(\frac{j}{\tau N^\alpha}\right)
 \underset{[-L,L]^2}{\iint} dxdy \left| \frac{f(x)-f(y)}{x-y}\right|^2  \times \\
& \notag  \sin\left[ \sqrt{N(N+k)} \left( F\left(\frac{\pi}{2} \frac{x}{N^\alpha} \right)- F\left(\frac{\pi}{2} \frac{y}{N^\alpha} \right) \right) \right]  \sin\left[ \sqrt{N(N+j)} \left( F\left(\frac{\pi}{2} \frac{x}{N^\alpha} \right)- F\left(\frac{\pi}{2} \frac{y}{N^\alpha} \right) \right) \right]   
 + \underset{L\to\infty}{O}(L^{-1}) \ ,
\end{align*}
where the support of the test function $f$ is included in $[-\frac{L}{2},\frac{L}{2}]$ and  $\Gamma_N=(\log N)^2$. Note that there is another error-term of order $N^{-\alpha}$ coming from (\ref{kernel_approximation_0}) that we dropped. 
Applying a trigonometric identity and the Riemann-Lebesgue lemma like in the proof of lemma~\ref{thm:variance_sigma}, we see that
\begin{align*}
&V_0(f_\alpha)
\simeq \frac{1}{4\pi^2}  \frac{1}{\tau^2 N^{2\alpha}} \underset{\begin{subarray}{c}   |k| \le \Gamma_N N^\alpha \\  |j| \le \Gamma_N N^\alpha \end{subarray}}{\sum}  \Phi\left(\frac{k}{\tau N^\alpha}\right)  \Phi\left(\frac{j}{\tau N^\alpha}\right)\times\\
&\hspace{3cm}\iint\limits_{[-L,L]^2}  \left| \frac{f(x)-f(y)}{x-y}\right|^2 
\cos\left( \frac{(k-j)\pi(x-y)}{2N^{\alpha}} \right) dxdy + \underset{L\to\infty}{O}(L^{-1}) 
\end{align*}
as $N\to\infty$.
The sums converge to some Riemann integrals and by the dominated convergence theorem, 
\begin{equation} \label{Phi_Fourier}
\limsup{N\to\infty} \left| V_0(f_\alpha)  -\frac{1}{4\pi^2}  \iint\limits_{[-L,L]^2}  \left| \frac{f(x)-f(y)}{x-y}\right|^2 \iint\limits_{\R^2} \Phi(u) \Phi(v)  \cos\left[ \frac{ \pi(u-v)(x-y)}{2 \tau^{-1}} \right] dudv dxdy \right| \le \frac{C(f)}{L} \ ,
\end{equation}
where the  constant $C(f)$ depends only on the test function $f$ by (\ref{variance_decay}).
To complete the proof, we can use that  for any $w\in\R$, 
\begin{equation} \label{variance_3}
  \iint_{\R^2} \Phi(u) \Phi(v)  \cos\left[ 2\pi(u-v)w \right] dudv =\big|\hat\Phi(w)\big|^2 \ ,
  \end{equation} 
  and  let~$L\to\infty$ in the inequality 
(\ref{Phi_Fourier}). \qed\\

As a consequence of lemma~\ref{thm:variance}, lemma~\ref{thm:variance_sigma} and proposition~\ref{thm:critical_Variance}, we obtain that
for any $0<\alpha<1$ and  $f\in H^{1/2}_0\cap L^\infty(\R)$,
\begin{equation*}
\lim_{N\to\infty}\Var_{K^N_{\Psi,\alpha}}[f_\alpha] = \frac{\tau}{2} \int_\R f(x)^2 dx  
+\frac{1}{4\pi^2}  \iint_{\R^2} \left| \frac{f(x)-f(y)}{x-y}\right|^2\left|\hat\Phi\left(\frac{(x-y)}{4 \tau^{-1}} \right)\right|^2  dxdy \ .
\end{equation*}

By theorem~\ref{thm:MNS_C^n}, this implies that for any $\Psi \in\F^*$ such that $\B^2_\Psi =1$ and $\tau>0$, 
\begin{equation}\label{variance_2} 
\Var[\Xi_{\Psi, \tau} f] = 2\tau \int_\R f(x)^2 dx  
+\frac{1}{4\pi^2}  \iint_{\R^2} \left| \frac{f(x)-f(y)}{x-y}\right|^2\big|\hat\Phi\big(\tau (x-y)\big)\big|^2  dxdy \ .
\end{equation}

\begin{remark} Such a direct computation seems possible only for the variance of linear statistics thanks to the special structure of the reproducing variance $V_0$, $(\ref{variance_reproducing})$. A technical difficulty to compute the limit of the higher-order cumulants  comes from the singularity of the correlation kernel $K_{\Psi,\alpha}^N$ along the diagonal. Therefore, it is better to exploit instead the fact that the kernel $L^N_{\Psi,\eta}$ given by $(\ref{kernel_L})$ is translation-invariant and use Soshnikov's method.
\end{remark}

We have seen in section~\ref{sect:C^n} that, up to a scaling, the modified ensembles have the same limit at the critical scale. Thus we have obtained three expressions for the variance of the random variable $\Xi_{\Psi,\tau} f$;  formula \eqref{MNS_C^n} and  formulae (\ref{modified_CUE_variance}) and (\ref{variance_2})  in the Circular and Gaussian case respectivel. Note that these formulae are well-defined for any shape $\Psi \in \F$ and  and we will now check that they are consistent. By \eqref{M}, $\M({\bf 1+1})= -1$ and by  (\ref{Lambda}), for any $u\in\R^2_0$,
$$\Lambda^{ {\bf 1}^2}(u)= \left(\begin{array}{cc} 0 & -u_2 \\ u_2 & 0 
\end{array}\right) \  .$$
Moroever, by definition (\ref{G}), for any $x_1<x_2$,
\begin{align} \notag \G^{{\bf 1}^2}_\tau(u,x)&= \max\left\{0, u_2-\tau(x_2-x_1)\right\}+  \max\left\{ -u_2-\tau(x_2-x_1),0\right\} \ .
 \end{align}

Hence, for any function $f \in H^{1/2}_0\cap L^\infty(\R)$, by formula \eqref{MNS_C^n} with $n=2$, 
\begin{equation} \label{C2} 
\Var\big[\Xi_{\Psi,\tau} f\big] 
= 2\tau  \int_\R f(x)^2 dx  
+4 \int_0^\infty \big|\hat{f}(u)\big|^2\iint\limits_{x_1<x_2} \Phi(x_1)\Phi(x_2)\big[u-\tau(x_2-x_1)\big]^+ dudx_1dx_2 \ .
\end{equation}

We first check that this formula matches with the RHS of (\ref{modified_CUE_variance}).
Using  that $\Phi=-\Psi'$ and the properties of $\Psi$, \eqref{class_F},  some integrations by parts  shows that for any $u>0$, 
\begin{equation*} 
\int_\R \Phi(t) \int_0^\infty \Phi(t+s) [u-s]^+ ds dt =
\frac{u}{2} -  \B^2_\Psi + \int_\R   \Psi(x+u)\big(1-\Psi(x)\big) dx \ .
 \end{equation*}
 
 Since, by convention $\B^2_\Psi=1$,  this implies that 
 \begin{equation} \label{variance_3} 
\Var\big[\Xi_{\Psi,\tau} f\big] 
=   2  \int_0^\infty \big|\hat{f}(u)\big|^2 u du
+4 \int_0^\infty \big|\hat{f}(u)\big|^2 
\int_\R   \Psi\left(\frac{x+u}{\tau}\right)\bigg(1-\Psi\left(\frac{x}{\tau}\right)\bigg) dx du \ .
\end{equation}

Moreover, we can check that for any $u\in\R$, 
\begin{equation*} 
 \int_\R \Psi\left(\frac{x}{\tau} \right) \bigg( 1-\Psi\left(\frac{x+u}{\tau}\right) \bigg) dx 
 =  u +  \int_\R  \Psi\left(\frac{x+u}{\tau}\right)  \bigg( 1-\Psi\left(\frac{x}{\tau}\right) \bigg) dx  
\end{equation*}
so that, according to formula \eqref{variance_3},

 \begin{equation} \label{variance_4} 
\Var\big[\Xi_{\Psi,\tau} f\big] 
=   2  \int_0^\infty \big|\hat{f}(u)\big|^2 
\int_\R  \bigg(  \Psi\left(\frac{x-u}{\tau}\right)+ \Psi\left(\frac{x+u}{\tau}\right) \bigg)\bigg(1-\Psi\left(\frac{x}{\tau}\right)\bigg) dx du \ .
\end{equation}

By \eqref{variance_5}, this establishes that the r.h.s.~of formulae \eqref{C2}
and  (\ref{modified_CUE_variance}) are equals.
It remains to  check that formula (\ref{C2}) also matches with  (\ref{variance_2}).  To do so, we  use an argument which is similar to the proof of the identity (\ref{norm_1/2}) for the variance of the sine process.

\begin{lemma} \label{thm:C2}
For any function $f\in H^{1/2}(\R)$ and any $w>0$, we have
 \begin{equation*} 
  \frac{1}{4\pi^2} \iint_{\R^2}  \left| \frac{f(x)-f(y)}{x-y} \right|^2 \cos\big(2\pi(x-y) w\big)  dxdy 
=  \int_\R \left|\hat{f}(u) \right|^2 \big[|u|-w\big]^+ du \ .
 \end{equation*}
\end{lemma}

\proof
By Plancherel's formula, for any $z\in\R$,
\begin{equation*}\int_\R \left|f(x)-f(x+z)\right|^2 dx = 4 \int_\R |\hat f(u)|^2 \sin^2(\pi u z) du \ .
 \end{equation*}
 Then, by Fubini's theorem, for any $w>0$,
 \begin{align} \notag
  \frac{1}{4\pi^2} \iint  \left| \frac{f(x)-f(y)}{x-y} \right|^2 \cos\big(2\pi(x-y) w\big)  dxdy 
&  = \frac{1}{4\pi^2} \iint  \left| \frac{f(x)-f(x+z)}{z}\right|^2 \cos\big(2\pi z w\big) dxdz\\
  &  \label{Variance_Plancherel_2} 
  = \frac{1}{\pi^2} \int |\hat f(u)|^2 \left( \int \frac{\sin^2(\pi u z)}{z^2} \cos\big(2\pi z w\big) dz\right) du \ .
\end{align}
Moreover, by Residue calculus,  one can show that for any $w>0$ and $u\in\R$,
\begin{equation} \label{V_identity_2}  
\int  \frac{\sin^2(\pi u z)}{z^2} \cos(2\pi w z) dz =  \pi^2 [|u| - w]^+ \ .
\end{equation}
The lemma follows by combining equations (\ref{Variance_Plancherel_2}) and (\ref{V_identity_2}). \qed\\

Lemma~\ref{thm:C2} implies that 
$$
\frac{1}{4\pi^2}  \iint\limits_{[-L,L]^2}  \left| \frac{f(x)-f(y)}{x-y}\right|^2 \iint\limits_{\R^2} \Phi(u) \Phi(v)
$$
\begin{align*}
  \frac{1}{4\pi^2}\iint_{\R^2} \iint_{ t<s}  \left|\frac{f(x)-f(y)}{x-y}\right|^2 &\Phi(t)\Phi(s) \cos\left[2 \pi \tau (t-s)(x-y) \right] dxdydsdt  \\
&= 2\int_0^\infty \left|\hat{f}(u)\right|^2\iint_{t<s} \Phi(t)\Phi(s)\big[u-\tau(s-t)\big]^+dudsdt\ .
\end{align*}
Then, if we take $w=\tau (s-t)$ in formula \eqref{variance_3}, by symmetry we obtain 
\begin{equation*}
\frac{1}{4\pi^2}  \iint_{\R^2} \left| \frac{f(x)-f(y)}{x-y}\right|^2\big|\hat\Phi\big(\tau (x-y)\big)\big|^2  dxdy 
= 4 \int_0^\infty \left|\hat{f}(u)\right|^2\iint_{t<s} \Phi(t)\Phi(s)\big[u-\tau(s-t)\big]^+dudsdt\ .
\end{equation*}
Hence, if we add up the Poisson contribution $ 2\tau\|f\|^2_{L^2}$, we conclude that the 
r.h.s.~of formulae \eqref{variance_2}
and  (\ref{C2}) are equals.
As an example, let us see what these formulae look like for the MNS ensemble.
The MNS shape is  $\psi(t)=(1+ e^t)^{-1}$ and an elementary integration gives for any $u\in\R$, 
\begin{equation*}
 \int_\R \psi(t+u)\big( 1-\psi(t) \big) dt = e^u \int_0^\infty \frac{1}{(1+s)(s+e^u)} ds =  \frac{u}{1-e^{-u}} \ .
\end{equation*}
Then, by formula \eqref{variance_4}, we get
\begin{equation*}  \label{MNS_variance_1}
\Var\big[\Xi_{\psi,\tau} f\big]  =  \int_\R  \left| \hat{f}(u) \right|^2 \frac{u}{\tanh(\frac{u}{2 \tau})} du \ .
  \end{equation*}
We can deduce the dual of this formula using equation (\ref{variance_2}). 
We have
\begin{equation*} 
\phi(t)= \frac{4}{\cosh[t/2] ^2} \ ,
\hspace{1cm} \hat{\phi}(u)= \frac{2\pi^2 u}{\sinh[2\pi^2 u]}  \ ,
\end{equation*}
so that
\begin{equation*} \label{MNS_variance_2}
\Var\big[\Xi_{\psi,\tau} f\big]  = 
2\tau \int f(x)^2 dx
+ \iint |f(x)-f(y)|^2 \left( \frac{ \pi \tau}{\sinh\big[2\pi^2\tau(x-y)\big]} \right)^2 dxdy \ .
\end{equation*}

\section{Proofs of lemmas~\ref{thm:G3} and~\ref{thm:G4}} \label{A:G}

To prove lemmas~\ref{thm:G3} and~\ref{thm:G4}, the strategy is to exploit the symmetries of formula (\ref{GG}) in order to simplify as much as possible the cumulants of the random variable $\Xi_\Psi f$. To this end, we will use the following convention. Given two functions, we write $f \equiv g$ if there exists a permutation $\sigma \in \Sy(n)$ such that    $f(u)= g(\sigma u)$ or if $f(u)= g(-u)$ for all $u\in \R^n$.
For any vector $(u_i)_{i=1}^n$ of real numbers, we also denote 
$$\max^+\{u_1,\dots, u_n \}=\max\{0,u_1,\dots, u_n \} \ .$$
Unfortunately the combinatorial structure behind the cumulants of the modified ensembles seems to be rather complicated and consequently the following computations are rather technical. \\

{\it Proof of lemma~\ref{thm:G3}.}
According to definition (\ref{M}), we have $\M({\bf 2+1})= -\frac{3}{2}$ and $\M({\bf1}^3)=2$, and by formula (\ref{GG}),
\begin{align} \label{C3_2}
\mathfrak G ^3_{\Psi}[f]
=\int_{\R^3_0} d^2u\ \Re\left\{\prod_{i=1}^n \hat{f}(u_i) \right\} \int_{\R^3_<}{d^3x}\prod_{i=1}^n \Phi(x_i) \sum_{|{\bf m}|=n} \left\{
3 \G^{\bf 2+ 1}(u,x) + 3 \G^{\bf 1+ 2}(u,x)  - 4  \G^{ {\bf 1}^3}(u,x)\right\} \ .
\end{align}
Hence, to prove formula (\ref{C3}), we need to compute the kernels 
$ \G^{\bf 1+ 2}$ and  $\G^{ {\bf 1}^3}$.
By definition (\ref{Lambda}), for any $u\in \R^n_0$, 
 $$ \Lambda^{\bf 2+1}(u)=  \left(\begin{array}{cc} 0 &-u_3 \\ u_3 & 0 \end{array}  \right)
 \equiv \left(\begin{array}{cc} 0 & -u_1\\ u_1& 0 \end{array}  \right) =  - \Lambda^{\bf 1+2}(u) \ .
 $$
Then, by definition (\ref{G}), 
\begin{align} \notag
 \G^{\bf 2+1}(u,x)\equiv  \G^{\bf 1+2}(u,x) 
 &\notag\equiv
  \left[u_1-x_2+x_1\right]^++\left[u_1-x_3+x_1\right]^++ \left[-u_1-x_2+x_1\right]^+\\
&\notag\hspace{.3cm}+\left[-u_1-x_3+x_1\right]^++\left[u_1-x_3+x_2\right]^+ +\left[-u_1-x_3+x_2\right]^+ \\
&\equiv \label{G_21}
 2\big( \left[u_1-x_2+x_1\right]^++\left[u_1+u_2-x_3+x_1\right]^+
+ \left[u_2-x_3+x_2\right]^+ \big) \ .
\end{align}
One can check that each term of the first line corresponds to a permutation in $\Sy(3)$ in the following order $123, 132, 213, 312, 231, 321$. At the second step, we used the symmetries and  the condition $u_1+u_2+u_3=0$.
By a  similar argument, the contribution of the composition ${\bf 1}^3=(1,1,1)$ is given by
 \begin{equation*}
  \Lambda^{{\bf 1}^3}(u) 
 = \left(\begin{array}{ccc} 0 & -u_2 & -u_2-u_3 \\ u_2 & 0 & -u_3 \\  u_2+u_3 & u_3 & 0 \end{array}  \right)
 \end{equation*}
and 
\begin{align*} \G^{{\bf 1}^3}(u,x)&=\max^+\left\{u_2-x_2+x_1,u_2+u_3-x_3+x_1\right\}+\max^+\left\{u_2-x_3+x_1,u_2+u_3-x_2+x_1\right\} \\
&+\max^+\left\{-u_2-x_2+x_1,\ u_3-x_3+x_1\right\}+\max^+\left\{-u_2-x_3+x_1,\ u_3-x_2+x_1\right\}\\
&+\max^+\left\{-u_2-u_3-x_2+x_1,-u_3-x_3+x_1\right\}+\max^+\left\{-u_2-u_3-x_3+x_1,-u_3-x_2+x_1\right\} \ .
\end{align*}
It is straightforward to see that the previous formula implies that
\begin{equation} \label{G_3}
\G^{{\bf 1}^3}(u,x)\equiv 6\ \max^+\left\{u_1-x_2+x_1,u_1+u_2-x_3+x_1\right\} \ .
\end{equation}

If we combine equations (\ref{G_21}) and (\ref{G_3}), we get
\begin{equation} \label{G^3} 
3 \G^{\bf 2+ 1}(u,x) + 3 \G^{\bf 1+ 2}(u,x)  - 4  \G^{ {\bf 1}^3}(u,x)
\equiv 12\ \varpi(u_1-x_2+x_1, u_2-x_3+x_2) \ .
\end{equation}
Then, if we make the change of variables $x_1=x$, $z_1= x_2-x_1$ and $z_2= x_3-x_2$ in (\ref{C3_2}), formula (\ref{G^3}) implies that
\begin{equation*}
\Cu^3\big[\Xi_\Psi f \big]= 12\int_{\R^2}d^2u\ \Re\left\{ \prod_i \hat{f}(u_i) \right\} \int_\R dx\ \Phi(x) \iint\limits_{(0,\infty)^2} \Phi(x+z_1)\Phi(x+z_1+z_2) \varpi(u_1-z_1,u_2-z_2) d^2z \ ,
\end{equation*}
where it is understood that $u_3=-u_1-u_2$ in the first integral.
\qed\\

 {\it Proof of lemma~\ref{thm:G4}.}  We fix $z\in \R_+^3$. We will proceed exactly as in the proof of lemma~\ref{thm:G3} except that we will not  give all the details. We will denote $\pm=+1$ or $-1$ and we let 
 \begin{align*}
 \zeta_1 &= \max^+\{1-z_1,2-z_1-z_2\} + \max^+\{1-z_1,2-z_1-z_2-z_3\} +  \max^+\{1-z_1-z_2, 2-z_1-z_2-z_3\} \\
 &\hspace{.3cm}+ 2[2-z_1 ]^+ + [2-z_1-z_2 ]^+ \ ,  \\
  \zeta_2 &=2\bigg( [1-z_1]^+ + [1-z_1-z_2 ]^+ + [1-z_1-z_2 -z_3]^+ \bigg) \ ,  \\
  \zeta_3 &=  \max^+\{1-z_2,2-z_2-z_3\}+   [2-z_2]^+  +  [1-z_2]^+ + [1-z_2-z_3 ]^+ \ , \\
 \zeta_4&= 4[1-z_1]^+ +2 [1-z_1-z_2 ]^+ \ ,\\
 \zeta_5 &=4[1-z_2]^+ +2 [1-z_2-z_3 ]^+ \ .
 \end{align*}
 We will compute the values of $\tilde \G^{\bf m}(v,z)$ for all compositions ${\bf m}$ of $4$ and all points $v=(\pm,\pm,\pm,\pm)$ such that $\sum v_i =0$. The computations are not difficult but there are many cases to check.
By definition (\ref{Lambda}),
\begin{equation} \label{Lambda_22}
 \Lambda^{\bf 2+2}(u)=  \left(\begin{array}{cc} 0 & -u_3-u_4 \\ u_3+u_4 & 0 \end{array}  \right)
 \equiv \left(\begin{array}{cc} 0 & -u_1-u_2 \\ u_1+u_2 & 0 \end{array}  \right) \ .
 \end{equation}
 Then, by defintion (\ref{G}), we can check that
\begin{align*}&\tilde \G^{\bf 2+2}(++--)\equiv\tilde \G^{\bf 2+2}(--++) \\
&\equiv2\big( [2-z_1 ]^+ +[2-z_2 ]^+ + [2-z_3 ]^+ + [2-z_1-z_2 ]^+ 
+ [2-z_2-z_3 ]^+  +  [2-z_1-z_2 -z_3]^+ \big) \ .
\end{align*}
We used that $v_1+v_2=0$ or $2$ and that $z_i \ge 0$ to check which terms are a priori non-zero. Moreover, for the same reasons,
\begin{equation*}\tilde \G^{\bf 2+2}(+-+-)= \tilde \G^{\bf 2+2}(+--+)+\tilde \G^{\bf 2+2}(-++-)= \tilde \G^{\bf 2+2}(-+-+) =0 \ .
 \end{equation*}

If we use the symmetry of the function $\Theta$, see (\ref{Theta}), under the change of variable $z_1\leftrightarrow z_3$ and that $\M({\bf 2+2})= -3$, we can conclude that
\begin{equation}  \label{G_22}
\M({\bf 2+2})\sum_{\begin{subarray}{c} v_1+\cdots+v_4=0\\ v_i\in\{-1,1\} \end{subarray} } \tilde \G^{\bf 2+2}(v,z) 
\leftarrow-24\bigg( [2-z_1 ]^+ + \frac{[2-z_2 ]^+}{2} + [2-z_1-z_2 ]^+  + \frac{ [2-z_1-z_2 -z_3]^+}{2} \bigg)
\end{equation}
in the sense that if we replace the l.h.s$.$ of equation (\ref{G_22}) by its RHS in formula  (\ref{C4}), it does not change of the value of the integral.

Let us continue with the compositions ${\bf 3 +1}$ and  ${\bf 1+3}$. We have
 $$ \Lambda^{\bf 3+1}(u)=  \left(\begin{array}{cc} 0 &-u_4 \\ u_4 & 0 \end{array}  \right)
 \equiv \left(\begin{array}{cc} 0 & -u_1\\ u_1& 0 \end{array}  \right) =  - \Lambda^{\bf 1+3}(u) \ .
 $$
This expression depends on a single variable (say $u_1$) and  collecting the non-zero terms yields
 \begin{align*} &\tilde \G^{\bf 3+1}(\pm,\pm,\pm,\pm)\equiv\G^{\bf 1+3}(\pm,\pm,\pm,\pm) \\
&\equiv2\big( [1-z_1 ]^+ +[1-z_2 ]^+ + [1-z_3 ]^+ + [1-z_1-z_2 ]^+ 
+ [1-z_2-z_3 ]^+  
 +  [1-z_1-z_2 -z_3]^+ \big) \ .
\end{align*}
We can again use the symmetry of formula (\ref{C4}) and, since $\M({\bf 3+1})=-2$, we get
\begin{align}  \notag
&\M({\bf 3+1})\sum_{\begin{subarray}{c} v_1+\cdots+v_4=0\\ v_i\in\{-1,1\} \end{subarray} } \tilde \G^{\bf 3+1}(v,z) +  \tilde \G^{\bf 1+3}(v,z)  \\
&\notag\hspace{1cm}\leftarrow-48\bigg( 2[1-z_1 ]^+ + [1-z_2 ]^++ 2[1-z_1-z_2 ]^+  + [1-z_1-z_2 -z_3]^+ \bigg)\\
&\label{G_31}
\hspace{1cm}=-24\big( \zeta_2 +2 [1-z_1]^+  +2 [1-z_2 ]^+ +2[1-z_1-z_2 ]^+  \big) \ .
 \\ \notag
\end{align}

Consider now the composition ${\bf 1}^4={\bf 1+1+1+1}$. By definition (\ref{Lambda}),
  $$ \Lambda^{{\bf 1}^4} 
 = \left(\begin{array}{cccc} 0 & -u_2 & -u_2-u_3 &-u_2-u_3 -u_4 \\u_2 & 0 &-u_3 &-u_3-u_4  \\u_2+u_3 & u_3 & 0 & -u_4 \\ u_2+ u_3+u_4 & u_3+ u_4 & u_4 & 0  \end{array}  \right) \ .
$$
If we look at all permutations in $\Sy(4)$ and use the symmetry under permutations of the $u_i$'s, we get
 \begin{align*}\tilde{G}^{{\bf 1}^4}(u,z)\equiv\ 4\big(& \max^+\{ u_1-z_1,\ u_1+u_2-z_1-z_2,\ u_1+u_2+u_3-z_1-z_2-z_3\}\\
&+ \max^+\{ u_1-z_1-z_2-z_3,\ u_1+u_2-z_1-z_2,\ u_1+u_2+u_3-z_1\} \\
&+\max^+\{ u_1-z_1-z_2,\ u_1+u_2-z_1,\ u_1+u_2+u_3-z_1-z_2-z_3\} \\
&+\max^+\{ u_1-z_1-z_2-z_3,\ u_1+u_2-z_1,\ u_1+u_2+u_3-z_1-z_2\} \\
&+ \max^+\{ u_1-z_1-z_2,\ u_1+u_2-z_1-z_2-z_3,\ u_1+u_2+u_3-z_1\} \\
&+\max^+\{ u_1-z_1,\ u_1+u_2-z_1-z_2-z_3,\ u_1+u_2+u_3-z_1-z_2\} \big) \ .
\end{align*}
So we can assume that $\tilde{G}^{{\bf 1}^4}$ is given by the RHS of this expression, then  it is straightforward to check that
\begin{align*}
\tilde{G}^{{\bf 1}^4}(++--) &=  8\bigg( \max^+\{1-z_1,2-z_1-z_2\} + \max^+\{1-z_1,2-z_1-z_2-z_3\}  + [2-z_1 ]^+ \bigg) \ , \\
\tilde{G}^{{\bf 1}^4}(-++-)&=\tilde{G}^{{\bf 1}^4}(+--+)= 4\zeta_2 \ , \\
\tilde{G}^{{\bf 1}^4}(+-+-)&=4\zeta_4 \ , \\
\tilde{G}^{{\bf 1}^4}(--++)&=\tilde{G}^{{\bf 1}^4}(-+-+)=0 \ .
\end{align*}
So that, since $\M({\bf 1}^4)=-6$,
\begin{align} \label{G_4}
 \M({\bf 1}^4) \sum_{\begin{subarray}{c} v_1+\cdots+v_4 =0\\ v_i\in\{-1,1\} \end{subarray} } \tilde \G^{{\bf 1}^4}(v,z) = -24 \bigg\{& 2\zeta_2 +\zeta_4  + 2 \max^+\{1-z_1,2-z_1-z_2\}\\
 &\notag + 2\max^+\{1-z_1,2-z_1-z_2-z_3\}  + 2[2-z_1 ]^+  \bigg\}  \ .
 \end{align}
 
 If we combine formulae (\ref{G_22}), (\ref{G_31}) and (\ref{G_4}),
  \begin{align} \label{G_1111}
  &\sum_{\begin{subarray}{c} v_1+\cdots+v_4=0\\ v_i\in\{-1,1\} \end{subarray} }  \sum_{\begin{subarray}{c}|{\bf m}|=4\\ \ell({\bf m})\neq 3 \end{subarray}} \M({\bf m})\tilde \G^{\bf m}(v,z) \\
&\hspace{1cm}\notag\leftarrow -24 \bigg\{ 3\zeta_2+\zeta_4 +2 \max^+\{1-z_1,2-z_1-z_2\}
 +2\max^+\{1-z_1, 2-z_1-z_2-z_3\}
+2 [1-z_1]^+ \\
    &\notag \hspace{1cm}+2 [1-z_2 ]^+ +2[1-z_1-z_2 ]^+ 
 3[2-z_1 ]^+  + \frac{[2-z_2 ]^+}{2} + [2-z_1-z_2 ]^+  + \frac{ [2-z_1-z_2 -z_3]^+}{2}\bigg\} \ .
\end{align}

Finally, we look at the composition ${\bf 2+1+1}$,
 \begin{equation*} \label{Lambda_211}
  \Lambda^{\bf 2+1+1}(u) 
 = \left(\begin{array}{ccc} 0 & -u_3 & -u_3-u_4 \\ u_3 & 0 & -u_4 \\  u_3+u_4 & u_4 & 0 \end{array}  \right)
 \equiv\left(\begin{array}{ccc} 0 & -u_1 & -u_1-u_2 \\ u_1 & 0 & -u_2 \\  u_1+u_2 & u_2 & 0 \end{array}  \right) \ ,
 \end{equation*}
and if we follow the same procedure, we can prove that
 \begin{align*}\tilde \G^{\bf 2+1+1}(++--) = \G^{\bf 2+1+1}(--++)& 
 = \zeta_1 + \zeta_2+\zeta_3 \\
\tilde \G^{\bf 2+1+1}(+-+-)+\tilde \G^{\bf 2+1+1}(-+-+)&= \tilde \G^{\bf 2+1+1}(+--+)+ \tilde \G^{\bf 2+1+1}(-++-) \\
&= 2\zeta_2+\zeta_4+\zeta_5
\end{align*}

  On the other hand
   $$ \Lambda^{1+2+1} 
 = \left(\begin{array}{ccc} 0 & -u_2-u_3 & -u_2-u_3-u_4\\ u_2+u_3 & 0 & -u_4 \\  u_2+u_3+u_4 & u_4 & 0 \end{array}  \right)
 \equiv\left(\begin{array}{ccc} 0 & u_1 +u_2 & u_2 \\ -u_1-u_2 & 0 & -u_1 \\  -u_2 & u_1 & 0 \end{array}  \right) \ .
 $$

 It is not difficult to see that, up to conjugation by a permutation matrix, we have $ \Lambda^{1+2+1}\equiv \Lambda^{2+1+1}$. This implies that $ \tilde \G^{1+2+1}\equiv\tilde \G^{2+1+1}$ because such conjugation only changes the order of the sum over $\Sy(4)$ in the definition (\ref{G}). Similarly, we can check that the matrix  $ \Lambda^{1+1+2}$ is also conjugated to $\Lambda^{2+1+1}$ by a permutation matrix, so that they give the same contribution to the 4$^{\text{th}}$ cumulant.
Since $\M({\bf 2+1+1})= 4$, putting all terms together, we conclude that
 \begin{equation*}\sum_{\begin{subarray}{c} v_1+\cdots+v_4=0\\ v_i\in\{-1,1\} \end{subarray} }  \sum_{\begin{subarray}{c}|{\bf m}|=4\\ \ell({\bf m})=3 \end{subarray}} \M({\bf m})\tilde \G^{\bf m}(v,z) 
= 24\left\{ \zeta_1 + 3\zeta_2+\zeta_3+\zeta_4+\zeta_5 \right\} \ .
\end{equation*}
 
 Observe that using the symmetry between $z_1$ and $z_3$ and the DHK formula (\ref{DHK_2}), we can show that
  $$ \zeta_1 + \zeta_3 \leftarrow 2 \bigg( [1-z_1]^++[1-z_2]^++  [1-z_1-z_2]^++ [2-z_1 ]^++ \frac{ [2-z_2]^+}{2} +[2-z_1-z_2]^+ +\frac{[2-z_1-z_2-z_3]^+}{2}  \bigg) $$
and we get
  \begin{align} \notag
&  \sum_{\begin{subarray}{c} v_1+\cdots+v_4=0\\ v_i\in\{-1,1\} \end{subarray} }  \sum_{\begin{subarray}{c}|{\bf m}|=4\\ \ell({\bf m})= 3 \end{subarray}} \M({\bf m})\tilde \G^{\bf m}(v,z) 
\leftarrow24 \bigg\{ 3\zeta_2+\zeta_4 +\zeta_5 +2 [1-z_1]^+  +2 [1-z_2 ]^+ +2[1-z_1-z_2 ]^+  \\
&\hspace{1cm}\label{G_211} 
+ 2[2-z_1 ]^+  + [2-z_2 ]^+ +2 [2-z_1-z_2 ]^+  + [2-z_1-z_2 -z_3]^+\bigg\} \ .
\end{align}

Finally, if we combine formulae (\ref{G_1111}) and (\ref{G_211}), many terms cancel but not all of them and we are left with
  \begin{align*}
  &\sum_{\begin{subarray}{c} v_1+\cdots+v_4=0\\ v_i\in\{-1,1\} \end{subarray} } 
  \sum_{|{\bf m}|=4} \M({\bf m})\tilde \G^{\bf m}(v,z)
\leftarrow 24\bigg( \zeta_5 - [2-z_1 ]^+  + \frac{[2-z_2 ]^+}{2} + [2-z_1-z_2 ]^+\\
&\hspace{1cm}  + \frac{ [2-z_1-z_2 -z_3]^+}{2}  -2 \max^+\{1-z_1,2-z_1-z_2\}  -2\max^+\{1-z_1, 2-z_1-z_2-z_3\}\bigg) \ .
\end{align*}
Finally, if we make the change of variable $z_1\leftrightarrow z_3$,
 $$  \zeta_5 \leftarrow 4[1-z_2]^+ +2 [1-z_1-z_2 ]^+  \ ,$$ 
and we have proved the formula of lemma~\ref{thm:G4}.\qed

\vspace{1cm}

{\large\bf Acknowledgement:} We thank Maurice Duits for helpful discussions regarding his related works on central limit theorems for orthogonal polynomial ensembles  and Erik Duse for suggesting the proof of proposition \ref{thm:MNS_property}. We thank the referees for helpful comments and suggestions on the exposition.


\begin{thebibliography}{10}

\bibitem{AHM_11}
{\sc Y.~Ameur, H.~Hedenmalm, and N.~Makarov}, {\em Fluctuations of eigenvalues
  of random normal matrices}, Duke Math. J. 159,  (2011), pp.~31--81.

\bibitem{ACQ11}
{\sc G.~Amir, I.~Corwin, J.~Quastel}, {\em Probability distribution of the free energy of the continuum directed random polymer in 1+1 dimensions}, Comm. Pure Appl. Math. 64 (2011), no. 4, pp.~466--537.


\bibitem{AGZ}
{\sc G.~W. Anderson, A.~Guionnet, and O.~Zeitouni}, {\em An Introduction to
  Random Matrices}, Cambridge Univ. Press, Cambridge, (2010).

\bibitem{Borodin_11}
{\sc A.~Borodin}, {\em Determinantal point processes}, in The {O}xford handbook
  of random matrix theory, ch.~11,
  Oxford Univ. Press, Oxford, (2011).

\bibitem{BEYY_14}
{\sc P.~Bourgade, L.~Erd\H{o}s, H.-T. Yau, and J.~Yin}, {\em Fixed energy
  universality for generalized {W}igner matrices}.
\newblock arXiv:1407.5606.

\bibitem{BK_99a}
{\sc A.~{Boutet de Monvel} and A.~Khorunzhy}, {\em Asymptotic distribution of
  smoothed eigenvalue density~{I}. {G}aussian random matrices}. Random Oper.
  Stochastic Equations 7,  (1999), pp.~1--22.

\bibitem{BK_99b}
\leavevmode\vrule height 2pt depth -1.6pt width 30pt, {\em Asymptotic
  distribution of smoothed eigenvalue density~{II}. {W}igner random matrices}.
  Random Oper. Stochastic Equations 7, no. 2 (1999), pp.~149--168.
  
  
  \bibitem{BD_14}
  {\sc J.~Breuer and M.~Duits}, 
 {\em Universality of mesoscopic fluctuations for orthogonal polynomial {E}nsembles}.
Comm. Math. Phy. 342,  no.~2  (2016), pp.~491--531.
  
  \bibitem{BD_15}
  \leavevmode\vrule height 2pt depth -1.6pt width 30pt, {\em {C}entral {L}imit {T}heorems for
  biorthogonal {E}nsembles and asymptotics of recurrence coefficients}, \newblock arXiv:1309.6224.


\bibitem{CL_95}
{\sc O.~Costin and J.~Lebowitz}, {\em Gaussian fluctuations in random
  matrices}, Phys. Rev. Lett. 75,  (1995), pp.~69--72.

\bibitem{DDMS_14}
{\sc D.~S. Dean, P.~L. Doussal, S.~N. Majumdar, and G.~Schehr}, {\em Finite
  temperature free fermions and the {K}ardar-{P}arisi-{Z}hang equation at
  finite time}, Phys. Rev. Lett. 114,  (2015).


\bibitem{Deift_al_99_a}
{\sc P.~Deift, T.~Kriecherbauer, K.~T.-R. McLaughlin, S.~Venakides, and
  X.~Zhou}, {\em Strong asymptotics of orthogonal polynomials with respect to
  exponential weights}, Comm. Pure Appl. Math. 52 (12),  (1999),
  pp.~1491--1552.

\bibitem{Deift_al_99_b}
\leavevmode\vrule height 2pt depth -1.6pt width 30pt, {\em Uniform asymptotics
  for polynomials orthogonal with respect to varying exponential weights and
  applications to universality questions in random matrix theory}, Comm. Pure
  Appl. Math. 52 (11),  (1999), pp.~1335--1425.


\bibitem{DJ_14}
{\sc M.~Duits and K.~Johansson}, {\em On mesoscopic equilibrium for linear
  statistics in {D}yson's {B}rownian motion}.
\newblock arXiv:1312.4295.

\bibitem{Durrett_10}
{\sc R.~Durrett}, {\em Probability: theory and examples}, Cambridge Series in
  Statistical and Probabilistic Mathematics, Cambridge Univ. Press,
  Cambridge, 4th~ed., (2010).

\bibitem{Erdos_14}
{\sc L.~Erd\H{o}s}, {\em Random matrices, log-gases and {H}\"older regularity},
  in Proceedings of ICM 2014.

\bibitem{EK_14a}
{\sc L.~Erd\H{o}s and A.~Knowles}, {\em The {A}ltshuler-{S}hklovskii formulas
  for random band matrices {I}: the unimodular case}, Comm. Math. Phy. 333,
  (2015), pp.~1365--1416.

\bibitem{EK_14b}
\leavevmode\vrule height 2pt depth -1.6pt width 30pt, {\em The
  {A}ltshuler-{S}hklovskii formulas for random band matrices {II}: the general
  case}, Ann. Henri Poincar\'{e} 16,  (2015), pp.~709--799.

\bibitem{EY_12}
{\sc L.~Erd\H{o}s and H.-T. Yau}, {\em Universality of local spectral
  statistics of random matrices}, Bull. Amer. Math. Soc. vol. 49,  (2012),
  pp.~377--414.

\bibitem{FKS_13}
{\sc Y.~V. Fyodorov, B.~A. Khoruzhenko, and N.~J. Simm}, {\em Fractional
  {B}rownian motion with {H}urst index {$H=0$} and the {G}aussian {U}nitary
  {E}nsemble},  Ann. Probab. 44 (2016), no. 4, pp.~2980--3031.

\bibitem{Gutzwiller}
{\sc M.~C. Gutzwiller}, {\em Chaos in Classical and Quantum Mechanics}, Interdisciplinary Applied Mathematics, 1. 
  Springer-Verlag, New-York (1990).

\bibitem{HKPV_06}
{\sc B.~Hough, M.~Krishnapur, Y.~Peres, and B.~Vir{\'a}g}, {\em Determinantal
  {P}rocesses and {I}ndependence}, Probab. Surv. 3,  (2006), pp.~206--229.


\bibitem{Johansson_05}
{\sc K.~Johansson}, {\em Random matrices and determinantal processes}, in
  Mathematical statistical physics, Elsevier B. V., Amsterdam, (2006), pp.~1--55.


\bibitem{Johansson_07}
\leavevmode\vrule height 2pt depth -1.6pt width 30pt, {\em From Gumbel to {T}racy-{W}idom}, Probab. Theory
  Relat. Fields 138,  (2007), pp.~75--112.



\bibitem{KSSV_14}
{\sc T.~Kriecherbauer, K.~Schubert, K.~Sch{\"u}ler, and M.~Venker}, {\em Global
  asymptotics for the {C}hristoffel-{D}arboux kernel of random matrix theory}.
\newblock arXiv:1401.6772.


\bibitem{L_15a}
{\sc G.~Lambert}, {\em Mesoscopic
  fluctuations for unitary invariant ensembles}.
\newblock arXiv:1510.03641

\bibitem{L_15b}
\leavevmode\vrule height 2pt depth -1.6pt width 30pt, {\em CLT for biorthogonal ensembles and related combinatorial identities}.
\newblock arXiv:1511.06121 

\bibitem{Lieb_Loss}
{\sc E.~H.~Lieb and M.~Loss}, {\em Analysis}, vol. 14 of Graduate Studies in Mathematics, Amer.
Math. Soc., Providence, RI, 2nd ed. (2001)




\bibitem{MNS_94}
{\sc M.~Moshe, H.~Neuberger, and B.~Shapiro}, {\em Generalized ensemble of
  random matrices}, Phys. Rev. Lett. 73,  (1994), pp.~1497--1500.

\bibitem{Pastur_Shcherbina}
{\sc L.~A. Pastur and M.~Shcherbina}, {\em Eigenvalue {D}istribution of {L}arge
  {R}andom {M}atrices}, Mathematical Surveys and Monographs 171, Amer. Math. Soc., Providence, RI, (2011).

\bibitem{PR_29}
{\sc M.~Plancherel and W.~Rotach}, {\em Sur les valeurs asymptotiques des
  polynomes d'{H}ermite}, Commentarii Mathematicii Helvetici
  1,  (1929), pp.~227--254.

\bibitem{RV_07b}
{\sc B.~Rider and B.~Vir{\'a}g}, {\em Complex determinantal processes and
  {$H^1$} noise}, Elect. J. Probab.~12,  (2007), pp.~1238--1257.

\bibitem{RV_07a}
\leavevmode\vrule height 2pt depth -1.6pt width 30pt, {\em The noise in the
  circular law and the {G}aussian free field}, Int. Math. Res. Not.,  (2007).

\bibitem{SS10a}
{\sc T.~Sasamoto, H.~Spohn}, {\em The crossover regime for the weakly asymmetric simple exclusion process}, J. Stat. Phys.~140 (2010), no. 2, pp.~209--231.

\bibitem{SS10b}
\leavevmode\vrule height 2pt depth -1.6pt width 30pt, {\em Exact height distributions for the KPZ equation with narrow wedge initial condition}, Nuclear Phys. B 834 (2010), no. 3, 523--542. 


\bibitem{Simon_04a}
{\sc B.~Simon}, {\em Orthogonal {P}olynomials on the {U}nit {C}ircle, Part 1:
  {C}lassical {T}heory}, Colloquium Publications vol. 54, Amer. Math. Soc., Providence, RI,  (2004).

\bibitem{Simon_05}
\leavevmode\vrule height 2pt depth -1.6pt width 30pt, {\em Trace {I}deals and
  {T}heir {A}pplications}, Mathematical Surveys and Monographs 120, Amer. Math. Soc., Providence, RI, 2nd~ed., (2005).

\bibitem{Soshnikov_00a}
{\sc A.~Soshnikov}, {\em The {C}entral
  {L}imit {T}heorem for local linear statistics in classical compact groups and
  related combinatorial identities}, Ann. Probab. 28,  (2000), pp.~1353--1370.

\bibitem{Soshnikov_00b}
\leavevmode\vrule height 2pt depth -1.6pt width 30pt, {\em Determinantal random point fields}, Russian Math.
  Surv. 55,  (2000), pp.~923--975.


\bibitem{Soshnikov_01}
\leavevmode\vrule height 2pt depth -1.6pt width 30pt, {\em Gaussian limit for
  determinantal random point fields}, Ann. Probab. 30,  (2001), pp.~1--17.

\end{thebibliography}
 \end{document}